\documentclass[a4paper,12pt,leqno]{article}

\usepackage{amsmath,amssymb,amsthm}
\newcommand{\alphabar}{{\overline{\alpha}}}
\newcommand{\gammabar}{{\overline{\gamma}}}

\newcommand{\afrak}{{\mathfrak{a}}}
\newcommand{\bfrak}{{\mathfrak{b}}}
\newcommand{\cfrak}{{\mathfrak{c}}}
\newcommand{\efrak}{{\mathfrak{e}}}
\newcommand{\hfrak}{{\mathfrak{h}}}
\newcommand{\kfrak}{{\mathfrak{k}}}
\newcommand{\mfrak}{{\mathfrak{m}}}
\newcommand{\pfrak}{{\mathfrak{p}}}
\newcommand{\qfrak}{{\mathfrak{q}}}
\newcommand{\rfrak}{{\mathfrak{r}}}
\newcommand{\sfrak}{{\mathfrak{s}}}
\newcommand{\tfrak}{{\mathfrak{t}}}
\newcommand{\wfrak}{{\mathfrak{w}}}

\newcommand{\afrakbar}{{\overline{\afrak}}}
\newcommand{\bfrakbar}{{\overline{\bfrak}}}
\newcommand{\cfrakbar}{{\overline{\cfrak}}}
\newcommand{\pfrakbar}{{\overline{\pfrak}}}
\newcommand{\sfrakbar}{{\overline{\sfrak}}}

\newcommand{\Afrak}{{\mathfrak{A}}}
\newcommand{\Dfrak}{{\mathfrak{D}}}
\newcommand{\Ffrak}{{\mathfrak{F}}}
\newcommand{\Gfrak}{{\mathfrak{G}}}
\newcommand{\Hfrak}{{\mathfrak{H}}}
\newcommand{\Jfrak}{{\mathfrak{J}}}
\newcommand{\Mfrak}{{\mathfrak{M}}}
\newcommand{\Nfrak}{{\mathfrak{N}}}
\newcommand{\Pfrak}{{\mathfrak{P}}}
\newcommand{\Rfrak}{{\mathfrak{R}}}
\newcommand{\Zfrak}{{\mathfrak{Z}}}

\newcommand{\Pssf}{{\textsf{\textit{P}}}}
\newcommand{\Rssf}{{\textsf{\textit{R}}}}

\newcommand{\Nm}{{N}}
\newcommand{\Tr}{{Sp}}
\renewcommand{\Im}{I}

\def\clap#1{\hbox to 0pt{\hss#1\hss}}
\def\mathclap{\mathpalette\mathclapinternal}
\def\mathclapinternal#1#2{%
  \clap{$\mathsurround=0pt#1{#2}$}}

\begin{document}

\title{{\bf
Class Fields of
Complex Multiplication
}}

\author{by\\
\scshape
Max Deuring
\\
\small\scshape
in G\"ottingen}

\date{\footnotesize
translated from German by
\scshape
CheeWhye Chin
}

\maketitle

\renewcommand{\contentsname}{\hfil Contents}
\let\stdcontentsline\contentsline
\renewcommand\contentsline[3]{\stdcontentsline{#1}{#2}{}}
\makeatletter
\renewcommand{\@dotsep}{10000}
\makeatother
\setcounter{tocdepth}{5}
\tableofcontents

\bigskip
\bigskip
\bigskip

\setlength{\parskip}{3ex plus 0.2ex}


\section*{A. Function-theoretic basics
\normalsize{\footnote{
For the tools from
function theory
used below,
please refer to
the literature.
}}
}
\addcontentsline{toc}{section}{\hfill
A. Function-theoretic basics
\hfill}


\paragraph{1. Modular functions.}

Let $\Mfrak^{(h)}$ denote
the \emph{homogeneous modular group},
i.e.~the group of all
integer-entry matrices
$M
 =
 \begin{pmatrix}
 a & b
 \\
 c & d
 \end{pmatrix}
$
with
$|M| = ad{\,-\,}bc = 1$.
Each $M \in \Mfrak$
gives rise to
a homogeneous linear transformation
$\begin{pmatrix}
 \omega_1
 \\
 \omega_2
 \end{pmatrix}
 \to
 M
 \,
 \begin{pmatrix}
 \omega_1
 \\
 \omega_2
 \end{pmatrix}
 =
 \begin{pmatrix}
 \omega_1'
 \\
 \omega_2'
 \end{pmatrix}
$
of a pair
$\omega_1,\omega_2$
of complex variables;
the group of these transformations
is isomorphic to
$\Mfrak^{(h)}$
and
will henceforth
also be denoted by
$\Mfrak^{(h)}$.
The fractional linear transformations
$\omega
 =
 \dfrac{\omega_1}{\omega_2}
 \to
 \omega'
 =
 \dfrac{\omega_1'}{\omega_2'}
 =
 \dfrac{a\,\omega{+}b}{c\,\omega{+}d}
 =
 M(\omega)
$
associated to
$M
 =
 \begin{pmatrix}
 a & b
 \\
 c & d
 \end{pmatrix}
$
from $\Mfrak^{(h)}$
form
the \emph{inhomogeneous modular group}
$\Mfrak$;
the kernel of the homomorphism
\[
  M
  \ \to\ 
  \left[\ 
    \omega
    \to
    M(\omega)
  \ \right]
\]
is
$\{E,-E\}$,
and
$\Mfrak
 \cong
 \Mfrak^{(h)}/\{E,-E\}
$.
We reserve 
the notation $\omega$
for the quotient
$\dfrac{\omega_1}{\omega_2}$;
here,
$\omega$
will always be
restricted to
the half-plane
$\Im\,\omega > 0$.

A function
$f(\omega)$
of the complex variable $\omega$
is called
a \emph{modular form of weight~$t$}
if:
\begin{itemize}
\itemsep0pt
\item[1.]
$f(\omega)$ is holomorphic
in the half-plane
$\Im\,\omega > 0$
except for poles,
\item[2.]
$f(M(\omega))
 =
 f(\omega)
 \,
 (c\,\omega{+}d)^{-t}
$
holds for every
$M
 =
 \begin{pmatrix}
 a & b
 \\
 c & d
 \end{pmatrix}
 \in 
 \Mfrak
$,
\item[3.]
the Fourier expansion
${\displaystyle
 f(\omega)
 =
 \sum_{\nu}
 a_\nu\,q^\nu
}$,
$q = e^{2\pi i \omega}$,
in a suitable half-plane
$\Im\,\omega > \alpha$,
which is possible
since
$f(\omega{+}1) = f(\omega)$,
has at most
finitely many terms
with negative exponents,
so that
$f(\omega)$
as a function of~$q$
is regular or polar
at $q=0$.
\end{itemize}

A modular form
is called
\emph{entire}
if it is holomorphic in
$\Im\,\omega > 0$.
Modular forms of level~0,
being functions
invariant under $\Mfrak$,
are called
\emph{modular functions};
they form a field
$k_{\Mfrak}$.
Here,
$k$ always denotes
the field of complex numbers.

The existence of
modular functions and modular forms
can be deduced
by means of
the existence theorem
in function theory,
but for the following
we will rely on
the explicit formulae
from the classical theory of
elliptic functions
and modular functions,
as these also allow one
to read off
the arithmetic properties of
the fundamental series expansion.

For the lattice
$\wfrak$
of all
$w = n_1\,\omega_1 + n_2\,\omega_2$,
where it is understood that
$n_\nu = 0,\pm1,\pm2,\ldots$,
and for any integer
$m \geqslant 2$,
one has
an entire modular form
of weight~$-2m$
given by
\[
  G_m(\omega)
  \ =\ 
  \omega_2^{2m}
  \,
  \sum_{w \in \wfrak,\ w \neq 0}
  w^{-2m}
  .
  \tag{1}
\]
Its Fourier expansion
reads as
\[
  G_m(\omega)
  \ =\ 
  \frac{(2\pi)^{2m}}{(2m)!}
  \ 
  \left[\ 
  B_m
  \ +\ 
  4m
  (-1)^m
  \ 
  \sum_{n=1}^{\infty}
  n^{2m-1}
  \ 
  \sum_{\nu=1}^{\infty}
  q^{n\,\nu}
  \ \right]
  ,
  \tag{2}
  \label{eqn:2}
\]
where
$B_m
 =
 \text{$m$-th Bernoulli number}
$.
Set
\[
  g_2(\omega)
  \ =\ 
  60
  \,
  G_2(\omega)
  ,
  \qquad
  g_3(\omega)
  \ =\ 
  140
  \,
  G_3(\omega)
  ;
  \tag{3}
\]
the \emph{discriminant}
\[
  \Delta(\omega)
  \ =\ 
  g_2(\omega)^3
  -
  27\,g_3(\omega)^2
  \tag{4}
\]
is an entire modular form
of weight~$-12$.
The \emph{Dedekind function}
\[
  \eta(\omega)
  \ =\ 
  q^{\frac{1}{24}}
  \ 
  \prod_{n=1}^{\infty}
  (1-q^n)
  \tag{5}
  \label{eqn:5}
\]
is holomorphic
in $\Im\,\omega > 0$,
and $\Delta(\omega)$
can also be expressed
in terms of it
as
\[
  \Delta(\omega)
  \ =\ 
  (2\pi)^{12}
  \ 
  \eta(\omega)^{24}
  .
  \tag{6}
  \label{eqn:6}
\]
The \emph{absolute invariant of the modular group}
\[
  j(\omega)
  \ =\ 
  2^6
  \,
  3^3
  \,
  g_2(\omega)^3
  \,
  \Delta(\omega)^{-1}
  \tag{7}
\]
is an entire modular function.
One has
\[
  j(i)
  \ =\ 
  2^6\,3^3
  ,
  \qquad
  j(e^{\textstyle\frac{2\pi i}{3}})
  \ =\ 
  0
  ,
  \qquad
  j(\infty)
  \ =\ 
  \infty
  .
  \tag{8}
  \label{eqn:8}
\]
Here,
$j(\omega)$
can be considered as
a meromorphic function
on the quotient space of
the half-plane
$\Im\,\omega > 0$
(including $\infty$)
mod~$\Mfrak$,
which is
a compact Riemann surface
$\Ffrak$
of genus~0;
it maps
this quotient space
bijectively onto
the Riemann sphere,
and is normalized
by~\eqref{eqn:8}.
From this
it follows that
$k_\Mfrak$
is the field of
rational functions
in $j(\omega)$
with complex number coefficients:
\ 
$k_\Mfrak = k(j(\omega))$.
The entire modular functions
are the polynomials
in $j(\omega)$.

Two lattices
$\wfrak$ and $\wfrak'$
are called \emph{equivalent}
if they differ only by
a scaling factor:
\ 
$\wfrak' = \varrho\,\wfrak$,
$\varrho \neq 0$.
We denote
an equivalence class of lattices
by $\kfrak$.
The classes $\kfrak$
correspond
in a uniquely reversible way
to the points~$\neq \infty$
of $\Ffrak$,
where
the class $\kfrak$
of a lattice
with basis $\omega_1,\omega_2$
is assigned to
the residue class of
$\omega = \dfrac{\omega_1}{\omega_2}$
modulo~$\Mfrak$.
Henceforth
we write
$j(\kfrak)
 =
 j
 \left(
 \dfrac{\omega_1}{\omega_2}
 \right)
$;
the value
$j(\kfrak)$
is called
the \emph{invariant of $\kfrak$},
and it determines
the class $\kfrak$
uniquely.

The Fourier expansions
\[
  j(\omega)
  \ =\ 
  q^{-1}
  \ +\ 
  c_0
  \ +\ 
  c_1\,q
  \ +\ 
  \cdots
  \tag{9}
  \label{eqn:9}
\]
and
\[
  \Delta(\omega)
  \ =\ 
  (2\pi)^{12}
  \,
  q
  \ 
  [\, 
  1
  \ +\ 
  D_1\,q
  \ +\ 
  \cdots
  \,]
  \tag{10}
\]
can be easily calculated
from~\eqref{eqn:2}
(that of $\Delta$
of course
also
from~\eqref{eqn:5}
and~\eqref{eqn:6}),
and from that follows
the fact that
the coefficients $c_\nu$ and $D_\nu$
are rational integers,
which is fundamental
for the number-theoretic properties
of the modular functions.
First
we use it
to get the following
conclusions:

For any number field $\Lambda$,
we denote by $\Lambda_\Mfrak$
the field
$\Lambda(j(\omega))$
of rational functions
in $j(\omega)$
with coefficients in $\Lambda$.
Then:
\ 
\textit{
$\Lambda_\Mfrak$
consists of
all modular functions
$f(\omega)$
whose Fourier expansion
$f(\omega)
 =
 \sum
 a_\nu\,q^\nu
$
has coefficients
$a_\nu$
in the field $\Lambda$.
}
Furthermore:
\ 
\textit{
if $\Lambda$ is
an algebraic number field,
then
an entire modular function
$f(\omega)
 =
 A_0
 +
 A_1\,j(\omega)
 +
 \cdots
 +
 A_N\,j(\omega)^N
$
has integers $A_\nu$ of $\Lambda$
as coefficients
if and only if
the Fourier coefficients
$a_\nu$ of $f(\omega)$
are integers of $\Lambda$.
}
Both follow easily
when
the Fourier expansion of $j(\omega)$
is inserted
into the expression of $f(\omega)$
as a rational function
(or polynomial)
of $j(\omega)$;
the first assertion
depends only on
the fact that
the $c_\nu$ in~\eqref{eqn:9}
are rational,
while
the second
depends on
the fact that
the $c_\nu$ are
rational integers
and that
the leading coefficient in~\eqref{eqn:9}
has value~1.


\paragraph{2. Modular functions for subgroups of $\Mfrak$.}

Let $\Nfrak$ be
a subgroup of $\Mfrak$
of finite index
$[\Mfrak{:}\Nfrak]$.
A \emph{modular form $f(\omega)$
of weight~$t$ for the group $\Nfrak$}
is defined by
the following
three requirements:
\begin{itemize}
\itemsep0pt
\item[1.]
$f(\omega)$ is regular
in
$\Im\,\omega > 0$
except for poles.
\item[2.]
For every
$N
 =
 \begin{pmatrix}
 a & b
 \\
 c & d
 \end{pmatrix}
 \in 
 \Nfrak
$,
one has
$f(N(\omega))
 =
 f(\omega)
 \,
 (c\,\omega{+}d)^{-t}
$.
\item[3.]
Let $M$ be
any element of $\Mfrak$.
If $\ell$ denotes
the smallest positive whole number
for which
$\begin{pmatrix}
 1 & 1
 \\
 0 & 1
 \end{pmatrix}^{\textstyle \ell}
$
lies in
$M^{-1}\,\Nfrak\,M$,
so that
$f(M(\omega))$
in condition~\textup{2}
has period~$\ell$,
then
in a suitable region
$\Im\,\omega > \varrho$,
$f(M(\omega))$
should have
a Fourier expansion
---
\emph{$q$-expansion}
---
of the form
\[
  f(M(\omega))
  \ =\ 
  \sum_{\nu=n}^{\infty}
  a_\nu
  \,
  q^{\textstyle \frac{\nu}{\ell}}
  ,
\]
in which
at most finitely many
coefficients $a_\nu$
of negative exponent~$\nu$
are different from~0.
Here
(and also always
in the following),
we set
\[
  q^{\textstyle \frac{1}{\ell}}
  \ =\ 
  e^{\textstyle \frac{2\pi i \omega}{\ell}}
  .
  \tag{11}
\]
If
$\Nfrak\,M_\nu$,
$\nu = 1,2,\ldots,[\Mfrak{:}\Nfrak]$,
are the residue classes of
$\Mfrak$ mod $\Nfrak$,
it suffices that
this requirement holds
for $M = M_1,\ldots,M_{[\Mfrak{:}\Nfrak]}$.
\end{itemize}

Forms of weight~0
are called \emph{functions}
for the group $\Nfrak$;
a form (function)
for $\Nfrak$
is called \emph{entire}
if it is regular
in $\Im\,\omega > 0$.

The functions
for the group $\Nfrak$
form a field extension $k_\Nfrak$
of $k_\Mfrak$.
For $M \in \Mfrak$,
an isomorphism $\lambda_M$
from $k_\Nfrak$
to $k_{M^{-1}\,\Nfrak\,M}$
is defined by
\[
  f(\omega)
  \,
  \lambda_M
  \ =\ 
  f(M(\omega))
  ,
  \tag{12}
\]
which leaves $k_\Mfrak$
elementwise fixed.
The isomorphism $\lambda_M$
depends only on
the residue class
$\Nfrak\,M$ of $\Nfrak$ in $\Mfrak$
to which $M$ belongs.
The coefficients of the polynomial
\[
  F(X)
  \ =\ 
  \prod_{\nu=1}^{[\Mfrak{:}\Nfrak]}
  (X - f(\omega)\lambda_{M_\nu})
  \ =\ 
  \prod_{\nu=1}^{[\Mfrak{:}\Nfrak]}
  (X - f(M_\nu(\omega)))
  \tag{13}
\]
remain unchanged
under each
modular substitution $M$,
and they also satisfy
conditions~\textup{1} and~\textup{3},
so
they lie in $k_\Mfrak$.
This proves that
$k_\Nfrak$ is algebraic
over $k_\Mfrak$
of degree~$\leqslant [\Mfrak{:}\Nfrak]$.
By the existence theorem
in function theory,
it follows that
\[
  [k_\Nfrak{:}k_\Mfrak]
  \ =\ 
  [\Mfrak{:}\Nfrak]
  ,
  \tag{14}
\]
so the functions
$f(\omega)\,\lambda_\nu
 =
 f(M_\nu(\omega))
$
are a full system
of conjugates of
$f(\omega)$ over $k_\Mfrak$,
and $F(X)$ is
the principal polynomial of
$f(\omega)$ over $k_\Mfrak$.

Later we will only need
the transformation subgroups
$\Nfrak = \Mfrak_S$
of $\Mfrak$.
In this case,
we will explicitly specify
functions of $k_\Nfrak$
of degree~$[\Mfrak{:}\Nfrak]$
over $k_\Mfrak$,
so that
here as well
we do not need to
rely on
the function-theoretic
existence theorem.

Instead of
a modular form $f(\omega)$
of weight~$t$,
it is sometimes convenient
to compute with
the associated
\emph{homogeneous modular form}
$f
 \begin{pmatrix}
 \omega_1
 \\
 \omega_2
 \end{pmatrix}
$,
which is defined by
\[
  f
  \begin{pmatrix}
  \omega_1
  \\
  \omega_2
  \end{pmatrix}
  \ =\ 
  f(\omega)
  \,
  \omega_2^t
  .
\]
Here,
$f
 \begin{pmatrix}
 \omega_1
 \\
 \omega_2
 \end{pmatrix}
$
is homogeneous
in $\omega_1,\omega_2$
of degree~$t$,
and the invariance property~\textup{2}
of $f(\omega)$
means for 
$f
 \begin{pmatrix}
 \omega_1
 \\
 \omega_2
 \end{pmatrix}
$
simply that
\[
  f
  \left(
  N
  \begin{pmatrix}
  \omega_1
  \\
  \omega_2
  \end{pmatrix}
  \right)
  \ =\ 
  f
  \begin{pmatrix}
  \omega_1
  \\
  \omega_2
  \end{pmatrix}
\]
holds
for every $N \in \Nfrak$.


\paragraph{3. Transformation of modular functions.}

Let $\Afrak$ be a group,
and $\Mfrak$ be
a subgroup of $\Afrak$.
To each $S \in \Mfrak$,
we assign the subgroup
\[
  \Mfrak_S
  \ =\ 
  \Mfrak
  \ \cap\ 
  S^{-1}\,\Mfrak\,S
\]
of $\Mfrak$;
here,
$\Mfrak_S$ depends only on
the residue class
$\Mfrak\,S$
of $S$ modulo $\Mfrak$.
Further,
to each pair of elements
$S,S'$ of $\Afrak$,
we assign
the complex
\[
  \Mfrak_{S,S'}
  \ =\ 
  \Mfrak
  \ \cap\ 
  S^{-1}\,\Mfrak\,S'
  .
\]
Then
$\Mfrak_{S,S'}$
is non-empty
if and only if
$S'$ lies in
$\Mfrak\,S\,\Mfrak$,
in which case
$\Mfrak_{S,S'}$
is a residue class
$\Mfrak_S\,M$
of $\Mfrak$ modulo $\Mfrak_S$.
Conversely,
each residue class
$\Mfrak_S\,M$
of $\Mfrak$ modulo $\Mfrak_S$
is also
a $\Mfrak_{S,S'}$,
namely
\[
  \Mfrak_S\,M
  \ =\ 
  \Mfrak_{S,SM}
  .
\]

One has
$\Mfrak_{S,S''}
 =
 \Mfrak_{S,S'}
$
if and only if
$\Mfrak\,S' = \Mfrak\,S''$,
whence
\[
  \Mfrak_{S,S'}
  \ \longleftrightarrow\ 
  \Mfrak\,S'
\]
is a one-one correspondence
from the residue classes of
$\Mfrak$ modulo $\Mfrak_S$
onto
those residue classes of
$\Afrak$ modulo $\Mfrak$
into which
$\Mfrak\,S\,\Mfrak$
is decomposed.

For $X$ in $\Mfrak$,
one has
\[
  X^{-1}
  \,
  \Mfrak_{S,S'}
  \,
  X
  \ =\ 
  \Mfrak_{SX,S'X}
  ,
\]
in particular
\[
  X^{-1}
  \,
  \Mfrak_S
  \,
  X
  \ =\ 
  \Mfrak_{SX}
  ;
\]
in other words:
\ 
for $S'$ in $\Mfrak\,S\,\Mfrak$,
$\Mfrak_{S'}$
is conjugate to
$\Mfrak_{S}$
in $\Mfrak$,
and every group
which is conjugate to
$\Mfrak_{S}$
in $\Mfrak$
is a $\Mfrak_{S'}$,
$S' \in \Mfrak\,S\,\Mfrak$.

Now let $\Mfrak$ be
the homogeneous modular group
and let $\Afrak$ be
the group of all
two-by-two
rational-entry
matrices
with positive determinants.
Since
for a rational number
$r \neq 0$
one always has
$\Mfrak_{rS} = \Mfrak_{S}$,
we can limit ourselves to
integer-entry $S$
with
relatively prime elements;
such an $S$
is called \emph{primitive}.

To the group-theoretic facts
just established,
two things
are to be be added
in the present case:
\begin{itemize}
\itemsep0pt
\item[1.]
For a given primitive $S$,
$\Mfrak\,S\,\Mfrak$
is the totality $\Afrak_s$
of \emph{all} primitive matrices
which have the same determinant $s$
as $S$;
in other words:
\ 
\textit{
For every pair
$S,S'$ of elements in $\Afrak_s$,
there exist unimodular
$M_{S,S'}, M_{S,S'}^*$
such that
\[
  S\,M_{S,S'}
  \ =\ 
  M_{S,S'}^*\,S'
  .
\]
}
\item[2.]
\textit{
The $\Afrak_s$
decomposes into
a finite number $\psi(s)$ of
residue classes $\Mfrak\,S'$;
a system of representatives
for these residue classes
is formed by
the triangular matrices
}
\[
  \begin{pmatrix}
  a & b
  \\
  0 & d
  \end{pmatrix}
  \qquad
  \begin{array}[c]{l}
  \textit{with}
  \quad 
  a > 0
  ,
  \quad 
  a\,d = s
  ,
  \quad 
  (a,b,d) = 1
  ,
  \quad 
  \textit{and}
  \\
  \textit{$b$ running over
  a full system
  of residues modulo~$d$,}
  \\
  \textit{say
  \ 
  $0 \leqslant b < d$}
  .
  \end{array}
\]
They are denoted by
$S_1,\ldots,S_{\psi(s)}$
in any order.
For a prime number
$s = p$,
one has in particular
$\psi(p) = p{+}1$,
and
we will then take
\[
  P_{\nu}
  \ =\ 
  \begin{pmatrix}
  1 & \nu
  \\
  0 & p
  \end{pmatrix}
  ;
  \quad
  \nu = 1,2,\ldots,p
  ;
  \qquad
  P_{p+1}
  \ =\ 
  \begin{pmatrix}
  p & 0
  \\
  0 & 1
  \end{pmatrix}
\]
as the system of representatives.
\end{itemize}

The triangular form of
the representatives
has the consequence that
the condition~\textup{3}
for a modular form
$f(\omega)$
for $\Mfrak_S$
needs to be checked
only for $M = E$.

The (conjugates of the)
subgroups $\Mfrak_S$ of $\Mfrak$
are called
\emph{transformation groups
of level~$s$};
we call
the classes $\Mfrak\,S$
which are
in one-one correspondence
with them
the \emph{equivalence classes
of level~$s$}.

All this also holds
when we take
$\Mfrak$
as the inhomogeneous modular group
and accordingly replace
$\Afrak$
by the factor group
$\Afrak/\{E,-E\}$,
since
$\{E,-E\}$
is contained in
every transformation group
$\Mfrak_S$.

Let us investigate
the \emph{transformations of the field}
$k_{\Mfrak_S}$
of functions
for the group $\Mfrak_S$.
The $\psi(s)$ isomorphisms
$\lambda_M$
of $k_{\Mfrak_S}$ over $k_\Mfrak$
are
\[
  \lambda_{M_{S,S_\nu}}
  \ =\ 
  \lambda_\nu
  ,
  \qquad
  \text{where
    $\lambda_\nu$
    maps
    $k_{\Mfrak_S}$
    onto
    $k_{\Mfrak_{S_\nu}}$
    .
  }
\]
Functions for the group
$\Mfrak_S$
can be constructed
as follows:
\ 
if
$H_1
 \begin{pmatrix}
 \omega_1
 \\
 \omega_2
 \end{pmatrix}
$
and
$H_2
 \begin{pmatrix}
 \omega_1
 \\
 \omega_2
 \end{pmatrix}
$
are two modular forms
(for the group $\Mfrak$)
of the same weight,
then
\[
  h_S(\omega)
  \ =\ 
  \dfrac{
    H_1
    \left(
    S
    \begin{pmatrix}
      \omega_1
      \\
      \omega_2
    \end{pmatrix}
    \right)
  }{
    H_2
    \begin{pmatrix}
      \omega_1
      \\
      \omega_2
    \end{pmatrix}
  }
\]
belongs to $\Mfrak_{S}$.
In particular,
one can choose
$H_1 = f(\omega)$
and
$H_2 = 1$
as modular function,
thus
\[
  h_S(\omega)
  \ =\ 
  f(S(\omega))
  .
\]
One has
\[
  h_S(\omega)
  \,
  \lambda_\nu
  \ =\ 
  h_{S_\nu}(\omega)
  ,
  \qquad
  f(S(\omega))
  \,
  \lambda_\nu
  \ =\ 
  f(S_\nu(\omega))
  .
\]
If the function
$f(\omega)$ is
a generator of
$k_\Mfrak$,
i.e.~if
$k_\Mfrak = k(f(\omega))$,
then
the $\psi(s)$~many functions
$f(S(\omega))\,\lambda_\nu
 =
 f(S_\nu(\omega))
$
are distinct from one another;
that proves that
\[
  [\,
  k_{\Mfrak_S}
  \,{:}\,
  k_{\Mfrak}
  \,]
  \ =\ 
  [\,
  \Mfrak
  \,{:}\,
  \Mfrak_S
  \,]
  \ =\ 
  \psi(s)
\]
and
\[
  k_{\Mfrak_S}
  \ =\ 
  k(f(\omega),f(S(\omega)))
  ,
\]
and in particular
this holds for
$f(\omega) = j(\omega)$:
\[
  k_{\Mfrak_S}
  \ =\ 
  k(j(\omega),j(S(\omega)))
  ,
\]
from which
it then also follows that
the $h_{S_\nu}(\omega)$ are
a full system of conjugates
of $h_S(\omega)$
over $k_\Mfrak$.

Of course,
$h_S(\omega)$
is entire
if $H_1$ and $H_2^{-1}$ are.

The $q$-expansion of $h_S(\omega)$
can be calculated
from the $q$-expansions
of $H_1$ and $H_2$.

We will now show that
there exists
an algebraic function field
$\Pssf_{\Mfrak_S}$
with the constant field
$\Pssf$
---
where
$\Pssf$ always denotes
the field of
rational numbers
---
from which
$k_{\Mfrak_S}$ arises by
extension of scalars
from $\Pssf$ to $k$,
namely
\[
  \Pssf_{\Mfrak_S}
  \ =\ 
  \Pssf(j(\omega),j(S(\omega)))
  .
\]
For that,
it must be shown that
the coefficients
$A_\nu^{(s)}(j(\omega))$
of the principal polynomial of
$j(S(\omega))$ over $k_\Mfrak$,
\[
  \begin{aligned}
  J_s(X,j(\omega))
  &
  \ =\ 
  \prod_{\nu=1}^{\psi(s)}
  (X - j(S_\nu(\omega)))
  \\
  &
  \ =\ 
  X^{\psi(s)}
  \ +\ 
  A_{\psi(s)-1}^{(s)}(j(\omega))
  \,
  X^{\psi(s)-1}
  \ +\ 
  \cdots
  \ +\ 
  A_{0}^{(s)}(j(\omega))
  \end{aligned}
\]
already lie in $P_{\Mfrak}$.
We derive it
from the fact that
the $q$-expansion of $j(S(\omega))$
has rational coefficients,
and indeed
we similarly prove in general:

\textit{
Let $f_S(\omega)$
be a function
from $k_{\Mfrak_S}$,
whose conjugate
$f_{\left(\begin{smallmatrix}
          s & 0
          \\
          0 & 1
          \end{smallmatrix}\right)}
$
has a $q$-expansion
with coefficients
belonging to
a field $\Lambda$.
Then
the absolutely irreducible polynomial
${\displaystyle
 F(X,Y)
 =
 \sum_{\nu=1}^{N}
 \sum_{\mu=1}^{M}
 c_{\nu\mu}
 X^\nu
 Y^\nu
}$
for which
$F(j(S(\omega)),j(\omega)) = 0$
likewise
has its coefficients
belonging to $\Lambda$
}
(up to
a free constant factor~$\neq 0$).

Namely, in
$F
 \left(\ 
 f_{\left(\begin{smallmatrix}
          s & 0
          \\
          0 & 1
          \end{smallmatrix}\right)}
 (\omega)
 \,,\,
 j(\omega)
 \ \right)
 =
 0
$,
let the $q$-expansions for
$j(\omega)$
and
$f_{\left(\begin{smallmatrix}
          s & 0
          \\
          0 & 1
          \end{smallmatrix}\right)}
 (\omega)
$
be inserted;
that yields
a system of
homogeneous linear equations
for the $c_{\nu\mu}$
with coefficients in $\Lambda$,
which
(up to a factor of
proportionality)
is uniquely solvable,
because
${\displaystyle
 \sum_{\nu=1}^{N}
 \sum_{\mu=1}^{M}
 c_{\nu\mu}
 \,
 f_S(\omega)^\nu
 \,
 j(\omega)^\nu
 =
 0
}$
is the only
linear relation
over $k$
among the
$f_S(\omega)^\nu\,j(\omega)^\mu$,
$\nu = 0,1,\ldots,N$,
$\mu = 0,1,\ldots,M$.
The solutions
$c_{\nu\mu}$
are therefore,
up to a factor of
proportionality,
elements of $\Lambda$.

\textit{
The functions
$f_S(\omega)$
from $\Pssf_{\Mfrak_S}$
are characterized
by the property that
the $q$-expansion of
$f_{\left(\begin{smallmatrix}
          s & 0
          \\
          0 & 1
          \end{smallmatrix}\right)}
 (\omega)
$
has rational coefficients.
}

That
the functions from 
$\Pssf_{\Mfrak_S}$
have this property
is clear,
because
they are
rational functions of
$j(\omega)$
and
$j(S(\omega))$
with rational coefficients.
Conversely,
suppose
$f_{\left(\begin{smallmatrix}
          s & 0
          \\
          0 & 1
          \end{smallmatrix}\right)}
 (\omega)
$
has rational $q$-coefficients.
Since
the irreducible equation
$F(f_S(\omega),j(\omega)) = 0$
for $f_S(\omega)$ in $k_\Mfrak$
has coefficients in
$\Pssf_\Mfrak$,
as shown above,
it follows that
$f_S(\omega)$
is algebraic over
$\Pssf_\Mfrak$.
The field
\[
 \Pssf_{\Mfrak_S}(f_S(\omega))
  = 
 \Pssf(j(\omega),j(S(\omega)),f_S(\omega))
\]
contains no constant
other than
the rational numbers,
because
the $q$-expansion
of a function from
$\Pssf
 \left(\ 
 j(\omega)
 \,,\,
 j(S(\omega))
 \,,\,
 f_{\left(\begin{smallmatrix}
          s & 0
          \\
          0 & 1
          \end{smallmatrix}\right)}
 (\omega)
 \ \right)
$
has rational coefficients
and
the $q$-expansion
of a constant
$\varrho$
reads as
$\varrho = \varrho \cdot q^0$.
One thus has
the degree equation
\[
  [\,
  \Pssf_{\Mfrak_S}(f_S(\omega))
  \,{:}\,
  \Pssf(j(\omega))
  \,]
  \ =\ 
  [\,
  k_{\Mfrak_S}(f_S(\omega))
  \,{:}\,
  k(j(\omega))
  \,]
  ,
\]
but since
$k_{\Mfrak_S}(f_S(\omega))
 =
 k_{\Mfrak_S}
$,
one has
\[
  [\,
  \Pssf_{\Mfrak_S}(f_S(\omega))
  \,{:}\,
  \Pssf(j(\omega))
  \,]
  \ =\ 
  [\,
  k_{\Mfrak_S}
  \,{:}\,
  k_{\Mfrak}
  \,]
  \ =\ 
  [\,
  \Pssf_{\Mfrak_S}
  \,{:}\,
  \Pssf_\Mfrak
  \,]
  ,
\]
whence
\[
  \Pssf_{\Mfrak_S}(f_S(\omega))
  \ =\ 
  \Pssf_{\Mfrak_S}
  ,
  \qquad
  f_S(\omega)
  \ \in\ 
  \Pssf_{\Mfrak_S}
  .
\]

Of course,
a function
from $k_{\Mfrak_S}$
of the form
\[
  h_S(\omega)
  \ =\ 
  \dfrac{
    H_1
    \left(
    S
    \begin{pmatrix}
      \omega_1
      \\
      \omega_2
    \end{pmatrix}
    \right)
  }{
    H_2
    \begin{pmatrix}
      \omega_1
      \\
      \omega_2
    \end{pmatrix}
  }
\]
lies in
$\Pssf_{\Mfrak_S}$
when the $q$-expansions
\[
  H_\nu(\omega)
  \ =\ 
  \sum_{n}
  h_n^{(\nu)}
  \,
  q^n
\]
have rational coefficients
$h_n^{(\nu)}$.

In any case,
the $q$-expansion
of a function
$f_S(\omega)$
from $\Pssf_{\Mfrak_S}$,
for arbitrary
$S
 =
 \begin{pmatrix}
 a & b
 \\
 0 & d
 \end{pmatrix}
$,
has coefficients
in the cyclotomic field
$\Pssf(\zeta_s)$,
$\zeta_s
 =
 e^{\textstyle \frac{2\pi i}{s}}
$,
because
that holds for
\[
 j(S(\omega))
 =
 j
 \left(
 \dfrac{a\,\omega{+}b}{d}
 \right)
 =
 {\displaystyle
 \sum_{n}
 c_n
 \,
 \zeta_s^{abn}
 \,
 q^{{\textstyle \frac{a}{d}}\,n}
 }
 .
\]

A function
$f_S(\omega)$
from
$\Pssf_{\Mfrak_S}$
will be called
\emph{integral}
if its conjugates
$f_{S_\nu}(\omega)$
have $q$-expansions
whose coefficients are
algebraic integers
(in $\Pssf(\zeta_s)$).

The coefficients of
the principal polynomial
\[
  F(X,j(\omega))
  \ =\ 
  \prod_{\nu=1}^{\psi(s)}
  (X - f_{S_\nu}(\omega))
\]
of an integral function
$f_S(\omega)$
from $\Pssf_{\Mfrak_S}$
have rational integer $q$-coefficients,
and from that,
it follows by~\S\textbf{1}
that:

\textit{
The principal polynomial
$F(X,j(\omega))$
of an entire integral function
$f_S(\omega)$
from $\Pssf_{\Mfrak_S}$
is a polynomial
in $X$ and $j(\omega)$
with rational integer coefficients.
}
In particular,

\textit{
$J_s(X,j(\omega))$
is a polynomial
in $X$ and $j(\omega)$
with rational integer coefficients.
}

When the level of transformation
$s = p$
is a prime number,
the principal polynomial of
an entire integral function
from $\Pssf_{\Mfrak_P}$
satisfies
a simple congruence modulo~$p$.
We derive it
first for
$J_s(X,j(\omega))$.

\textit{
As a polynomial
in $X$ and $j(\omega)$,
$J_p(X,j(\omega))$
satisfies
the congruence}
\[
  J_p(X,j(\omega)
  \ \equiv\ 
  (X^p-j(\omega))
  \,
  (X-j(\omega)^p)
  \mod p
  .
\]
By~\S\textbf{1}
it suffices to show that
corresponding coefficients of
the expansions in $X$ and $q$
on both sides
are congruent modulo~$p$.
One has
\[
  J_p(X,j(\omega))
  \ =\ 
  (X - j(p\,\omega))
  \,
  \prod_{\nu=0}^{p-1}
  \left(\,
  X
  -
  j\left(\dfrac{\omega{+}\nu}{p}\right)
  \,\right)
  .
\]
From
\[
  j(\omega)
  \ =\ 
  \sum
  c_n
  \,
  q^n
  ,
\]
it follows that
\[
  j\left(\dfrac{\omega{+}\nu}{p}\right)
  \ =\ 
  \sum
  c_n
  \,
  \zeta_p^{\nu n}
  \,
  q^{\textstyle \frac{n}{p}}
  ,
\]
thus
the coefficient-wise congruence
means
\[
  \begin{gathered}
  j\left(\dfrac{\omega{+}\nu}{p}\right)
  \ \equiv\ 
  j\left(\dfrac{\omega}{p}\right)
  \mod
  (1{-}\zeta_p)
  ,
  \\
  J_p(X,j(\omega))
  \ \equiv\ 
  (X-j(p\,\omega))
  \,
  \left(\,
  X^p
  -
  j\left(\dfrac{\omega}{p}\right)^p
  \,\right)
  \mod
  (1{-}\zeta_p)
  .
  \end{gathered}
\]
Further,
one has
\[
  j(p\,\omega)
  \ =\ 
  \sum
  c_n
  \,
  q^{pn}
  ,
\]
hence
\[
  j(p\,\omega)
  \ \equiv\ 
  j(\omega)^p
  \mod
  p
\]
and
\[
  j(\omega)
  \ \equiv\ 
  j\left(\dfrac{\omega}{p}\right)
  \mod
  p
  ,
\]
thus
\[
  J_p(X,j(\omega))
  \ \equiv\ 
  (X - j(\omega)^p)
  \,
  (X^p - j(\omega))
  \mod
  (1{-}\zeta_p)
  ;
\]
but as the coefficients of
the expansions in $X$ and $q$
from both sides
are rational integers,
this also holds
modulo~$p$.

Now let
$f_P(\omega)$
be an arbitrary
entire integral function
from $\Pssf_{\Mfrak_P}$,
and let
$F(X,j(\omega))$
be its principal polynomial
over $\Pssf_\Mfrak$.

For the $q$-expansion of
$f_{\left(\begin{smallmatrix}
          1 & \nu
          \\
          0 & p
          \end{smallmatrix}\right)}
 (\omega)
$,
$\nu = 0,1,2,\ldots,p{-}1$,
one again has
\[
  f_{\left(\begin{smallmatrix}
          1 & \nu
          \\
          0 & p
          \end{smallmatrix}\right)}
  (\omega)
  \ \equiv\ 
  f_{\left(\begin{smallmatrix}
          1 & 0
          \\
          0 & p
          \end{smallmatrix}\right)}
  (\omega)
  \mod
  (1{-}\zeta_p)
  ,
\]
whence
\[
  F(X,j(\omega))
  \ \equiv\ 
  \left(\,
  X
  -
  f_{\left(\begin{smallmatrix}
          p & 0
          \\
          0 & 1
          \end{smallmatrix}\right)}
  (\omega)
  \,\right)
  \,
  \left(\,
  X^p
  -
  f_{\left(\begin{smallmatrix}
          1 & 0
          \\
          0 & p
          \end{smallmatrix}\right)}
  (\omega)^p
  \,\right)
  \mod
  (1{-}\zeta_p)
\]
in the sense of
the $q$-coefficients.
Let us write
$F(X,j(\omega))$
out as
\[
  F(X,j(\omega))
  \ =\ 
  X^{p+1}
  \ +\ 
  Q_p(j(\omega))
  \,
  X^p
  \ +\ 
  \cdots
  \ +\ 
  Q_0(j(\omega))
  ,
\]
from which
it follows that
\[
  F(X,j(\omega))
  \ \equiv\ 
  X^{p+1}
  \ +\ 
  Q_p(j(\omega))
  \,
  X^p
  \ +\ 
  Q_1(j(\omega))
  \,
  X
  \ +\ 
  Q_0(j(\omega))
  \mod
  p
  ,
\]
i.e.~the polynomials
$Q_\mu(j(\omega))$
of $j(\omega)$,
$2 \leqslant \mu \leqslant p{-}1$,
have coefficients
which are integers
divisible by~$p$.

In addition to $j(S(\omega))$,
we still need
another primitive element of
$\Pssf_{\Mfrak_S}/\Pssf_\Mfrak$,
namely
the function formed by
$H_1 = s^{12}\,\Delta$
and
$H_2 = \Delta$:
\[
  \varphi_S(\omega)
  \ =\ 
  s^{12}
  \,
  \dfrac{
    \Delta
    \left(
    S
    \begin{pmatrix}
      \omega_1
      \\
      \omega_2
    \end{pmatrix}
    \right)
  }{
    \Delta
    \begin{pmatrix}
      \omega_1
      \\
      \omega_2
    \end{pmatrix}
  }
  .
  \tag{15}
\]
That
$\varphi_S(\omega)$
belongs to
$\Pssf_{\Mfrak_S}$
follows from the fact that
the $q$-expansion coefficients
of $\Delta$
are rational.
The $q$-expansions of
the conjugates of
$\varphi_S(\omega)$
yield
\[
  \varphi_{\left(\begin{smallmatrix}
          a & b
          \\
          0 & d
          \end{smallmatrix}\right)}
  (\omega)
  \ =\ 
  a^{12}
  \,
  \zeta_s^{a\,b}
  \,
  q^{\textstyle\frac{a-d}{d}}
  \,
  \dfrac
  { 1
    +
    D_1
    \,
    \zeta_s^{a\,b}
    \,
    q^{\textstyle\frac{a}{d}}
    +
    \cdots
    +
    D_n
    \,
    \zeta_s^{a\,b\,n}
    \,
    q^{\textstyle{\frac{a}{d}}\,n}
    +
    \cdots
  }
  { 1
    +
    D_1
    \,
    q
    +
    \cdots
    +
    D_n
    \,
    q^n
    +
    \cdots
  }
  .
  \tag{16}
  \label{eqn:16}
\]
They show that
the $\psi(s)$~conjugates
$\varphi_{S_\nu}(\omega)$
are distinct from one another;
whence
$\varphi_S(\omega)$
is indeed
a primitive element of
$\Pssf_{\Mfrak_S}/\Pssf_\Mfrak$.

Since
$\Delta$ and $\Delta^{-1}$
are both entire modular forms,
and since
the $q$-expansions~\eqref{eqn:16}
have obviously
integer coefficients,
it follows that
$\varphi_S(\omega)$
is an entire integral function
from $\Pssf_{\Mfrak_S}$,
and similarly for
$j(S(\omega))$.
Its principal polynomial
\[
  \begin{aligned}
  \Phi_s(X,j(\omega))
  &
  \ =\ 
  \prod_S
  (X - \varphi_S(\omega))
  \\
  &
  \ =\ 
  X^{\psi(s)}
  \ +\ 
  B_{\psi(s)-1}^{(s)}(j(\omega))
  \,
  X^{\psi(s)-1}
  \ +\ 
  \cdots
  \ +\ 
  B_{0}^{(s)}(j(\omega))
  \end{aligned}
\]
is therefore also
a polynomial in $X$ and $j(\omega)$
with rational integer coefficients.

The coefficient
$B_{0}^{(s)}(j(\omega))$
can be calculated.
As $\varphi_S(\omega)$
obviously has
no zero
for any finite value of $\omega$,
it follows that
$B_{0}^{(s)}(j(\omega))$
has no zero either,
and thus
$B_0^{(s)}$
as a polynomial in $j(\omega)$
must be a constant.
Its value
is obtained as
the product of
the leading coefficients
of the $q$-series
for all
$\varphi_S(\omega)$:
\[
  \prod_{S}
  \varphi_S(\omega)
  \ =\ 
  (-1)^{\psi(s)}
  \,
  B_0^{(s)}(j(\omega))
  \ =\ 
  \prod_{S}
  a^{12}
  \,
  \zeta_s^{ab}
  ,
\]
and since
this must be
a rational number,
\[
  B_0^{(s)}(j(\omega))
  \ =\ 
  \pm
  \prod_{S}
  a^{12}
  .
  \tag{17}
  \label{eqn:17}
\]

When $s = p$ is
a prime number,
this becomes
\[
  B_0^{(s)}(j(\omega))
  \ =\ 
  (-1)^{p+1}
  \,
  p^{12}
  \,
  \prod_{b=0}^{p-1}
  \zeta_p^b
  \ =\ 
  p^{12}
  .
  \tag{18}
  \label{eqn:18}
\]

To express
a function
$f_S(\omega)$
from $k_{\Mfrak_S}$
in terms of
the basis
$1,
 \varphi_S(\omega),
 \ldots,
 \varphi_S(\omega)^{\psi(s)-1}
$
of $\Pssf_{\Mfrak_S}/\Pssf_\Mfrak$,
we use
the \textit{Lagrange} interpolation formula
in the well-known manner:
\[
  f_S(\omega)
  \,
  \Phi_s'(\varphi_S(\omega),j(\omega))
  \ =\ 
  \Tr_{k_{\Mfrak_S}/k_\Mfrak}
  \left[\ 
    f_S(\omega)
    \,
    \dfrac
    {\Phi_s(X,j(\omega))}
    {X - \varphi_S(\omega)}
  \ \right]_{X = \varphi_S(\omega)}
\]
or
\[
  f_S(\omega)
  \,
  \Phi_s'(\varphi_S(\omega),j(\omega))
  \ =\ 
  a_0(j(\omega))
  \ +\ 
  a_1(j(\omega))
  \,
  \varphi_S(\omega)
  \ +\ 
  \cdots
  \ +\ 
  a_{\psi(s)-1}(j(\omega))
  \,
  \varphi_S(\omega)^{\psi(s)-1}
\]
with
\[
  a_\nu(j(\omega))
  \ =\ 
  \sum_{\mu = \nu+1}^{\psi(s)}
  B_\mu(j(\omega))
  \,
  \Tr_{k_{\Mfrak_S}/k_\Mfrak}
  \left[\ 
    f_S(\omega)
    \,
    \varphi_S(\omega)^{\mu-(\nu+1)}
  \ \right]
  .
\]

From this
we read off:
\begin{itemize}
\itemsep0pt
\item[1.]
\textit{
If $f_S(\omega)$ is
entire,
the
$a_\nu(j(\omega))$
are the polynomials.
}
\item[2.]
\textit{
If $f_S(\omega)$ is
an integral function
from $\Pssf_{\Mfrak_S}$,
the $q$-coefficients
of all $a_\nu(j(\omega))$
are integers.
}
\\
and thus furthermore
\item[3.]
\textit{
If $f_S(\omega)$ is
an entire integral function
from $\Pssf_{\Mfrak_S}$,
the polynomials
$a_\nu(j(\omega))$
in $j(\omega)$
have rational integer coefficients.
}
\end{itemize}

We apply this
to the case when
$s = p$
is a prime number
to obtain
a proposition on
the divisibility of
the coefficients of
$a_0(j(\omega))$
by~$p$:

\textit{
Let $f_P(\omega)$ be
an entire integral function
from $\Pssf_{\Mfrak_P}$,
where $P$ is of
determinant~$p$.
Furthermore,
suppose the coefficients of
the $q$-expansion of
$\smash[b]{
 f_{\left(\begin{smallmatrix}
         p & 0
         \\
         0 & 1
         \end{smallmatrix}\right)}
 (\omega)
}$
are divisible by~$p$.
Then
the coefficients of
$a_0(j(\omega))$
are divisible by~$p$.
}
\footnote{
This theorem
summarizes
an argument
given by
\textit{Hasse},
J.~f.~Math.
\textbf{157},
125--126 (1927).
}

Namely,
one has
\[
  f_P(\omega)
  \,
  \Phi_p'(\varphi_P(\omega),j(\omega))
  \ =\ 
  a_0(j(\omega))
  \ +\ 
  a_1(j(\omega))
  \,
  \varphi_P(\omega)
  \ +\ 
  \cdots
  \ +\ 
  a_p(j(\omega))
  \,
  \varphi_P(\omega)^p
\]
for
$P
 =
 \begin{pmatrix}
   p & 0
   \\
   0 & 1
 \end{pmatrix}
$;
since
$f_{\left(\begin{smallmatrix}
         p & 0
         \\
         0 & 1
         \end{smallmatrix}\right)}
 (\omega)
$
and
\[
  \varphi_{\left(\begin{smallmatrix}
         p & 0
         \\
         0 & 1
         \end{smallmatrix}\right)}
  (\omega)
  \ =\ 
  p^{12}
  \,
  q^{p-1}
  \,
  \dfrac
  {1 + D_1\,q^p + \cdots}
  {1 + D_1\,q + \cdots}
\]
are divisible by~$p$
and
the $a_\nu(j(\omega))$
have integer $q$-coefficients,
it follows that
the $q$-coefficients of
$a_0(j(\omega))$
are divisible by~$p$.


\paragraph{4. The functions $\sqrt[24]{\varphi_S(\omega)}$.}

The twenty-fourth root of
$\Delta(\omega)$,
that is
\[
  \sqrt[24]{\Delta(\omega)}
  \ =\ 
  \sqrt{2\,\pi}
  \,
  \eta(\omega)
  ,
\]
is a function of
the variable $\omega$;
it is regular in
$\Im\,\omega > 0$.
It is actually
not a modular form
for $\Mfrak$,
but it is well-known
how it transforms
when a modular substitution
is performed on $\omega$.
\footnote{
This transformation formula
for the $\eta$-function
was in essence
already stated by
\textit{Dedekind}.
}
Namely,
for
$M
 =
 \begin{pmatrix}
   \alpha & \beta
   \\
   \gamma & \delta
 \end{pmatrix}
$
with $\gamma \geqslant 0$
from $\Mfrak$,
one has
\[
  \eta(M(\omega))
  \ =\ 
  \varepsilon(M)
  \,
  \sqrt{-i\,(\gamma\,\omega{+}\delta)}
  \,
  \eta(\omega)
  ,
\]
where the root
is taken with
positive real part
and $\varepsilon(M)$ denotes
the following
24-th root of unity:
\ 
let
$\gamma
 =
 2^\lambda
 \,
 \gamma_1
$
with $\gamma_1$ odd
and
$\gamma_1 = 1$
in the case when
$\gamma = 0$;
then
\[
  \varepsilon 
  \begin{pmatrix}
    \alpha & \beta
    \\
    \gamma & \delta
  \end{pmatrix}
  \ =\ 
  \left(
  \dfrac{\alpha}{\gamma_1}
  \right)
  \,
  \zeta_{24}^{
    \left(\,
    \beta\,\delta\,(1-\gamma^2)
    +
    \gamma\,(\alpha+\delta)
    +
    3\,(1-\gamma_1)
    +
    3\,\alpha\,(\gamma-\gamma_1)
    +
    \lambda
    \cdot
    \frac{3}{2}
    \,
    (\alpha^2-1)
    \,\right)
  }
  ,
  \qquad
  \zeta_{24}
  \ =\ 
  e^{\textstyle\frac{2\pi i}{24}}
  .
\]

We now define
the 24-th root of
$\varphi_{\left(\begin{smallmatrix}
        s & 0
        \\
        0 & 1
        \end{smallmatrix}\right)}
 (\omega)
$
by
\[
  \sqrt[24]{
    \varphi_{\left(\begin{smallmatrix}
          s & 0
          \\
          0 & 1
          \end{smallmatrix}\right)}
    (\omega)
  }
  \ =\ 
  \dfrac
  {\sqrt[24]{
    s^{12}
    \,
    \Delta
    \begin{pmatrix}
      s\,\omega_1
      \\
      \omega_2
    \end{pmatrix}
  }}
  {\sqrt[24]{
    \Delta
    \begin{pmatrix}
      \omega_1
      \\
      \omega_2
    \end{pmatrix}
  }}
  \ =\ 
  \sqrt{s}
  \,
  \dfrac{\eta(s\,\omega)}{\eta(\omega)}
  ,
  \qquad
  \sqrt{s} > 0
  .
\]
How does
$\sqrt[24]{
 \varphi_{\left(\begin{smallmatrix}
       s & 0
       \\
       0 & 1
       \end{smallmatrix}\right)}
 (\omega)
 }
$
transform
when an element $M$ of
$\Mfrak_{\left(\begin{smallmatrix}
       s & 0
       \\
       0 & 1
       \end{smallmatrix}\right)}
$
is applied to $\omega$?
Since
\[
  \begin{pmatrix}
  s & 0
  \\
  0 & 1
  \end{pmatrix}
  \,
  \begin{pmatrix}
  \alpha & \beta
  \\
  \gamma & \delta
  \end{pmatrix}
  \,
  \begin{pmatrix}
  s & 0
  \\
  0 & 1
  \end{pmatrix}^{-1}
  \ =\ 
  \begin{pmatrix}
  \alpha & s\,\beta
  \\
  s^{-1}\,\gamma & \delta
  \end{pmatrix}
  ,
\]
the elements of
$\Mfrak_{\left(\begin{smallmatrix}
       s & 0
       \\
       0 & 1
       \end{smallmatrix}\right)}
$
are
the modular substitutions
$\begin{pmatrix}
 \alpha & \beta
 \\
 \gamma & \delta
 \end{pmatrix}
$
with $\gamma \equiv 0 \mod s$.
Therefore,
\[
  \sqrt[24]{
    \varphi_{\left(\begin{smallmatrix}
          s & 0
          \\
          0 & 1
          \end{smallmatrix}\right)}
    \left(\,
      \begin{pmatrix}
        \alpha & \beta
        \\
        \gamma & \delta
      \end{pmatrix}
      (\omega)
    \,\right)
  }
  \ =\ 
  \sqrt{s}
  \ 
  \dfrac
  {\eta
    \left(\,
      \begin{pmatrix}
        \alpha & s\,\beta
        \\
        s^{-1}\,\gamma & \delta
      \end{pmatrix}
      (s\,\omega)
    \,\right)
  }
  {\eta
    \left(\,
      \begin{pmatrix}
        \alpha & \beta
        \\
        \gamma & \delta
      \end{pmatrix}
      (\omega)
    \,\right)
  }
  \ =\ 
  \dfrac
  {\varepsilon
      \begin{pmatrix}
        \alpha & s\,\beta
        \\
        s^{-1}\,\gamma & \delta
      \end{pmatrix}
  }
  {\varepsilon
      \begin{pmatrix}
        \alpha & \beta
        \\
        \gamma & \delta
      \end{pmatrix}
  }
  \,
  \sqrt[24]{
    \varphi_{\left(\begin{smallmatrix}
          s & 0
          \\
          0 & 1
          \end{smallmatrix}\right)}
    (\omega)
  }
  .
\]

From this
it is now easy to conclude:
\ 
if $s = t^2$
is a square
relatively prime to~6,
then for
$M
 \in
 \Mfrak_{\left(\begin{smallmatrix}
       t^2 & 0
       \\
       0 & 1
       \end{smallmatrix}\right)}
$
one has:
\[
    \sqrt[24]{
    \varphi_{\left(\begin{smallmatrix}
          t^2 & 0
          \\
          0 & 1
          \end{smallmatrix}\right)}
    (M\,\omega)
  }
  \ =\ 
  \sqrt[24]{
    \varphi_{\left(\begin{smallmatrix}
          t^2 & 0
          \\
          0 & 1
          \end{smallmatrix}\right)}
    (\omega)
  }
  .
\]
One has to prove
\[
  \varepsilon
  \begin{pmatrix}
    \alpha & t^2\,\beta
    \\
    t^{-2}\,\gamma & \delta
  \end{pmatrix}
  \ =\ 
  \varepsilon
  \begin{pmatrix}
    \alpha & \beta
    \\
    \gamma & \delta
  \end{pmatrix}
  .
  \tag{19}
  \label{eqn:19}
\]
For both $\varepsilon$-values,
$\lambda$
is the same number,
since $s$ is odd.
The exponent of
$\zeta_{24}$
in
$\varepsilon
 \begin{pmatrix}
   \alpha & t^2\,\beta
   \\
   t^{-2}\,\gamma & \delta
 \end{pmatrix}
$
thus arises from
that in
$\varepsilon
 \begin{pmatrix}
   \alpha & \beta
   \\
   \gamma & \delta
 \end{pmatrix}
$
by replacing
$\beta$ by $t^2\,\beta$,
$\gamma$ by $t^{-2}\,\gamma$,
and
$\gamma_1$ by $t^{-2}\,\gamma_1$.
Since
$t^2 \equiv 1 \mod 24$,
the same power of
$\zeta_{24}$
appears in
$\varepsilon
 \begin{pmatrix}
   \alpha & t^2\,\beta
   \\
   t^{-2}\,\gamma & \delta
 \end{pmatrix}
$
as in
$\varepsilon
 \begin{pmatrix}
   \alpha & \beta
   \\
   \gamma & \delta
 \end{pmatrix}
$.
Since
$\left(
 \dfrac{\alpha}{t^{-2}\,\gamma_1}
 \right)
 =
 \left(
 \dfrac{\alpha}{\gamma_1}
 \right)
$
as well,
\eqref{eqn:19}
follows.

Still keeping in mind that
the $q$-expansion of
$\eta(\omega)$
and hence also that of
$\sqrt[24]{
  \varphi_{\left(\begin{smallmatrix}
        t^2 & 0
        \\
        0 & 1
        \end{smallmatrix}\right)}
  (\omega)
 }
$
have rational coefficients,
we thus have:
\ 
\textit{
For $t$ relatively prime to~6,
$\sqrt[24]{
  \varphi_{\left(\begin{smallmatrix}
        t^2 & 0
        \\
        0 & 1
        \end{smallmatrix}\right)}
  (\omega)
 }
$
is a function from
$\Pssf_{\Mfrak_{\left(\begin{smallmatrix}
        t^2 & 0
        \\
        0 & 1
        \end{smallmatrix}\right)}}
$.
}

For any primitive $S$
of determinant $t^2$,
we define
$\sqrt[24]{
  \varphi_S(\omega)
 }
$
as the function
from the field
$\Pssf_{\Mfrak_S}$
which is
conjugate to
$\sqrt[24]{
  \varphi_{\left(\begin{smallmatrix}
        t^2 & 0
        \\
        0 & 1
        \end{smallmatrix}\right)}
  (\omega)
 }
$.
The function
$\sqrt[24]{
  \varphi_{\left(\begin{smallmatrix}
        t^2 & 0
        \\
        0 & 1
        \end{smallmatrix}\right)}
  (\omega)
 }
$
is entire and integral,
since that holds for
$\varphi_{\left(\begin{smallmatrix}
        t^2 & 0
        \\
        0 & 1
        \end{smallmatrix}\right)}
 (\omega)
$.


\paragraph{5. Elliptic functions.}

We need
the Weierstass $\wp$-function
\[
  \wp
  \left(
    z
    \,,\,
    \begin{pmatrix}
      \omega_1
      \\
      \omega_2
    \end{pmatrix}
  \right)
  \ =\ 
  z^{-2}
  \ +\ 
  \sum_{\substack{
    n_1,n_2=-\infty
    \\
    n_1,n_2 \neq 0,0
    }}^{+\infty}
  [\,
    (z - n_1\,\omega_1 - n_2\,\omega_2)^{-2}
    \ -\ 
    (n_1\,\omega_1 + n_2\,\omega_2)^{-2}
  \,]
\]
as a function
of the three complex variables
$z,\omega_1,\omega_2$,
in which
we assume that
$\Im(\omega) > 0$,
which amounts to
no restriction.
Since $\wp$
as a function of
$\omega_1,\omega_2$
depends only on
the lattice
$\wfrak$
with the basis
$\omega_1,\omega_2$,
we also
occasionally write
$\wp(z,\wfrak)$.
In
$\wp
 \left(
   z
   \,,\,
   \begin{pmatrix}
     \omega_1
     \\
     \omega_2
   \end{pmatrix}
 \right)
$,
$\begin{pmatrix}
  \omega_1
  \\
  \omega_2
  \end{pmatrix}
$
is understood as
a column vector.
The function $\wp$
is homogeneous in
$z,\omega_1,\omega_2$
of degree~$-2$.
Since $\wp$
as a function of
$\omega = \dfrac{\omega_1}{\omega_2}$
has period~1,
$\wp$ has a Fourier expansion
in powers of
$q = e^{2\pi i \omega}$;
with
$U = e^{\textstyle \frac{2\pi i z}{\omega_2}}$,
it reads as
\[
  \wp
  \left(
    z
    \,,\,
    \begin{pmatrix}
      \omega_1
      \\
      \omega_2
    \end{pmatrix}
  \right)
  \ =\ 
  \,-\,
  \dfrac
  {\left(\frac{2\pi}{\omega_2}\right)^2}
  {12}
  \,
  \left[\ 
  1
  +
  \dfrac
  {12}
  {\left(
   U^{\frac{1}{2}}
   -
   U^{-\frac{1}{2}}
   \right)^2
  }
  -
  24
  \,
  \sum_{n,m=1}^{\infty}
  n
  \,
  q^{nm}
  +
  12
  \,
  \sum_{n,m=1}^{\infty}
  n
  \,
  q^{nm}
  \,
  (U^n + U^{-n})
  \ \right]
  ,
\]
valid for
$\Im\,\omega > 0$
(i.e.~$|q| < 1$)
and
$|q| < U < \dfrac{1}{|q|}$,
i.e.
$\left|\,
 \Im
 \left(\dfrac{z}{\omega_2}\right)
 \,\right|
 >
 \Im\,\omega
$.


\paragraph{6. The Weber function.}

For number-theoretic purposes
the $\wp$-function
must be replaced by
a suitably normalized
elliptic function
which is adapted to
an imaginary quadratic number field
$\Sigma$
and an order $\Rssf$ in it.
To the order $\Rssf$
is assigned
an elliptic function
$\tau_\Rssf
 \left(
   z
   \,,\,
   \begin{pmatrix}
     \omega_1
     \\
     \omega_2
   \end{pmatrix}
 \right)
$,
the \emph{Weber function of
the order $\Rssf$},
namely
\[
  \tau_\Rssf
  \left(
    z
    \,,\,
    \begin{pmatrix}
      \omega_1
      \\
      \omega_2
    \end{pmatrix}
  \right)
  \ =\ 
  g^{(e)}
  \begin{pmatrix}
    \omega_1
    \\
    \omega_2
  \end{pmatrix}
  \,
  \wp
  \left(
    z
    \,,\,
    \begin{pmatrix}
      \omega_1
      \\
      \omega_2
    \end{pmatrix}
  \right)^{\frac{1}{2}\,e}
  ,  
  \tag{20}
\]
in which
$e$ denotes
the number of units
of $\Rssf$;
thus
$e = 2$
except when
$\Rssf$ is the principal order of
$\Pssf(\sqrt{-1})$
where $e = 4$,
or when
$\Rssf$ is the principal order of
$\Pssf(\sqrt{-3})$
where $e = 6$;
and $g^{(e)}$
is a certain
entire modular form
of weight~$e$,
namely
\[
  \left\{
  \qquad
  \begin{aligned}
  g^{(2)}
  \begin{pmatrix}
    \omega_1
    \\
    \omega_2
  \end{pmatrix}
  &
  \ =\ 
  -
  2^7 \cdot 3^5
  \,
  \dfrac{g_2(\omega)\,g_3(\omega)}{\Delta(\omega)}
  \,
  \omega_2^2
  \\
  g^{(4)}
  \begin{pmatrix}
    \omega_1
    \\
    \omega_2
  \end{pmatrix}
  &
  \ =\ 
  \phantom{-}
  2^8 \cdot 3^4
  \,
  \dfrac{g_2^2(\omega)}{\Delta(\omega)}
  \,
  \omega_2^4
  \\
  g^{(6)}
  \begin{pmatrix}
    \omega_1
    \\
    \omega_2
  \end{pmatrix}
  &
  \ =\ 
  -
  2^9 \cdot 3^6
  \,
  \dfrac{g_3(\omega)}{\Delta(\omega)}
  \,
  \omega_2^6
  .
  \end{aligned}
  \qquad
  \right.
  \tag{21}
\]
The function
$\tau_\Rssf
 \left(
   z
   \,,\,
   \begin{pmatrix}
     \omega_1
     \\
     \omega_2
   \end{pmatrix}
 \right)
$
is
\emph{homogenous of degree~0}
as a function of
the three complex variables
$z,\omega_1,\omega_2$,
and as a function of
$\omega_1,\omega_2$,
it is invariant
with respect to
modular substitutions.
The $g^{(e)}$ is chosen
so that
$\tau_\Rssf$ is not
identically~0
when a basis of
an ideal of $\Rssf$
is used for
$\omega_1,\omega_2$;
the numerical factor in
$g^{(e)}$
causes
the \textit{Fourier} expansion of
$\tau_\Rssf$
to have the form
\[
  \begin{aligned}
  &
  \tau_\Rssf
  \left(
   z
   \,,\,
   \begin{pmatrix}
     \omega_1
     \\
     \omega_2
   \end{pmatrix}
  \right)
  \ =\ 
  \\
  &
  \left[\ 
  q^{\textstyle \frac{e}{2}}
  \ +\ 
  \cdots
  \ +\ 
  t_\nu\,q^\nu
  \ +\ 
  \cdots
  \ \right]
  \,
  \left[\ 
  1
  \ +\ 
  \dfrac{12\,U}{(1-U)^2}
  \ +\ 
  12
  \,
  \sum_{n,m=1}^{\infty}
  n
  \,
  q^{nm}
  \,
  (U^n + U^{-n} - 2)
  \ \right]^{\textstyle \frac{e}{2}}
  ,
  \end{aligned}
  \tag{22}
  \label{eqn:22}
\]
where the coefficients
$t_\nu$
are rational integers.

We can also write
$\tau_\Rssf(z,\wfrak)$
instead of
$\tau_\Rssf
 \left(
   z
   \,,\,
   \begin{pmatrix}
     \omega_1
     \\
     \omega_2
   \end{pmatrix}
 \right)
$
where $\wfrak$ denotes
the lattice
with the basis
$\omega_1,\omega_2$,
since
$\tau_\Rssf$
depends only on $\wfrak$
and not on
the particular basis.


\paragraph{7. The division values of the Weber function.}
\footnote{
\textit{Hasse},
loc. cit.
{\ \ \llap{${}^{2)}$}}
131 to 136.
}
Let $N$ be
a natural number.
By the
\emph{$N$-th division values}
of
$\tau_\Rssf
 \left(
   z
   \,,\,
   \begin{pmatrix}
     \omega_1
     \\
     \omega_2
   \end{pmatrix}
 \right)
$,
we mean
the functions
\[
  \tau_\Rssf
  \left(
    \dfrac{x_1\,\omega_1 + x_2\,\omega_2}{N}
    \,,\,
    \begin{pmatrix}
      \omega_1
      \\
      \omega_2
    \end{pmatrix}
  \right)
  \ =\ 
  \tau_\Rssf
  \left(
    \dfrac{1}{N}
    \,
    (x_1,x_2)
    \begin{pmatrix}
      \omega_1
      \\
      \omega_2
    \end{pmatrix}
    \,,\,
    \begin{pmatrix}
      \omega_1
      \\
      \omega_2
    \end{pmatrix}
  \right)
  ,
\]
where $x_1,x_2$ are integers
with
\[
  (x_1,x_2)
  \ \not\equiv\ 
  (0,0)
  \ 
  \mod
  N
  .
\]
It only depends on
$x_1,x_2$
modulo~$N$.

If
$(x_1,x_2,N) = 1$,
then
$\tau_\Rssf
 \left(
   \dfrac{x_1\,\omega_1 + x_2\,\omega_2}{N}
   \,,\,
   \begin{pmatrix}
     \omega_1
     \\
     \omega_2
   \end{pmatrix}
 \right)
$
is called a \emph{proper}
$N$-th division value
of $\tau_\Rssf$.
Obviously
$\tau_\Rssf
 \left(
   \dfrac{x_1\,\omega_1 + x_2\,\omega_2}{N}
   \,,\,
   \begin{pmatrix}
     \omega_1
     \\
     \omega_2
   \end{pmatrix}
 \right)
$
is a proper $N_1$-th division value,
for
$N_1 = N/(x_1,x_2,N)$.

\textit{
The division values of
$\tau_\Rssf$
are regular functions of $\omega$
in $\Im\,\omega > 0$.
}
Let $M$ be
a modular substitution.
We apply $M$ to $\omega$;
then
$\tau_\Rssf
 \left(
   \dfrac{x_1\,\omega_1 + x_2\,\omega_2}{N}
   \,,\,
   \begin{pmatrix}
     \omega_1
     \\
     \omega_2
   \end{pmatrix}
 \right)
$
becomes
\[
  \tau_\Rssf
  \left(
    \dfrac{1}{N}
    \,
    (x_1,x_2)
    \cdot
    M
    \begin{pmatrix}
      \omega_1
      \\
      \omega_2
    \end{pmatrix}
    \,,\,
    M
    \begin{pmatrix}
      \omega_1
      \\
      \omega_2
    \end{pmatrix}
  \right)
  \ =\ 
  \tau_\Rssf
  \left(
    \dfrac{1}{N}
    \,
    (x_1,x_2)
    M
    \cdot
    \begin{pmatrix}
      \omega_1
      \\
      \omega_2
    \end{pmatrix}
    \,,\,
    \begin{pmatrix}
      \omega_1
      \\
      \omega_2
    \end{pmatrix}
  \right)
  .
\]
Every $M$ in $\Mfrak$
thus
interchanges
the proper $N$-th division values
of $\tau_\Rssf$
with one another.

The $N$-th division values
of $\tau_\Rssf$
can be expanded
in powers of
$q^{\textstyle \frac{1}{N}}$;
we just have to
substitute
$U
 =
 e^{\textstyle
   \frac{2\pi i}{N}
   \,
   \scriptstyle
   (x_1\,\omega + x_2)
 }
 =
 q^{\textstyle
   \frac{x_1}{N}
 }
 \,
 \zeta_N^{x_2}
$
in~\eqref{eqn:22},
whereby
we choose
$x_1$ in its residue class
modulo~$N$
satisfying
$0 \leqslant x_1 < N$
so as to fulfill
the convergence condition
of~\eqref{eqn:22}.
That gives
\[
  \begin{gathered}[t]
  \tau_\Rssf
  \left(
   \dfrac{x_1\,\omega_1 + x_2\,\omega_2}{N}
   \,,\,
   \begin{pmatrix}
     \omega_1
     \\
     \omega_2
   \end{pmatrix}
  \right)
  \ =\ 
  \qquad\qquad\qquad\qquad\qquad
  \qquad\qquad\qquad\qquad\qquad
  \\
  \mathclap{
  \left[\ 
  q^{\textstyle -\frac{e}{2}}
  \ +\ 
  \cdots
  \ \right]
  \,
  \left[\ 
  1
  \ +\ 
  \dfrac
  {12\,\zeta_N^{x_2}\,q^{\textstyle \frac{x_1}{N}}}
  {\left(
     1 - \zeta_N^{x_2}\,q^{\textstyle \frac{x_1}{N}}
   \right)^2
  }
  \ +\ 
  12
  \,
  \sum_{n,m=1}^{\infty}
  n
  \,
  q^{nm}
  \,
  \left(\,
    \zeta_N^{x_2\,n}\,q^{\textstyle \frac{x_1}{N}\,n}
    +
    \zeta_N^{-x_2\,n}\,q^{\textstyle -\frac{x_1}{N}\,n}
    -
    2
  \,\right)
  \ \right]^{\textstyle \frac{e}{2}}
  .
  }
  \end{gathered}
  \tag{23}
  \label{eqn:23}
\]

In this series
only finitely many
powers of $q$
with negative exponents
occur.
The coefficients
lie in $\Pssf(\zeta_N)$.
When we apply
the automorphism
$\zeta_N \to \zeta_N^r$,
$(r,N)=1$,
of the field
$\Pssf(\zeta_N)$
to the coefficients of
the $q$-expansion of
$\tau_\Rssf
 \left(
  \dfrac{x_1\,\omega_1 + x_2\,\omega_2}{N}
  \,,\,
  \begin{pmatrix}
    \omega_1
    \\
    \omega_2
  \end{pmatrix}
 \right)
$,
it becomes
the $q$-expansion of
$\tau_\Rssf
 \left(
  \dfrac{x_1\,\omega_1 + r\,x_2\,\omega_2}{N}
  \,,\,
  \begin{pmatrix}
    \omega_1
    \\
    \omega_2
  \end{pmatrix}
 \right)
$.

It now follows
from the above
that:
\[
  \prod_{\substack{
    x_1,x_2 \bmod N
    \\
    (x_1,x_2,N) = 1
  }}
  \left(\ 
    X
    \ -\ 
    \tau_\Rssf
    \left(
     \dfrac{x_1\,\omega_1 + x_2\,\omega_2}{N}
     \,,\,
     \begin{pmatrix}
       \omega_1
       \\
       \omega_2
     \end{pmatrix}
    \right)
  \ \right)
\]
is a polynomial
$T_N(X,j(\omega))$
in $X$ and $j(\omega)$
with rational coefficients.
The polynomial
$T_N(X,j(\omega))$
is called
the \emph{$N$-th order division polynomial
of the function
$\tau_\Rssf$}.

We further note that
the series~\eqref{eqn:23}
has integer coefficients
when
$x_1 \not\equiv 0
 \pmod{N}
$,
whereas
it becomes integral
only after multiplication by
$(1-\zeta_N^{x_2})^e$
when
$x_1 \equiv 0
 \pmod{N}
$.
Now one has
\[
  \prod_{\substack{
    x \bmod N
    \\
    (x,N) = 1
  }}
  (1 - \zeta_N^{x})
  \ =\ 
  \begin{cases}
  1
  ,
  &
  \quad
  \text{if $N$ is not a prime power}
  ,
  \\
  \ell
  ,
  &
  \quad
  \text{if $N$ is a power of the prime number $\ell$}
  .
  \end{cases}
\]

\textit{
Consequently
the coefficients of the
divison polynomials
$T_N(X,j(\omega))$
of order~$N$
are rational integers
when $N$ is not a prime power;
and
if $N$ is a power of
the prime number $\ell$,
then
they become integers
after multiplication by
$\ell^e$.
}

Let $p$ be
a prime number
which does not divide~$N$.
With a matrix $P$
of determinant~$p$,
we form the function
\[
  \delta_P((x_1,x_2);\omega)
  \ =\ 
  \tau_\Rssf
  \left(
   \dfrac{1}{N}
   \,
   (x_1,x_2)
   \begin{pmatrix}
     \omega_1
     \\
     \omega_2
   \end{pmatrix}
   \,,\,
   \begin{pmatrix}
     \omega_1
     \\
     \omega_2
   \end{pmatrix}
  \right)^p
  \ -\ 
  \tau_\Rssf
  \left(
   \dfrac{p}{N}
   \,
   (x_1,x_2)
   \begin{pmatrix}
     \omega_1
     \\
     \omega_2
   \end{pmatrix}
   \,,\,
   P
   \begin{pmatrix}
     \omega_1
     \\
     \omega_2
   \end{pmatrix}
  \right)
  .
\]
The function $\delta_P$
depends only on
the equivalence class of $P$
and is regular in
$\Im\,\omega > 0$.
Note that
\[
  \tau_\Rssf
  \left(
   \dfrac{p}{N}
   \,
   (x_1,x_2)
   \begin{pmatrix}
     \omega_1
     \\
     \omega_2
   \end{pmatrix}
   \,,\,
   P
   \begin{pmatrix}
     \omega_1
     \\
     \omega_2
   \end{pmatrix}
  \right)
  \ =\ 
  \tau_\Rssf
  \left(
   \dfrac{1}{N}
   \cdot
   (x_1,x_2)
   \,
   p
   \,
   P^{-1}
   \cdot
   P
   \begin{pmatrix}
     \omega_1
     \\
     \omega_2
   \end{pmatrix}
   \,,\,
   P
   \begin{pmatrix}
     \omega_1
     \\
     \omega_2
   \end{pmatrix}
  \right)
\]
is an $N$-th division value of
$\tau_\Rssf
 \left(
  z
  \,,\,
  P
  \begin{pmatrix}
    \omega_1
    \\
    \omega_2
  \end{pmatrix}
 \right)
$.

If $M$ is
a modular matrix,
one has
\[
  \delta_P((x_1,x_2)\,{;}\,M(\omega))
  \ =\ 
  \delta_{MP}((x_1,x_2)M\,{;}\,\omega)
  .
  \tag{24}
\]

If $M$ in fact
belongs to
the group $\Mfrak_P$,
then $M$ only interchanges
those $\delta_P((x_1,x_2);\omega)$
with
$(x_1,x_2,N) = 1$
with one another.
We note that
$\delta_P((x_1,x_2);\omega)$
has a $q$-expansion
with only finitely many
negative exponents of $q$,
because that holds
for the $N$-th division value of
$\tau_\Rssf
 \left(
  z
  \,,\,
  \begin{pmatrix}
    \omega_1
    \\
    \omega_2
  \end{pmatrix}
 \right)
$
and also for that of
$\tau_\Rssf
 \left(
  z
  \,,\,
  P
  \begin{pmatrix}
    \omega_1
    \\
    \omega_2
  \end{pmatrix}
 \right)
$;
hence we have:

\textit{
The coefficients
$D_P^{(\nu)}(\omega)$
of the polynomial
\[
  S_P(X,\omega)
  \ =\ 
  \prod_{\substack{
    x_1,x_2 \bmod N
    \\
    (x_1,x_2,N) = 1
  }}
  \left(\,
    X
    \ -\ 
    \delta_P((x_1,x_2);\omega)
  \,\right)
  \ =\ 
  \sum_{\nu}
  D_P^{(\nu)}(\omega)
  \,
  X^\nu
\]
are entire functions
from $k_{\Mfrak_P}$.
}
By the way,
as already expressed
in the notation,
$D_{P_\mu}^{(\nu)}(\omega)$
is conjugate to
$D_P^{(\nu)}(\omega)$
as a function from
$k_{\Mfrak_{P_\mu}}$,
for this conjugate
can be calculated as
\[
  D_P^{(\nu)}(\omega)
  \,
  \lambda_\mu
  \ =\ 
  D_P^{(\nu)}(M_{P,P_\mu}(\omega))
  \ =\ 
  D_{P\,M_{P,P_\mu}}^{(\nu)}(\omega)
  \ =\ 
  D_{M_{P,P_\mu}^*\,P_\mu}^{(\nu)}(\omega)
  \ =\ 
  D_{P_\mu}^{(\nu)}(\omega)
  .
\]

We now prove:
\textit{
$D_P^{(\nu)}(\omega)$
in fact
lies in
$\Pssf_{\Mfrak_P}$.
}
\\
For that
we investigate
the $q$-expansion coefficients
of
\[
  \begin{aligned}
  &
  \delta_{\left(\begin{smallmatrix}
         p & 0
         \\
         0 & 1
         \end{smallmatrix}\right)}
  ((x_1,x_2)\,{;}\,\omega)
  \ =\
  \\
  &
  \qquad
  \qquad
  \tau_\Rssf
  \left(
   \dfrac{1}{N}
   \,
   (x_1,x_2)
   \begin{pmatrix}
     \omega_1
     \\
     \omega_2
   \end{pmatrix}
   \,,\,
   \begin{pmatrix}
     \omega_1
     \\
     \omega_2
   \end{pmatrix}
  \right)^p
  \ -\ 
  \tau_\Rssf
  \left(
   \dfrac{1}{N}
   \,
   (x_1,x_2\,p)
   \begin{pmatrix}
     p\,\omega_1
     \\
     \omega_2
   \end{pmatrix}
   \,,\,
   \begin{pmatrix}
     p\,\omega_1
     \\
     \omega_2
   \end{pmatrix}
  \right)
  .
  \end{aligned}
\]
These coefficients
obviously lie in
the cyclotomic field
$\Pssf(\zeta_N)$.
We apply to them simultaneously
the automorphisms
$\zeta_N \to \zeta_N^r$,
$(r,N) = 1$,
of $\Pssf(\zeta_N)$;
then
the $q$-expansion of
$\delta_{\left(\begin{smallmatrix}
        p & 0
        \\
        0 & 1
        \end{smallmatrix}\right)}
 ((x_1,x_2);\omega)
$
becomes that of
$\delta_{\left(\begin{smallmatrix}
        p & 0
        \\
        0 & 1
        \end{smallmatrix}\right)}
 ((x_1,r\,x_2);\omega)
$.
The $q$-expansion coefficients of
$D_{\left(\begin{smallmatrix}
        p & 0
        \\
        0 & 1
        \end{smallmatrix}\right)}
 ^{(\nu)}
 (\omega)
$
are therefore
invariant under all
automorphisms of $\Pssf(\zeta_N)$,
and so they are
rational numbers.
But from that
the claim
$D_P^{(\nu)}(\omega)
 \in
 \Pssf_{\Mfrak_P}
$
follows by~\S\textbf{3}.

\textit{
The function
$D_P^{(\nu)}(\omega)$
has $p$-integral
$q$-expansion coefficients,
}
because
that holds for
the $N$-th division value
of $\tau_\Rssf$.
But furthermore,
we have:
\textit{
The $q$-expansion coefficients of
$D_{\left(\begin{smallmatrix}
        p & 0
        \\
        0 & 1
        \end{smallmatrix}\right)}
 ^{(\nu)}
 (\omega)
$
are divisible by $p$.
}
Namely,
one has
\[
  \delta_{\left(\begin{smallmatrix}
         p & 0
         \\
         0 & 1
         \end{smallmatrix}\right)}
  ((x_1,x_2)\,{;}\,\omega)
  \ =\
  \tau_\Rssf
  \left(
   \dfrac{1}{N}
   \,
   (x_1,x_2)
   \begin{pmatrix}
     \omega_1
     \\
     \omega_2
   \end{pmatrix}
   \,,\,
   \begin{pmatrix}
     \omega_1
     \\
     \omega_2
   \end{pmatrix}
  \right)^p
  \ -\ 
  \tau_\Rssf
  \left(
   \dfrac{1}{N}
   \,
   (x_1,x_2\,p)
   \begin{pmatrix}
     p\,\omega_1
     \\
     \omega_2
   \end{pmatrix}
   \,,\,
   \begin{pmatrix}
     p\,\omega_1
     \\
     \omega_2
   \end{pmatrix}
  \right)
  ,
\]
and since
the $q$-expansion of
\[
  \tau_\Rssf
  \left(
   \dfrac{1}{N}
   \,
   (x_1,x_2\,p)
   \begin{pmatrix}
     p\,\omega_1
     \\
     \omega_2
   \end{pmatrix}
   \,;\,
   \begin{pmatrix}
     p\,\omega_1
     \\
     \omega_2
   \end{pmatrix}
  \right)
\]
arises from that of
\[
  \tau_\Rssf
  \left(
   \dfrac{1}{N}
   \,
   (x_1,x_2)
   \begin{pmatrix}
     \omega_1
     \\
     \omega_2
   \end{pmatrix}
   \,;\,
   \begin{pmatrix}
     \omega_1
     \\
     \omega_2
   \end{pmatrix}
  \right)
\]
by replacing
$q$ by $q^p$
and
$\zeta_N$ by $\zeta_N^p$,
it follows that
the $p$-th power of
the $q$-expansion of
$\tau_\Rssf
 \left(
  \dfrac{1}{N}
  \,
  (x_1,x_2)
  \begin{pmatrix}
    \omega_1
    \\
    \omega_2
  \end{pmatrix}
  \,;\,
  \begin{pmatrix}
    \omega_1
    \\
    \omega_2
  \end{pmatrix}
 \right)
$
is termwise congruent to it
modulo~$p$,
whence
the $q$-expansion of
$\delta_{\left(\begin{smallmatrix}
        p & 0
        \\
        0 & 1
        \end{smallmatrix}\right)}
 ((x_1,x_2);\omega)
$
indeed
has coefficients
divisible by~$p$.


\section*{B. Number-theoretic basics}
\addcontentsline{toc}{section}{\hfill
B. Number-theoretic basics
\hfill}


\paragraph{8. Orders in quadratic number fields.}
\footnote{
\textit{Dirichlet-Dedekind},
Vorlesungen \"uber Zahlentheorie
(Lectures on number theory),
XI.~Supplement,
cf.~also
this Encyclopedia I~2, 19.
}
Let $\Sigma$ be
a number field;
the principal order of $\Sigma$
will be denoted by
$\Rssf_1$,
while
$\Rssf$
stands for
an arbitrary order of $\Sigma$.
Ideals of $\Rssf$
will be denoted by
$\afrak_\Rssf, \bfrak_\Rssf, \ldots$;
an ideal of $\Rssf$
is called
\emph{integral in $\Rssf$}
if it lies in $\Rssf$.
In particular,
let $\Ffrak_\Rssf$ be
the \emph{conductor of $\Rssf$},
that is to say
the ring-quotient
$\Rssf{:}\Rssf_1$,
or,
what amounts to the same thing,
the largest ideal of $\Rssf$,
integral in $\Rssf$,
which is at the same time
an ideal of $\Rssf_1$.

As usual,
two integral ideals
$\afrak_\Rssf$ and $\bfrak_\Rssf$
of $\Rssf$
are called \emph{coprime}
if
$\afrak_\Rssf + \bfrak_\Rssf
 =
 \Rssf
$.
An arbitrary ideal
$\afrak_\Rssf$ of $\Rssf$
is called
\emph{prime to
the integral ideal
$\bfrak_\Rssf$ of $\Rssf$}
if there exists
an element $\beta$
of $\Sigma$
such that
$\beta\,\Rssf$ and $\beta\,\afrak_\Rssf$
are integral in $\Rssf$
and
are coprime to $\bfrak_\Rssf$.
An ideal
$\afrak_\Rssf$ is called
a \emph{characteristic ideal of $\Rssf$}
or
\emph{from the order $\Rssf$}
if the module-quotient
$\afrak_\Rssf{:}\afrak_\Rssf$
equals $\Rssf$,
that is to say,
if $\Rssf$ is
maximal among
the orders of $\Sigma$
of which
$\afrak_\Rssf$
is an ideal.
An ideal
$\afrak_\Rssf$ is called
\emph{invertible}
if there exists an ideal
$\afrak_\Rssf^{-1}$
such that
$\afrak_\Rssf\,\afrak_\Rssf^{-1}
 =
 \Rssf
$;
then
$\afrak_\Rssf$
is necessarily
a characteristic ideal of $\Rssf$.
The invertible ideals of $\Rssf$
form a group
$\Jfrak_\Rssf$.
The ideals of $\Rssf$
prime to $\Ffrak_\Rssf$
are invertible
(and hence characteristic in $\Rssf$);
they thus form
a subgroup
$\Jfrak_\Rssf^{(0)}$
of
$\Jfrak_\Rssf$.

\textit{
The group
$\Jfrak_\Rssf^{(0)}$
is a canonical isomorphic image of
the group
$\Dfrak_\Rssf$
of divisors of $\Sigma$
prime to $\Ffrak_\Rssf$;
}
we denote 
the ideal
associated to
the divisor $\afrak$
by $\afrak_\Rssf$.
More precisely,
we can describe
this isomorphism
as follows.
Let $S_\Rssf$ be
the ring of all elements
$\dfrac{\alpha}{\beta}$
where $\alpha$
and $\beta \neq 0$
lie in $\Rssf$
and
$\beta\,\Rssf$ is prime to
$\Ffrak_\Rssf$.
Then for any given divisor
$\afrak$,
the associated ideal
$\afrak_\Rssf$
is the collection of
the elements of $S_\Rssf$
divisible by $\afrak$;
conversely
$\afrak$ is
the greatest common divisor of
all elements of
$\afrak_\Rssf$.
For integral $\afrak$,
the index
$[\Rssf{:}\afrak_\Rssf]$
is equal to
the norm
$\Nm\,\afrak$
of the divisor;
in general
one has
\[
  [\,
  \bfrak_\Rssf
  \,{:}\,
  \bfrak_\Rssf\,\afrak_\Rssf
  \,]
  \ =\ 
  \Nm\,\afrak  
\]
for an invertible ideal
$\bfrak_\Rssf$ of $\Rssf$.
Let $\Rssf'$ be
an order
contained in $\Rssf$;
then
one has
$\Ffrak_{\Rssf'}
 \subseteq
 \Ffrak_\Rssf
$.
For divisors $\afrak$
prime to $\Ffrak_{\Rssf'}$
(and hence to $\Ffrak_\Rssf$),
one has
\[
  \afrak_{\Rssf'}
  \ =\ 
  \afrak_\Rssf
  \ \cap\ 
  S_{\Rssf'}
  ,
  \qquad
  \text{thus}
  \quad
  \afrak_{\Rssf'}
  \ =\ 
  \afrak_\Rssf
  \ \cap\ 
  R
  \quad
  \text{for integral $\afrak$}
\]
and
\[
  \afrak_\Rssf
  \ =\ 
  \afrak_{\Rssf'}
  \,
  R
  .
\]

Let $\Hfrak_\Rssf$ be
the group of principal ideals
$\alpha\,\Rssf$
($\alpha \neq 0$ from $\Sigma$)
of $\Rssf$,
and let
$\Rfrak_\Rssf = \Jfrak_\Rssf / \Hfrak_\Rssf$
be the
\emph{ideal class group} of $\Rssf$.
In every ideal class of $\Rssf$
there exists
an ideal prime to $\Ffrak_\Rssf$,
i.e.~one has
$\Jfrak_\Rssf
 =
 \Jfrak_\Rssf^{(0)}
 \,
 \Hfrak_\Rssf
$,
whence setting
$\Jfrak_\Rssf^{(0)}
 \cap
 \Hfrak_\Rssf
 =
 \Hfrak_\Rssf^{(0)}
$,
one sees that
$\Rfrak_\Rssf$
is canonically isomorphic with
$\Jfrak_\Rssf^{(0)} / \Hfrak_\Rssf^{(0)}$.
Here,
$\Hfrak_\Rssf^{(0)}$
is the group of all ideals
$\alpha\,\Rssf$,
$\alpha \in G_\Rssf$,
where $G_\Rssf$ is
the group of all quotients
$\dfrac{\beta}{\gamma}$
of elements
$\beta \neq 0$,
$\gamma \neq 0$,
of $\Rssf$
with
$\beta\,\Rssf,
 \gamma\,\Rssf
$
prime to $\Ffrak_\Rssf$.
The caonical isomorphism
from $\Dfrak_\Rssf$
onto $\Jfrak_\Rssf^{(0)}$
thus induces
a canonical isomorphism
from the factor group of
$\Dfrak_\Rssf$
modulo
the group $(G_\Rssf)$ of
principal divisors $(\alpha)$,
$\alpha \in G_\Rssf$
---
which shall be called
the \emph{divisor class group of $\Rssf$}
---
onto
the ideal class group of $\Rssf$.
The divisor class group of $\Rssf$
shall henceforth also be
denoted by
$\Rfrak_\Rssf$.

For a \emph{quadratic} number field
$\Sigma$,
the following holds:

\textit{
Any characteristic ideal
$\afrak_\Rssf$
of an order $\Rssf$
is invertible.
The conductor
$\Ffrak_\Rssf$
of an order $\Rssf$
is of the form
$\Ffrak_\Rssf
 =
 f\,\Rssf_1
$
with $f$ being
a natural number;
for every natural number $f$,
there exists
exactly one order
with conductor
$f\,\Rssf$
}
(we say briefly:
\emph{with conductor $f$}),
namely,
the ring $R_f$ of
all integral elements of $\Sigma$
which are congruent
modulo~$f$
to rational numbers.
Then
$S_{\Rssf_f}$ is
the ring of
all $f$-integral elements of $\Sigma$
which are congruent
modulo~$f$
to rational numbers,
and
$G_{\Rssf_f}$
consists of
all elements of $\Sigma$
prime to $f$
(as divisor)
which are congruent
modulo~$f$
to rational numbers.

Having
$\Rssf_{f'} \subseteq \Rssf_f$
is equivalent to
having
$f \mid f'$.
If $\afrak$ is
a divisor prime to $f'$,
and
$f \mid f'$,
and if
$\alpha_1,\alpha_2$
is a basis of $\afrak_{R_f}$
while
$\alpha_1',\alpha_2'$
is a basis of $\afrak_{R_{f'}}$,
then
$\begin{pmatrix}
 \alpha_1'
 \\
 \alpha_2'
 \end{pmatrix}
 =
 S
 \begin{pmatrix}
 \alpha_1
 \\
 \alpha_2
 \end{pmatrix}
$
for a \emph{primitive} matrix $S$
of determinant
$s = f'/f$.

To prove that,
we choose $\beta$
so that
$\beta\,\afrak$ is integral
and prime to $f'$;
thus
\[
  \beta\,\afrak_{\Rssf_{f'}}
  \ +\ 
  f'\,\Rssf_1
  \ =\ 
  \Rssf_{f'}
  ,
\]
from which
it follows,
on multiplication by $\Rssf_f$,
that
\[
  \beta\,\afrak_{R_{f'}}
  \ +\ 
  f'\,\Rssf_1
  \ =\ 
  \Rssf_f
  ,
\]
and as
$f'\,\Rssf
 \subseteq
 \Rssf_{f'}
$,
one has
\[
  \beta\,\afrak_{R_{f'}}
  \ +\ 
  \Rssf_{f'}
  \ =\ 
  \Rssf_f
  .
\]
Thus we have
the following relation
between
additive factor groups:
\[
  \Rssf_f
  \,/\,
  \Rssf_{f'}
  \ =\ 
  (\beta\,\afrak_{R_f} + \Rssf_{f'})
  \,/\,
  \Rssf_{f'}
  \ \cong\ 
  \beta\,\afrak_{\Rssf_f}
  \,/\,
  \beta\,\afrak_{\Rssf_f} \cap \Rssf_{f'}
  \ =\ 
  \beta\,\afrak_{\Rssf_f}
  \,/\,
  \beta\,\afrak_{\Rssf_{f'}}
  \ \cong\ 
  \afrak_{\Rssf_f}
  \,/\,
  \afrak_{\Rssf_{f'}}
  ;
\]
but it is easy to see that
$\Rssf_{f'}$ has index
$s = f'/f$ in $\Rssf_f$,
whence
$\begin{pmatrix}
 \alpha_1'
 \\
 \alpha_2'
 \end{pmatrix}
 =
 S
 \begin{pmatrix}
 \alpha_1
 \\
 \alpha_2
 \end{pmatrix}
$
for an integer-entry matrix $S$
of determinant~$s$.
The fact that
$S$ is primitive
arises from
the following fact:

\textit{
If $\alpha_1,\alpha_2$ is
a basis of
a characteristic ideal of
$\Rssf_f$,
and $S$ is
a primitive matrix
of determinant~$s$,
then
$S
 \begin{pmatrix}
 \alpha_1
 \\
 \alpha_2
 \end{pmatrix}
$
is a basis of
a characteristic ideal
$\bfrak_{\Rssf_{f'}}$
of an order
$\Rssf_{f'}$
whose conductor $f'$
satisfies:
$f' \mid fs$
and
$fs^{-1} \mid f'$.
}

As $S$ and $s\,S^{-1}$
are primitive
and of determinant~$s$,
it suffices to prove
that
$f' \mid fs$,
i.e.~that
for any
$\xi \in \Rssf_{fs}$,
one always has
$\bfrak_{\Rssf_{f'}}\,\xi
 \subseteq
 \bfrak_{\Rssf_{f'}}
$,
or that
one has
\[
  S
  \begin{pmatrix}
  \alpha_1
  \\
  \alpha_2
  \end{pmatrix}
  \,
  \xi
  \ =\ 
  D_\xi'
  \,
  S
  \begin{pmatrix}
  \alpha_1
  \\
  \alpha_2
  \end{pmatrix}
\]
for an \emph{integer}-entry matrix
$D_\xi'$.
Now if
$\xi = r + fs\eta$
with
$r$ being
a rational integer
and
$\eta$ being
an integer in $\Sigma$,
then
$f\eta$ is in $\Rssf_f$,
and
$\begin{pmatrix}
 \alpha_1
 \\
 \alpha_2
 \end{pmatrix}
 f
 \eta
 =
 D_{f\eta}
 \begin{pmatrix}
 \alpha_1
 \\
 \alpha_2
 \end{pmatrix}
$
with an integer-entry
$D_{f\eta}$,
thus
$S
 \begin{pmatrix}
 \alpha_1
 \\
 \alpha_2
 \end{pmatrix}
 \xi
 =
 (r\,E + S\,D_{f\eta}\,s\,S^{-1})
 \,
 S
 \begin{pmatrix}
 \alpha_1
 \\
 \alpha_2
 \end{pmatrix}
$,
and
$D_\xi'
 =
 (r\,E + S\,D_{f\eta}\,s\,S^{-1})
$
has integer-entries.

When $p$ is a prime number,
one has the sharper result:

\textit{
If $\alpha_1,\alpha_2$ is
a basis of
a characteristic ideal
$\afrak_{\Rssf_f}$ of $\Rssf_f$,
where $f$ is
divisible by~$p$,
then
for a primitive matrix $P$
of determinant~$p$,
$P
 \begin{pmatrix}
 \alpha_1
 \\
 \alpha_2
 \end{pmatrix}
$
is a basis of
a characteristic ideal
$\bfrak_{\Rssf'}$
in
$\Rssf' = \Rssf_{fp^{-1}}$
or in
$\Rssf' = \Rssf_{fp}$
but not in
$\Rssf_f$.
}

Proof.
One has
$p\,\afrak_{\Rssf_f}
 \subset
 \bfrak_{\Rssf'}
 \subset
 \afrak_{\Rssf_f}
$
with
$\bfrak_{\Rssf'}
 \neq
 \afrak_{\Rssf_f}
 \ \,\text{or}\,\ 
 p\,\afrak_{\Rssf_f}
$;
thus
\[
  p_{\Rssf_f}
  \ \subset\ 
  \bfrak_{\Rssf'}
  \,
  \afrak_{\Rssf_f}^{-1}
  \ \subset\ 
  \Rssf_f
  ,
  \qquad
  \bfrak_{\Rssf'}
  \,
  \afrak_{\Rssf_f}^{-1}
  \ \neq\ 
  \Rssf_f
  \ \ \text{or}\ \ 
  p\,\Rssf_f
  ,
  \qquad
  [\,
  \Rssf_f
  \,{:}\,
  \bfrak_{\Rssf'}\,\afrak_{\Rssf_f}^{-1}
  \,]
  \ =\ 
  [\,
  \bfrak_{\Rssf'}\,\afrak_{\Rssf_f}^{-1}
  \,{:}\,
  p\,\Rssf_f
  \,]
  \ =\ 
  p
\]
because
$[
 \Rssf_f
 \,{:}\,
 p\,\Rssf_f
 ]
 =
 p^2
$.

For the module-quotient
$[
 \Rssf_f
 \,{:}\,
 \Rssf_{fp^{-1}}
 ]
 =
 \qfrak
$
(the conductor of
$\Rssf_f$
with respect to
$\Rssf_{fp^{-1}}$)
one has
correspondingly
\[
  p
  \,
  \Rssf_f
  \ \subset\ 
  \qfrak
  \ \subset\ 
  \Rssf_f
  ,
  \qquad
  [\,
  \Rssf_f
  \,{:}\,
  \qfrak
  \,]
  \ =\ 
  [\,
  \qfrak
  \,{:}\,
  p\,\Rssf_f
  \,]
  \ =\ 
  p
  .
  \tag{25}
  \label{eqn:25}
\]
This is because
one has
\[
  \qfrak
  \ =\ 
  1
  \cdot
  \qfrak
  \ \subseteq\ 
  \Rssf_{fp^{-1}}
  \cdot
  \qfrak
  \ \subseteq\ 
  \Rssf_f
  ,
\]
and
\[
  p
  \,
  \Rssf_f
  \cdot
  \Rssf_{fp^{-1}}
  \ =\ 
  p
  \,
  \Rssf_{fp^{-1}}
  \ \subseteq\ 
  \Rssf_f
  ,
  \qquad
  \text{i.e.}
  \quad
  p
  \,
  \Rssf_f
  \ \subseteq\ 
  \qfrak
  ,
\]
since
for $\xi \in \Rssf_{fp^{-1}}$,
there exists
a rational number $r$
with
$\xi \equiv r \mod f\,p^{-1}$,
so
$p\,\xi \equiv p\,r \mod f$,
which proves that
$p\,\xi \in \Rssf_f$.
Here
$\qfrak$ is
an ideal of $\Rssf_{fp^{-1}}$,
consequently
$\neq
 \Rssf_f
 \ \ \text{or}\ \ 
 p\,\Rssf_f
$,
because
$\Rssf_{fp^{-1}}
 \cdot
 \Rssf_{fp^{-1}}
 \,
 \qfrak
 =
 \Rssf_{fp^{-1}}
 \,
 \qfrak
 \subseteq
 \Rssf_f
$
proves that
$\Rssf_{fp^{-1}}
 \,
 \qfrak
 =
 \qfrak
$.
Now suppose
$\bfrak_{\Rssf'}$
and thus also
$\bfrak_{\Rssf'}\,\afrak_{\Rssf_f}^{-1}$
are characteristic in
$\Rssf_f$;
then
one must have
$\qfrak
 \neq
 \bfrak_{\Rssf'}\,\afrak_{\Rssf_f}^{-1}
$,
and since
$\qfrak$ is prime
by~\eqref{eqn:25},
one has
$\qfrak + \bfrak_{\Rssf'}\,\afrak_{\Rssf_f}^{-1}
 =
 \Rssf_f
$,
and consequently
$\qfrak\,\bfrak_{\Rssf'}\,\afrak_{\Rssf_f}^{-1}
 =
 \qfrak \cap \bfrak_{\Rssf'}\,\afrak_{\Rssf_f}^{-1}
 =
 p\,\Rssf_f
$,
whence
it follows that:
$\qfrak\,\bfrak_{\Rssf'}\,\afrak_{\Rssf_f}^{-1}
 =
 \Rssf_{fp^{-1}}
 \,
 \qfrak_{\Rssf'}
 \,
 \afrak_{\Rssf_f}^{-1}
 =
 p
 \,
 \Rssf_{f\,p^{-1}}
$,
which is to say
$\Rssf_f
 =
 \Rssf_{fp^{-1}}
$,
a contradiction.

\textit{
Let $\alpha_1,\alpha_2$ be
a basis of
a characteristic ideal
$\afrak_{\Rssf_f}$
of $\Rssf$,
and
let $S$ be primitive 
of determinant~$s$
such that
$S
 \begin{pmatrix}
 \alpha_1
 \\
 \alpha_2
 \end{pmatrix}
$
is a basis of
a characteristic ideal
$\bfrak_{\Rssf_{fs}}$
of $\Rssf_{fs}$.
Then
one has
$\bfrak_{\Rssf_{fs}}
 \,
 \Rssf_f
 =
 \afrak_{\Rssf_f}
$.
If a second
primitive matrix $S'$
of determinant~$s$
also makes
$S'
 \begin{pmatrix}
 \alpha_1
 \\
 \alpha_2
 \end{pmatrix}
$
a basis of
a characteristic ideal
$\bfrak'_{\Rssf_{fs}}$
of
$\Rssf_{fs}$,
and if
$\bfrak'_{\Rssf_{fs}}$
is equivalent to
$\bfrak_{\Rssf_{fs}}$,
then
$S'$ is equivalent to
$S\,D_\varepsilon$,
where
$D_\varepsilon$ is
the representing matrix of
a unit $\varepsilon$
of $\Rssf_f$
for the basis
$\alpha_1,\alpha_2$:
one has
$\begin{pmatrix}
 \alpha_1
 \\
 \alpha_2
 \end{pmatrix}
 \varepsilon
 =
 D_\varepsilon
 \begin{pmatrix}
 \alpha_1
 \\
 \alpha_2
 \end{pmatrix}
$.
}

Proof.
Let $\beta$ in $\Sigma$
be chosen
so that
$\beta\,\bfrak_{\Rssf_{fs}}$
is integral in $\Rssf_{fs}$
and prime to
$fs\,\Rssf_1$,
so that
$\beta\,\bfrak_{\Rssf_{fs}}
 +
 fs\,\Rssf_1
 =
 \Rssf_{fs}
$.
Since
it follows
by multiplication with $\Rssf_f$
that
\[
  \beta\,\bfrak_{\Rssf_{fs}}\,\Rssf_f
  \ +\ 
  fs\,\Rssf_1
  \ =\ 
  \Rssf_f
  ,
\]
and since
$\beta\,\bfrak_{\Rssf_{fs}}\,\Rssf_f
 \subseteq
 \beta\,\afrak_{\Rssf_f} \cap \Rssf_f
 \subseteq
 \Rssf_f
$,
one has
\[
  (\,
  \beta\,\afrak_{\Rssf_f}
  \cap
  \Rssf_f
  \,)
  \ +\ 
  fs\,\Rssf_1
  \ =\ 
  \Rssf_f
  ,
\]
i.e.~$\beta\,\afrak_{\Rssf_f} \cap \Rssf_f$
is integral in $\Rssf_f$
and prime to
$fs\,\Rssf_1$.
By the above proof,
the canonically associated ideal
$\beta\,\afrak_{\Rssf_f}
 \cap
 \Rssf_f
 \cap
 \Rssf_{fs}
 =
 \beta\,\afrak_{\Rssf_f}
 \cap
 \Rssf_{fs}
$
of
$\Rssf_{fs}$
thus has index~$s$
in
$\beta\,\afrak_{\Rssf_f}
 \cap
 \Rssf_f
$,
and because
\[
  \beta\,\bfrak_{\Rssf_{fs}}
  \ \subseteq\ 
  \beta\,\afrak_{\Rssf_f} \cap \Rssf_{fs}
  \ \subseteq\ 
  \beta\,\afrak_{\Rssf_f} \cap \Rssf_f
  \ \subseteq\ 
  \beta\,\afrak_{\Rssf_f}
\]
and by assumption
one has
$[
 \beta\,\afrak_{\Rssf_f}
 \,{:}\,
 \beta\,\bfrak_{\Rssf_{fs}}
 ]
 =
 s
$,
it follows that
\[
  \beta\,\afrak_{\Rssf_f} \cap \Rssf_{fs}
  \ =\ 
  \beta\,\bfrak_{\Rssf_{fs}}
  ,
  \qquad
  \beta\,\afrak_{\Rssf_f} \cap \Rssf_f
  \ =\ 
  \beta\,\afrak_{\Rssf_f}
  ,
\]
and hence
one must have
\[
  \beta
  \,
  \bfrak_{\Rssf_{fs}}
  \,
  \Rssf_f
  \ =\ 
  \beta
  \,
  \afrak_{\Rssf_f}
  ,
  \qquad
  \bfrak_{\Rssf_{fs}}
  \,
  \Rssf_f
  \ =\ 
  \afrak_{\Rssf_f}
  ,
\]
as claimed.

The fact that
$\bfrak_{\Rssf_{fs}}'$
is equivalent to
$\bfrak_{\Rssf_{fs}}$
means that
there exists
an element
$\varepsilon$ of $\Sigma$
with
$\bfrak_{\Rssf_{fs}}'
 =
 \bfrak_{\Rssf_{fs}}
 \,
 \varepsilon
$.
Multiplication by $\Rssf_f$
yields
$\afrak_{\Rssf_f}
 =
 \afrak_{\Rssf_f}
 \,
 \varepsilon
$;
thus
$\varepsilon$ is a unit
in $\Rssf_f$.

Then
$S'
 \begin{pmatrix}
 \alpha_1
 \\
 \alpha_2
 \end{pmatrix}
$
and
$S
 \begin{pmatrix}
 \alpha_1
 \\
 \alpha_2
 \end{pmatrix}
 \,
 \varepsilon
 =
 S
 \,
 D_\varepsilon
 \begin{pmatrix}
 \alpha_1
 \\
 \alpha_2
 \end{pmatrix}
$
are bases of
$\bfrak'_{\Rssf_{fs}}
 =
 \bfrak_{\Rssf_{fs}}
 \,
 \varepsilon
$;
consequently
$S'$ and $S\,D_\varepsilon$
are equivalent.

\textit{
Suppose
the conductor $f$
is divisible by~$s$.
If $\alpha_1,\alpha_2$ is
a basis of
a characteristic ideal
$\afrak_{\Rssf_f}$ of $\Rssf_f$,
then
there exists
exactly one class of
primitive matrices $S$
of determinant~$s$
for which
$S
 \begin{pmatrix}
 \alpha_1
 \\
 \alpha_2
 \end{pmatrix}
$
is a basis of
a characteristic ideal of
$\Rssf_{fs^{-1}}$.
When
$\afrak_{\Rssf_f}$
is prime to $f$,
this ideal of $\Rssf_{fs^{-1}}$
is equal to
$s\,\afrak_{\Rssf_{fs^{-1}}}$.
}

Proof.
Since
the result
does not depend on
a factor of
$\afrak_{\Rssf_f}$,
one can assume
$\afrak_{\Rssf_f}$
to be prime to $f$.
A basis of
$\afrak_{\Rssf_f}$
canonically associated to
the ideal
$\afrak_{\Rssf_{fs^{-1}}}$
of $\Rssf_{fs^{-1}}$
is of the form
$S_0^{-1}
 \begin{pmatrix}
 \alpha_1
 \\
 \alpha_2
 \end{pmatrix}
$,
where $S_0$ is primitive,
$|S_0| = s$.
Then
$S = s\,S_0^{-1}$
is primitive
of determinant~$s$
and
$S
 \begin{pmatrix}
 \alpha_1
 \\
 \alpha_2
 \end{pmatrix}
$
is a basis of
the characteristic ideal
$s\,\afrak_{\Rssf_{fs^{-1}}}$
of $\Rssf_{fs^{-1}}$.

Let $S'$ be any
primitive matrix
of determinant~$s$,
for which
$S'
 \begin{pmatrix}
 \alpha_1
 \\
 \alpha_2
 \end{pmatrix}
$
is a basis of
a characteristic ideal
$\bfrak_{\Rssf_{fs^{-1}}}$
of $\Rssf_{fs^{-1}}$.
The primitive matrix
$s\,{S'}^{-1}$
of determinant~$s$
sends
the basis
$S'
 \begin{pmatrix}
 \alpha_1
 \\
 \alpha_2
 \end{pmatrix}
$
of
$\bfrak_{\Rssf_{fs^{-1}}}$
to the basis
$s\,\alpha_1, s\,\alpha_2$
of
$s\,\afrak_{\Rssf_f}$;
thus one has
\[
  s
  \,
  \afrak_{\Rssf_f}
  \,
  \Rssf_{fs^{-1}}
  \ =\ 
  \bfrak_{\Rssf_{fs^{-1}}}
  ,
\]
and for $S'=S$
this becomes
\[
  s
  \,
  \afrak_{\Rssf_f}
  \,
  \Rssf_{fs^{-1}}
  \ =\ 
  s
  \,
  \afrak_{\Rssf_{fs^{-1}}}
  ,
\]
whence
$\bfrak_{\Rssf_{fs^{-1}}}
 =
 s
 \,
 \afrak_{\Rssf_{fs^{-1}}}
$,
so
$S'
 \begin{pmatrix}
 \alpha_1
 \\
 \alpha_2
 \end{pmatrix}
$
and
$S
 \begin{pmatrix}
 \alpha_1
 \\
 \alpha_2
 \end{pmatrix}
$
are bases of
$s\,\afrak_{\Rssf_{fs^{-1}}}$,
and therefore
$S'$ and $S$
are equivalent.


\section*{C. The first main theorem}
\addcontentsline{toc}{section}{\hfill
C. The first main theorem
\hfill}


\paragraph{9. Formulation of the first main theorem.}

Let $\Sigma$ be
an imaginary quadratic number field,
and let
$\Rssf = \Rssf_f$
be the order
of conductor~$f$
in $\Sigma$.
The ideals in
a given ideal class
$\kfrak_\Rssf$
of $\Rssf$
are lattices in
an equivalence class
in the sense of~\S\textbf{1}.
Let
$j(\kfrak_\Rssf)$
denote
the invariant of
this lattice class;
thus
$j(\kfrak_\Rssf)
 =
 j
 \left(
 \dfrac
 {\alpha_1}
 {\alpha_2}
 \right)
$
where
$\alpha_1,\alpha_2$
is a basis of
an ideal
$\afrak_\Rssf$
from
$\kfrak_\Rssf$
with
$\Im
 \left(
 \dfrac
 {\alpha_1}
 {\alpha_2}
 \right)
 >
 0
$.
This
\emph{invariant
of the ideal class
$\kfrak_\Rssf$}
determines
$\kfrak_\Rssf$
uniquely,
because
if $\kfrak'_{\Rssf'}$
is an ideal class
in a field $\Sigma'$
and
$j(\kfrak'_{\Rssf'})
 =
 j(\kfrak_\Rssf)
$,
then
with
$\alpha'_1,\alpha'_2$
being a basis of
an ideal
$\afrak'_{\Rssf'}$
from $\kfrak'_{\Rssf'}$,
we have
$\begin{pmatrix}
 \alpha'_1
 \\
 \alpha'_2
 \end{pmatrix}
 =
 M
 \begin{pmatrix}
 \alpha_1
 \\
 \alpha_2
 \end{pmatrix}
 \,
 \beta
$
with
$\beta \neq 0$,
$M \in \Mfrak$,
so
$\left(
 \dfrac{\alpha'_1}{\alpha'_2}
 \right)
 =
 M
 \left(
 \dfrac{\alpha_1}{\alpha_2}
 \right)
$,
therefore
$\Sigma' = \Sigma$,
$\beta \in \Sigma$,
and
$\afrak'_{\Rssf'}
 =
 \afrak_\Rssf
 \,
 \beta
$
lies in
$\kfrak_\Rssf$,
whence
$\kfrak'_{\Rssf'}
 =
 \kfrak_\Rssf
$.

Let $h_\Rssf$ be
the
\emph{ideal class number}
of $\Rssf$,
i.e.~the order of
the divisor class group
$\Dfrak_\Rssf / (G_\Rssf)$
of $\Rssf$,
and let
$\kfrak_\Rssf^{(1)},
 \ldots,
 \kfrak_\Rssf^{(h_\Rssf)}
$
be the ideal classes of $\Rssf$.

\textit{First main theorem.}
\ 
\textit{
The class invariants
$j(\kfrak_\Rssf^{(1)}),
 \ldots,
 j(\kfrak_\Rssf^{(h_\Rssf)})
$
of $\Rssf$
are algebraic numbers.
The number field
$\Omega_\Rssf
 =
 \Sigma(j(\kfrak_\Rssf))
$
generated over $\Sigma$
by one of the invariants
$j(\kfrak_\Rssf)$
is the class field of $\Sigma$
for the ring divisor class group
$\Dfrak_\Rssf / (G_\Rssf)
 =
 \Rfrak_\Rssf
$.
The field
$\Omega_\Rssf$
is called
the ring class field of $\Sigma$
for the order $\Rssf$.
}

We give
two proofs
of this theorem.
The first
uses the general theory of
abelian number fields
and is accordingly
short.
The central idea
of the second,
after
the class invariants
are identified as
algebraic numbers,
is to determine
the structure of the field
$\Omega_\Rssf$
by means of
the theory of modular functions,
i.e.~in essence,
to see that
the Galois group of
$\Omega_\Rssf/\Sigma$
is isomorphic to
the divisor class group of $\Rssf$
and to prove
the reciprocity law
for $\Omega_\Rssf/\Sigma$,
and in the process
to use only
the elementary parts of
the theory of algebraic numbers.

We can view
the statement of
the first main theorem
in a different way.

If $f(\omega)$ denotes
a modular function
associated to
the full modular group
or\
to one of its subgroups
of finite index,
we call
the complex numbers
$f(\xi)$
for imaginary quadratic argument $\xi$
(with positive imaginary part)
the \emph{singular values
of the function $f(\omega)$}.
The first main theorem
identifies
the number field
generated by
the singular values of
the function $j(\omega)$
(the determination of
$\Pssf(j(\xi))$
inside
$\Sigma(j(\xi))$
is done in~\S\textbf{16}),
because
the lattice
with basis
$\xi,1$
is an ideal in
an order of
the quadratic number field
generated from $\Sigma$.


\paragraph{10. Proof of the first main theorem by means of 
the general theory of abelian number fields.}

It suffices to prove
the following:

\textit{
The complex number
$j(\kfrak_\Rssf)$
is an algebraic number,
and
the degree-one prime divisors
$\pfrak$ of $\Sigma$
which split completely
in $\Sigma(j(\kfrak_\Rssf))$
are,
apart from
finitely many exceptions,
exactly those
for which
$\pfrak_\Rssf$
is a principal ideal.
}

For that,
we need from
the theory of modular functions
only the fact that
for a prime number degree~$p$,
the transformation polynomial
$J_p(X,j(\omega))$
of $j(\omega)$
has rational integer coefficients
and
satisfies the congruence
(cf.~\S\textbf{3})
\[
  J_p(X,j(\omega))
  \ \equiv\ 
  (X - j(\omega)^p)
  \,
  (X^p - j(\omega))
  \mod
  p
  .
  \label{eqn:26}
  \tag{26}
\]

Let $\pfrak$ be
a degree-one
prime divisor of $\Sigma$
with
$\Nm_{\Sigma/\Pssf}\pfrak = p$.
Then
one has
\[
  J_p(j(\kfrak_\Rssf),j(\kfrak_\Rssf \pfrak_\Rssf^{-1}))
  \ =\ 
  0
  ,
  \tag{27}
\]
because
if
$\alpha_1,\alpha_2$
is a basis of
an ideal $\afrak_\Rssf$
of $\kfrak_\Rssf \pfrak_\Rssf^{-1}$,
so that
$j(\kfrak_\Rssf \pfrak_\Rssf^{-1})
 =
 j
 \left(
 \dfrac
 {\alpha_1}
 {\alpha_2}
 \right)
$,
then
there exists
an integer-entry matrix $P$
of determinant~$p$
such that
\[
  P
  \begin{pmatrix}
  \alpha_1
  \\
  \alpha_2
  \end{pmatrix}
\]
is a basis of
$\afrak_\Rssf\,\pfrak_\Rssf
 \in
 \kfrak_\Rssf
$,
whence
$j(\kfrak_\Rssf)
 =
 j
 \left(
 P
 \left(
 \dfrac
 {\alpha_1}
 {\alpha_2}
 \right)
 \right)
$.

From~\eqref{eqn:26}
it follows that
$J_p(X,X)
 \equiv
 -(X^p - X)^2
 \mod
 p
$,
whence
\[
  J_p(X,X)
  \ \neq\ 
  0
  .
  \tag{28}
\]
Let us choose $\pfrak$
so that
$\pfrak_\Rssf$ is
a principal ideal;
then it follows from~\eqref{eqn:26}
that
\[
  J_p(j(\kfrak_\Rssf),j(\kfrak_\Rssf))
  \ =\ 
  0
  ,
\]
which shows that
$j(\kfrak_\Rssf)$
is an algebraic number.

If we exclude
the finitely many $\pfrak$
which are not
prime to the denominators
of all $j(\kfrak_\Rssf)$,
it follows from~\eqref{eqn:26}
that
\[
  \bigl(\,
  j(\kfrak_\Rssf)
  -
  j(\kfrak_\Rssf \pfrak_\Rssf^{-1})^p
  \,\bigr)
  \,
  \bigl(\,
  j(\kfrak_\Rssf)^p
  -
  j(\kfrak_\Rssf \pfrak_\Rssf^{-1})
  \,\bigr)
  \ \equiv\ 
  0
  \mod
  p
  \tag{29}
\]
and thus
for each prime divisor
$\Pfrak$ of $\pfrak$
in $\Sigma(j(\kfrak_\Rssf))$,
\[
  j(\kfrak_\Rssf)^p
  \ \equiv\ 
  j(\kfrak_\Rssf \pfrak_\Rssf^{-1})
  \qquad
  \text{or}
  \qquad
  j(\kfrak_\Rssf)
  \ \equiv\ 
  j(\kfrak_\Rssf \pfrak_\Rssf^{-1})^p
  \mod
  \Pfrak
  .
  \tag{30}
  \label{eqn:30}
\]
As $\pfrak_\Rssf$ is
a principal ideal,
this simplifies to
\[
  j(\kfrak_\Rssf)^p
  \ \equiv\ 
  j(\kfrak_\Rssf)
  \mod
  \Pfrak
  ;
\]
thus
if we now
further suppose that
$\pfrak$
is prime to
the discriminant of
$j(\kfrak_\Rssf)$
over $\Sigma$,
it follows that
in fact
\[
  \alpha^p
  \ \equiv\ 
  \alpha
  \mod
  \Pfrak
\]
for all $\pfrak$-integers
$\alpha$
of $\Sigma(j(\kfrak_\Rssf))$;
therefore
$\pfrak$
splits completely
in $\Sigma(j(\kfrak_\Rssf))$.

Conversely,
suppose
$\pfrak$
splits completely
in $\Sigma(j(\kfrak_\Rssf))$.
Then
one has
\[
  j(\kfrak_\Rssf)^p
  \ \equiv\ 
  j(\kfrak_\Rssf)
  \mod
  \Pfrak
  ,
\]
whence
by~\eqref{eqn:30},
it follows that
\[
  j(\kfrak_\Rssf)
  \ \equiv\ 
  j(\kfrak_\Rssf \pfrak_\Rssf^{-1})
  \mod
  \Pfrak
  .
\]
We now exclude
the finitely many $\pfrak$
which are not prime to
all the differences
$j(\kfrak_\Rssf)
 -
 j(\kfrak_\Rssf')
$,
$\kfrak_\Rssf
 \neq
 \kfrak_\Rssf'
$,
and it follows that
$j(\kfrak_\Rssf)
 =
 j(\kfrak_\Rssf \pfrak_\Rssf^{-1})
$;
thus
$\pfrak_\Rssf$
is a principal ideal.


\paragraph{11. Overview of
the proof of the first main theorem
not relying on
general class field theory.}

Deviating from~\S\textbf{9},
we define
the field $\Omega_\Rssf$
by
\[
  \Omega_\Rssf
  \ =\ 
  \Sigma
  \bigl(\,
  j(\kfrak_\Rssf^{(1)}),
  \ldots,
  j(\kfrak_\Rssf^{(h_\Rssf)})
  \,\bigr)
  .
\]
The proof
splits into
several steps.

First,
it will be shown that
the class invariants
$j(\kfrak_\Rssf)$
are integral algebraic numbers,
and
to that end
the proof
carried out in~\S\textbf{3}
for the algebraicity of
$j(\kfrak_\Rssf)$
will be expanded;
$j(\kfrak_\Rssf)$
will be seen as
a zero of
a suitable
$J_s(X,X)$.

Second,
a closer examination on
the multiplicity of
$j(\kfrak_\Rssf)$
as a zero of
$J_s(X,X)$
shows:
all (algebraic) numbers
which are
conjugate over $\Pssf$ to
one of the class invariants
$j(\kfrak_\Rssf)$
of the order $\Rssf$
are in fact
class invariants
$j(\kfrak_\Rssf')$
of the same order $\Rssf$;
in other words:
the class polynomial
of the order $\Rssf$
\[
  H_\Rssf(X)
  \ =\ 
  \prod_{\nu=1}^{h_\Rssf}
  (X - j(\kfrak_\Rssf^{(\nu)})
  \tag{31}
\]
has (integral) rational coefficients.
From this
it follows that
$\Omega_\Rssf$
is normal
over $\Pssf$
(all the more so
over $\Sigma$).

Third,
an isomorphism
$\kfrak_\Rssf
 \to
 \sigma(\kfrak_\Rssf)
$
from the class group
$\Rfrak_\Rssf$ of $\Rssf$
\emph{into}
the Galois group
$\Gfrak^{\Omega_\Rssf}_{\Sigma}$
of $\Omega_\Rssf$ over $\Sigma$
will be established.

This isomorphism
$\sigma$
has the property that
\[
  j(\hfrak_\Rssf)
  \,
  \sigma(\kfrak_\Rssf)
  \ =\ 
  j(\hfrak_\Rssf\,\kfrak_\Rssf^{-1})
  ;
  \tag{32}
\]
the $j(\hfrak_\Rssf)$
are thus
a full system of (algebraic) numbers
conjugate over $\Sigma$:
the class polynomial
$H_\Rssf(X)$
is irreducible over $\Sigma$.
The basis for that,
aside from
the theory of
decomposition group,
in particular
the \textit{Frobenius}-congruence,
is a sharpening of
the statement~\eqref{eqn:30},
namely
the congruence
\[
  j(\kfrak_\Rssf)^p
  \ \equiv\ 
  j(\kfrak_\Rssf \pfrak_\Rssf^{-1})
  \mod
  \Pfrak
  \tag{33}
  \label{eqn:33}
\]
for any
degree-one prime divisor
$\pfrak$
of norm~$p$ over $\Sigma$
and its prime factor
$\Pfrak$
in $\Omega_\Rssf$;
for its proof
the function
$\varphi_S(\omega)$
is needed.

Fourth,
it will be shown that
\textit{
$\sigma(\Rfrak_\Rssf)$
is equal to
$\Gfrak^{\Omega_\Rssf}_{\Sigma}$
}
(and not only
to a subgroup of
$\Gfrak^{\Omega_\Rssf}_{\Sigma}$).
In other words:

One has
$\Omega_\Rssf
 =
 \Sigma(j(\kfrak_\Rssf))
$
for every class
$\kfrak_\Rssf$
of the order $\Rssf$
(whence
the deviating definition
of $\Omega_\Rssf$
is justified),
or also:
$\Sigma(j(\kfrak_\Rssf))$
is already
normal over $\Sigma$.
For that
we have to compute
the \textit{Galois} group
$\Gfrak^{\Omega_\Rssf}_{\Pssf}$
of $\Omega_\Rssf$
over $\Pssf$,
which is done
again
by means of~\eqref{eqn:33}.

Fifth:
$\sigma$
is the canonical isomorphism
between
the class group
$\Rfrak_\Rssf$
and the Galois group
$\Gfrak^{\Omega_\Rssf}_{\Sigma}$
asserted by
the \textit{Artin} reciprocity law.
Namely,
it will be shown that
the decomposition assertion
of the reciprocity law
holds:

$\sigma(\efrak_\Rssf\,\pfrak_\Rssf)$
is the \textit{Frobenius} automorphism
of the prime divisor $\pfrak$
of $\Sigma$ in $\Omega_\Rssf$.
For the degree-one prime divisors
(except finitely many)
this will follow directly
from the definition of $\sigma$;
for the degree-two prime divisors
it follows from
the structure of
$\Gfrak^{\Omega_\Rssf}_{\Pssf}$.


\paragraph{12. The class polynomial.}

By the 
\emph{class polynomial
$H_\Rssf(X)$
of the order $\Rssf$},
we mean
the polynomial
with the roots
$j(\kfrak_\Rssf^{(\nu)})$;
thus
\[
  H_\Rssf(X)
  \ =\ 
  \prod_{\nu=1}^{h_\Rssf}
  (X - j(\kfrak_\Rssf^{(\nu)}))
  .
\]
The invariant
$j(\kfrak_\Rssf)$
can be written in the form
\[
  j(\kfrak_\Rssf)
  \ =\ 
  j(\alpha)
  ,
  \qquad
  \alpha
  \ =\ 
  \dfrac{\alpha_1}{\alpha_2}
  ,
\]
where
$\alpha_1,\alpha_2$
is a basis of
an ideal from $\kfrak_\Rssf$.
Other expressions for
$j(\kfrak_\Rssf)$
follow from that.
For any $\mu \neq 0$ in $\Rssf$,
one has
\[
  \begin{pmatrix}
  \alpha_1
  \\
  \alpha_2
  \end{pmatrix}
  \mu
  \ =\ 
  D_\mu
  \begin{pmatrix}
  \alpha_1
  \\
  \alpha_2
  \end{pmatrix}
\]
for some
rational-integer-entry
representation matrix
$D_\mu$,
and hence
\[
  j(\kfrak_\Rssf)
  \ =\ 
  j(D_\mu(\alpha)).
\]
The greatest common divisor
of the entries of $D_\mu$
is at the same time
the largest rational part
of $\mu$ in $\Rssf$,
which is
the largest rational integer $t$
for which
$\mu\,t^{-1}$
still lies in $\Rssf$.
Hence
$D_\mu$
is primitive
exactly for
those $\mu$
whose largest rational part
in $\Rssf$
is equal to~$1$;
these $\mu$
are called
\emph{primitive} elements
of $\Rssf$.
Since
$\Nm\,\mu
 =
 |\,D_\mu\,|
$,
one has
for primitive $\mu$
\[
  J_{\Nm\mu}
  (
  j(\kfrak_\Rssf)
  ,
  j(\kfrak_\Rssf)
  )
  \ =\ 
  0
  ,
  \tag{34}
  \label{eqn:34}
\]
or also
\[
  H_\Rssf(X)
  \quad
  \text{is a divisor of}
  \quad
  J_{\Nm\mu}(X,X)
  .
  \tag{35}
\]

Now
in $J_s(X,X)$
for non-square $s$,
the highest power of $X$
in $J_s(X,X)$
has coefficient~$\pm1$:
this coefficient is also
the leading coefficient
of the $q$-expansion
of the function
$J_s(j(\omega),j(\omega))$.
By means of
the system of representatives
$\begin{pmatrix}
 a & b \\
 0 & d
 \end{pmatrix}
$
of the equivalence classes
of level~$s$,
it amounts to
\[
  \begin{aligned}
  J_s(j(\omega),j(\omega))
  &
  \ =\ 
  \prod
  \ 
  \left(
    j
    \Bigl(
    \dfrac{a\omega{+}b}{d}
    \Bigr)
    -
    j(\omega)
  \right)
  \\
  &
  \ =\ 
  \prod
  \ 
  \left(\ 
  \zeta_s^{-ab}
  \,
  q^{\textstyle -\frac{a}{d}}
  \ +\ 
  \cdots
  \ -\ 
  q^{-1}
  \ +\ 
  \cdots
  \ \right)
  .
  \end{aligned}
\]
The leading coefficient
of the series
in the brackets
is either
$\zeta_s^{-ab}$
or~$-1$,
according to
whether
$a>d$
or
$a<d$;
the case
$a=d$
is not possible,
because
$s = ad$
should not be
a square.
The leading coefficient of
$J_s(j(\omega),j(\omega))$
is consequently
a root of unity,
which
as a rational number
must be equal to~$\pm1$.
Thus~\eqref{eqn:34}
proves our claim
when $\Nm\,\mu$ is not
a square.
A primitive element $\mu$
with non-square norm $\Nm\,\mu$
in the order
$\Rssf = \Rssf_f(\sqrt{-d})$
of conductor~$f$
in the field
$\Sigma = \Pssf(\sqrt{-d})$
with $d$ square-free
is
\[
  \begin{aligned}
  &
  \mu
  \ =\ 
  f
  \,
  \sqrt{-d}
  ,
  &
  &
  \Nm\,\mu
  \ =\ 
  d
  \,
  f^2
  &
  &
  \text{for $d>1$}
  ,
  \\
  &
  \mu
  \ =\ 
  f
  \,
  (1{+}\sqrt{-d})
  ,
  &
  &
  \Nm\,\mu
  \ =\ 
  2
  \,
  f^2
  &
  &
  \text{for $d=1$}
  .
  \end{aligned}
\]

\textit{
$H_\Rssf(X)$
has rational integer coefficients.
}
\\
It suffices to see that
the coefficients of
$H_\Rssf(X)$
are rational numbers;
that they are integers
follows from the fact that
the roots of $H_\Rssf(X)$
are integral.

For a given $s$,
we first ask for
all the zeros of
$J_s(X,X)$
and their multiplicities;
as
$J_1(X,X) =0$,
let us henceforth
assume that
$s>1$.

Obviously
the (complex) number
$\xi = j(\varrho)$
is a zero of
$J_s(X,X)$,
or
$\varrho$
is a zero of
the function
\[
  J_s(j(\omega),j(\omega))
  \ =\ 
  \prod_{\nu=1}^{\psi(s)}
  (\,
  j(S_\nu(\omega)) - j(\omega)
  \,)
  ,
\]
if and only if
there exists
a primitive matrix $S$
of determinant $s$
with
\[
  S(\varrho)
  \ =\ 
  \varrho
  .
\]
This results in
a quadratic equation
for $\varrho$
with rational coefficients
which are not all~$0$
since $S \neq E$;
thus
$\varrho$
is an imaginary quadratic (algebraic) number:
every zero of
$J_s(X,X)$
is a singular class invariant.

The multiplicity of
a zero
$\xi = j(\varrho)$
of $J_s(X,X)$
is equal to
the number~$t$ of solutions
\[
  S^{(1)}
  \,,\,
  \ldots
  \,,\,
  S^{(t)}
\]
of the equation
$S(\varrho) = \varrho$
which are
inequivalent modulo~$\Mfrak$.
That will be proved
if it is shown that
\[
  J_s(X,X)
  \,
  (X-\xi)^{-t}
  \qquad
  \text{is regular but~$\neq 0$
    at $X = \xi$}
  ,
\]
or that
\[
  J_s(j(\omega),j(\omega))
  \,
  (j(\omega)-j(\varrho))^{-t}
  \qquad
  \text{is regular but~$\neq 0$
    at $\omega = \varrho$}
  .
\]
Let us write
this function
in the form
\[
  J_s(j(\omega),j(\omega))
  \,
  (j(\omega)-j(\varrho))^{-t}
  \ =\ 
  \prod_{\mu=1}^{t}
  \ 
  \dfrac
  {j(S^{(\mu)}(\omega)) - j(\omega)}
  {j(\omega)-j(\varrho)}
  \qquad
  \prod_{\mathclap{\text{
      $S_\nu$ inequiv. to all $S^{(\mu)}$
      }}}
  \ 
  \bigl(\,
  j(S_\nu(\omega))
  -
  j(\omega)
  \,\bigr)
  ,
\]
from which
we see that
it suffices to prove:
\[
  j(\omega)
  -
  j(\varrho)
  \qquad
  \text{and}
  \qquad
  j(S^{(\mu)}(\omega))
  -
  j(\omega)
\]
have the same order
at the zero
$\omega = \varrho$.
Now one has
\[
  \begin{gathered}
  j(\omega)
  -
  j(\varrho)
  \ =\ 
  \sum_{n=1}^{\infty}
  \dfrac{j^{(n)}(\varrho)}{n!}
  \,
  (\omega - \varrho)^n
  ,
  \\
  j(S^{(\mu)}(\omega))
  -
  j(\omega)
  \ =\ 
  \sum_{n=1}^{\infty}
  \dfrac{j^{(n)}(\omega)}{n!}
  \,
  (S^{(\mu)}(\omega) - \omega)^n
  \\
  \ =\ 
  \sum_{n=1}^{\infty}
  \dfrac
  {j^{(n)}(\varrho) + j^{(n+1)}(\varrho)(\omega-\varrho) + \cdots}
  {n!}
  \,
  (S^{(\mu)}(\omega) - \omega)^n
  ,
  \end{gathered}  
\]
and these expansions
show that
the orders of
the two functions
at $\omega = \varrho$
differ by
the order of
the zero $\varrho$ of
$S^{(\mu)}(\omega) - \omega$
as a factor.
But the latter
is equal to~$1$,
because
$S^{(\mu)}(\omega) - \omega$
has the two \emph{distinct} zeros
$\varrho$ and $\overline{\varrho}$.

Now let us ask
conversely:
how many times
is a \emph{given}
class invariant
$j(\kfrak_\Rssf)$
a zero of
$J_s(X,X)$?
Again let
$\alpha_1,\alpha_2$
be a basis of
an ideal from $\kfrak_\Rssf$;
thus
$j(\kfrak_\Rssf)
 =
 j
 \left(
 \dfrac
 {\alpha_1}
 {\alpha_2}
 \right)
$.
The desired multiplicity
is the number of
inequivalent primitive $S$
of determinant~$s$
with
$S
 \left(
 \dfrac
 {\alpha_1}
 {\alpha_2}
 \right)
 =
 \dfrac{\alpha_1}{\alpha_2}
$,
or
\[
  S
  \begin{pmatrix}
  \alpha_1
  \\
  \alpha_2
  \end{pmatrix}
  \ =\ 
  \begin{pmatrix}
  \alpha_1
  \\
  \alpha_2
  \end{pmatrix}
  \,
  \mu
\]
for a certain element
$\mu$.
Obviously
$\mu$ lies in $\Rssf$,
and
$S$ is equal to
the representation matrix
$D_\mu$
determined by
the basis
$\alpha_1,\alpha_2$.
Let
$S' = D_{\mu'}$
be a second primitive matrix
of determinant~$s$
with
$S'
 \left(
 \dfrac
 {\alpha_1}
 {\alpha_2}
 \right)
 =
 \dfrac{\alpha_1}{\alpha_2}
$.
Then one has
\[
  S^{-1}\,S'
  \ =\ 
  D_{\mu}^{-1}\,D_{\mu'}
  \ =\ 
  D_{\mu^{-1}\mu'}
  ,
\]
whence
$S$ and $S'$
are equivalent,
i.e.~$D_{\mu^{-1}\mu'}$
is unimodular
if and only if
$\mu^{-1}\mu'$
is a unit in $\Rssf$:

The multiplicity of
the zero $j(\kfrak_\Rssf)$
of $J_s(X,X)$
is equal to
the number of
primitive elements
of norm~$s$
in $\Rssf$
which are
non-associated in $\Rssf$.

Now let $\Rssf$
be given.
We want to determine:
for which $s$
does there exist
only one single class
associated to
a primitive $\mu_0$ in $\Rssf$
with $\Nm\,\mu_0 = s$?

Now
$s = 1$
always has
this property,
for then
$\mu_0$
must be a unit in $\Rssf$.
Henceforth
let us assume
$s > 1$.
As with $\mu_0$,
$\overline{\mu_0}$
is also primitive
of norm~$s$;
thus
the desired $\mu_0$
must satisfy
$\overline{\mu_0}
 = 
 \varepsilon
 \,
 \mu_0
$
for a unit
$\varepsilon$ in $\Rssf$.
As $\mu_0$
cannot be rational
(or else
$\mu_0$
would have
the rational part
$\mu_0 = s^{\frac{1}{2}} > 1$
in $\Rssf$),
one has
$\varepsilon \neq 1$.
Now let $\Rssf$ be
the order
$\Rssf
 =
 \Rssf_f(\sqrt{-d})
$
of conductor~$f$
in the field
$\Sigma
 =
 \Pssf(\sqrt{-d})
$.
A basis of $\Rssf$ is
\[
  \begin{aligned}
  &
  1
  \,,\,
  f\,\sqrt{-d}
  &
  \qquad
  &
  \text{for}
  \quad
  -d \not\equiv 1
  \bmod
  4
  \\
  &
  1
  \,,\,
  f\,\dfrac{1{+}\sqrt{-d}}{2}
  &
  \qquad
  &
  \text{for}
  \quad
  -d \equiv 1
  \bmod
  4
  .
  \end{aligned}
\]
Let us leave aside
the two cases
$\Rssf = \Rssf_1(\sqrt{-1})$
and
$\Rssf = \Rssf_1(\sqrt{-3})$
for now.
Then one has
\[
  \overline{\mu_0}
  \ =\ 
  -
  \mu_0
  ,
\]
thus
\[
  \mu_0
  \ =\ 
  b
  \,
  f
  \,
  \sqrt{-d}
\]
with a rational~$b$.
In the case of
$-d \not\equiv 1
 \bmod
 4
$,
$b$ must be an integer,
and since
$\mu_0$ is primitive,
one has
$b = \pm1$,
hence
\[
  \mu_0
  \ =\ 
  \pm
  f
  \,
  \sqrt{-d}
  ,
  \qquad
  \Nm\,\mu_0
  \ =\ 
  d
  \,
  f^2
  .
\]
In the case of
$-d \equiv 1
 \bmod
 4
$,
\[
  \mu_0
  \ =\ 
  -
  b
  f
  \ +\ 
  2b
  \cdot
  f
  \,
  \dfrac
  {1{+}\sqrt{-d}}{2}
\]
so
$2b$ and $b f$
must be integers,
and since
$\mu_0$ is primitive
in $\Rssf$,
it follows that
\[
  \begin{aligned}
  &
  \mu_0
  \ =\ 
  \pm
  \phantom{\tfrac{1}{2}}
  \,
  f
  \,
  \sqrt{-d}
  ,
  &
  \qquad
  &
  \Nm\,\mu_0
  \ =\ 
  d
  \,
  f^2
  &
  \qquad
  &
  \text{for odd $f$}
  ,
  \\
  &
  \mu_0
  \ =\ 
  \pm
  \tfrac{1}{2}
  \,
  f
  \,
  \sqrt{-d}
  ,
  &
  \qquad
  &
  \Nm\,\mu_0
  \ =\ 
  \tfrac{1}{4}
  d
  \,
  f^2
  &
  \qquad
  &
  \text{for even $f$}
  .
  \end{aligned}
\]

In the case of
$\Rssf = \Rssf_1(\sqrt{-1})$
with the basis
$1,\sqrt{-1}$,
we have either
\[
  \overline{\mu_0}
  \ =\ 
  -
  \mu_0
  ,
\]
which gives
$\mu_0
 =
 \pm\sqrt{-1}
$,
whence
the excluded case of
$\Nm\,\mu_0 = 1$,
or
\[
  \begin{gathered}
  \overline{\mu_0}
  \ =\ 
  \pm
  \sqrt{-1}
  \,
  \mu_0
  ,
  \\
  \mu_0
  \ =\ 
  \pm
  1
  \pm
  \sqrt{-1}
  ,
  \qquad
  \Nm\,\mu_0
  \ =\ 
  2
  .
  \end{gathered}
\]
The four elements
$\mu_0$
with
$\Nm\,\mu_0 = 2$
are associated in $\Rssf$.

In the case of
$\Rssf
 =
 \Rssf_1(\sqrt{-3})
$
with the basis
${\displaystyle
 1
 \,,\,
 \zeta_3
 = 
 \tfrac{1}{2}
 (-1{+}\sqrt{-3})
}$,
one has either
\[
  \overline{\mu_0}
  \ =\ 
  -
  \mu_0
  ,
\]
thus
\[
  \mu_0
  \ =\ 
  \pm
  \sqrt{-3}
  ,
  \qquad
  \Nm\,\mu_0
  \ =\ 
  3
  ,
\]
or
\[
  \overline{\mu_0}
  \ =\ 
  \pm
  \zeta_3^\nu
  \,
  \mu_0
  ,
  \qquad
  \nu
  \ =\ 
  1,2
  ;
\]
and
from these
it follows that
either
\[
  \mu_0
  \ =\ 
  \pm
  \zeta_3^\nu
  ,
\]
which is the excluded case of
$\Nm\,\mu_0 = 1$,
or
\[
  \mu_0
  \ =\ 
  \pm
  \zeta_3^\nu
  \sqrt{-3}
  ,
  \qquad
  \Nm\,\mu_0
  \ =\ 
  3
  .
\]
The six elements
with norm~$3$
are associated in $\Rssf$.

Summarizing,
we can say:

To each $\Rssf$,
there exists
exactly one value
$s = s_\Rssf > 1$
with the property that
there exists
exactly one class
in $\Rssf$
which is
associated in $\Rssf$
to primitive elements
$\mu$
of norm
$\Nm\,\mu = s_\Rssf$,
namely
\[
  \begin{aligned}
  &
  s_\Rssf
  \ =\ 
  d
  \,
  f^2
  &
  \qquad
  &
  \text{for}
  \quad
  \Rssf
  \ =\ 
  \Rssf_f(\sqrt{-d})
  ,
  &
  \quad
  &
  -d \not\equiv 1
  \bmod
  4
  ,
  \\
  &
  &
  &
  &
  &
  \qquad
  \text{excluding}
  \quad
  \Rssf
  =
  \Rssf_1(\sqrt{-1})
  ,
  \\
  &
  s_\Rssf
  \ =\ 
  2
  &
  &
  \text{for}
  \quad
  \Rssf
  \ =\ 
  \Rssf_1(\sqrt{-1})
  ,
  &
  &
  \\
  &
  s_\Rssf
  \ =\ 
  d
  \,
  f^2
  &
  &
  \text{for}
  \quad
  \Rssf
  \ =\ 
  \Rssf_f(\sqrt{-d})
  ,
  &
  &
  -d \equiv 1
  \bmod
  4
  ,
  \quad
  \text{$f$ odd}
  ,
  \\
  &
  s_\Rssf
  \ =\ 
  \tfrac{1}{4}
  d
  \,
  f^2
  &
  &
  \text{for}
  \quad
  \Rssf
  \ =\ 
  \Rssf_f(\sqrt{-d})
  ,
  &
  &
  -d \equiv 1
  \bmod
  4
  ,
  \quad
  \text{$f$ even}
  .
  \end{aligned}
\]

What possibilities are there
that
for two distinct orders
$\Rssf$ and $\Rssf'$,
the elements
$s_\Rssf$ and $s_{\Rssf'}$
are equal?

\noindent
When the two orders
$\Rssf = \Rssf_f(\sqrt{-d})$
and
$\Rssf' = \Rssf_{f'}(\sqrt{-d})$
belong to
the same field
$\Pssf(\sqrt{-d})$,
with
$f < f'$,
the above table shows:

One always has
$s_\Rssf < s_{\Rssf'}$
excluding the case
$-d \equiv 1
 \bmod
 4
$,
$f$~odd
and
$f' = 2\,f$,
in which case
$s_\Rssf
 =
 s_{\Rssf'}
 =
 d
 \,
 f^2
$.

As $s_\Rssf$
differs only by
a square factor
from
the discriminant of the field
to which
the order $\Rssf$ belongs,
with the one exception
$\Rssf = \Rssf_1(\sqrt{-1})$,
$s_\Rssf = 2$,
there exists
only one possibility
for two orders
$\Rssf, \Rssf'$
in different fields
to have
$s_\Rssf = s_{\Rssf'}$,
namely
\[
  \Rssf
  \ =\ 
  \Rssf_1(\sqrt{-1})
  \qquad
  \text{and}
  \qquad
  \Rssf'
  \ =\ 
  \Rssf_1(\sqrt{-2})
  \qquad
  \text{with}
  \qquad
  s_\Rssf = s_{\Rssf'} = 2
  .
\]
Hence
there does not at all exist
three distinct orders
with the same
$s_\Rssf$.

Next
let $\Rssf$ be
an order
for which
there does not exist
an $\Rssf'$
different from it
with
$s_\Rssf = s_{\Rssf'}$.
Then
the simple zeros of
$J_{s_\Rssf}(X,X)$
are exactly
the class invariants
$j(\kfrak_\Rssf)$
of $\Rssf$,
and so
$H_\Rssf(X)$
is the product of
those factors of
$J_{s_\Rssf}(X,X)$
which are
irreducible over $\Pssf$
and which divide
$J_{s_\Rssf}(X,X)$
only to the first power;
thus
$H_\Rssf$
has rational coefficients.

Now let us consider
two distinct orders
$\Rssf, \Rssf'$
with
$s_\Rssf = s_{\Rssf'}$.
The simple zeros of
$J_{s_\Rssf}(X,X)$
are then exactly
the class invariants
$j(\kfrak_\Rssf)$
of $\Rssf$
and those
$j(\kfrak_{\Rssf'})$
of $\Rssf'$,
and it follows as above that
$H_\Rssf(X)\,H_{\Rssf'}(X)$
is a polynomial
with rational coefficients.

First case:
$\Rssf
 =
 \Rssf_f(\sqrt{-d})
$,
$\Rssf'
 =
 \Rssf_{2f}(\sqrt{-d})
$,
$-d \equiv 1
 \bmod
 4
$,
but
$\Rssf
 \neq
 \Rssf_1(\sqrt{-3})
$.
Then
there exists
a primitive element
$\lambda$
in $\Rssf$
with norm
$\dfrac{d{+}1}{4}\,f^2$,
namely
$\lambda
 =
 f
 \,
 \dfrac{1{+}\sqrt{-d}}{2}
$,
but none
in $\Rssf'$,
because
the norm of
a primitive element
\[
  \lambda'
  \ =\ 
  a
  +
  b
  \cdot
  2f
  \,
  \dfrac{1{+}\sqrt{-d}}{2}
  \ =\ 
  a'
  +
  b
  \,
  f
  \,
  \sqrt{-d}
\]
of $\Rssf'$,
when $b = 0$
and thus
$\lambda' = \pm 1$,
satisfies
\[
  \Nm\,\lambda'
  \ =\ 
  1
  \ <\ 
  \dfrac{d{+}1}{4}
  \,
  f^2
  ,
\]
and
when $b \neq 0$,
satisfies
\[
  \Nm\,\lambda'
  \ =\ 
  {a'}^2
  +
  b^2
  \,
  f^2
  \,
  d
  \ \geqslant\ 
  d
  \,
  f^2
  \ >\ 
  \dfrac{d{+}1}{4}
  \,
  f^2
  .
\]
It follows that
$J_{\frac{d+1}{4} f^2}(X,X)$
is divisible by
$H_\Rssf(X)$,
but relatively prime to
$H_{\Rssf'}(X)$,
whence
\[
  H_\Rssf(X)
  \ =\ 
  \bigl(\,
  H_\Rssf(X)\,H_{\Rssf'}(X)
  \,,\,
  J_{\frac{d+1}{4} f^2}(X,X)
  \,\bigr)
\]
has rational coefficients,
and then
the same holds for
the second factor
$H_{\Rssf'}(X)$
as well.

Second case:
$\Rssf
 =
 \Rssf_1(\sqrt{-3})
$,
$\Rssf'
 =
 \Rssf_2(\sqrt{-3})
$.
Here
the first factor
$H_\Rssf(X)$
can be easily
determined:
$\Rssf$
has only one class,
and its invariant
is calculated
from the basis
$1
 ,
 \zeta_3
$
to be
\[
  j(\zeta_3)
  \ =\ 
  0
  ,
  \qquad
  \text{so}
  \quad
  H_\Rssf(X)
  \ =\ 
  X
  .
\]

Third case:
$\Rssf
 =
 \Rssf_1(\sqrt{-1})
$,
$\Rssf'
 =
 \Rssf_1(\sqrt{-2})
$.
Here
$\Rssf_1(\sqrt{-1})$
has only one class,
whose invariant
is calculated
via the basis
$1
 ,
 i
$
to be
\[
  j(i)
  \ =\ 
  2^6 \cdot 3^3
  ,
  \qquad
  \text{so}
  \quad
  H_\Rssf(X)
  \ =\ 
  X - 2^6 \cdot 3^3
  .
\]

We are getting to
the proof of
the congruence~\eqref{eqn:33},
but before that,
let us insert
an investigation on


\paragraph{13. The singular values of
the functions $\varphi_S(\omega)$.
First part.}

The singular values
$\varphi_S(\alpha)$,
with
$\alpha$
imaginary quadratic,
are algebraic integers.
This follows from
the equation
\[
  \Phi_s
  (\,
  \varphi_S(\alpha)
  \,,\,
  j(\alpha)
  \,)
  \ =\ 
  0
  ,
  \qquad
  s = |S|
  ,
\]
because
$j(\alpha)$
as a singlar class invariant
is an (algebraic) integer,
and the coefficients of
$\Phi_s(X,j(\omega))$
are rational integers.

From~\eqref{eqn:17}
it follows that
the product
\[
  \prod_{\nu=1}^{\psi(s)}
  \varphi_{S_\nu}(\alpha)
\]
extended over
all classes of level~$s$
is a rational integer,
which increases as
a power of~$s$;
for a prime number
$s = p$,
one has
by~\eqref{eqn:18}
more precisely
\[
  \prod_{\nu=1}^{p+1}
  \varphi_{P_\nu}(\alpha)
  \ =\ 
  p^{12}
  .
  \tag{36}
  \label{eqn:36}
\]

The construction of
the functions
$\varphi_S(\omega)$
allows us to
determine
the divisor
generated by
$\varphi_S(\omega)$;
but 
we only do this
for a few
special cases
here.

The determinant of
$S = P$
is a prime number~$p$,
which firstly
does not remain prime
in the quadratic field
$\Sigma = \Pssf(\alpha)$
generated by $\alpha$,
thus
becoming
the product of
two equal or distinct
degree-one
prime divisors
$\pfrak$ and $\pfrakbar$,
and secondly
does not divide
the conductor of the order
$\Rssf$ of $\Sigma$
for which
$\alpha$ is
the ideal basis quotient.

\textit{
Let
$\alpha
 =
 \dfrac{\alpha_1}{\alpha_2}
$,
where
$\alpha_1,\alpha_2$
is a basis of
an ideal
$\afrak_\Rssf$
of $\Rssf$,
and let
$P_\pfrak$
be a matrix
of determinant~$p$
---
which is
uniquely determined
up to equivalence
---
for which
$P_\pfrak
 \begin{pmatrix}
 \alpha_1
 \\
 \alpha_2
 \end{pmatrix}
$
is a basis of
$\afrak_\Rssf
 \,
 \pfrak_\Rssf
$.
Then one has
\[
  \varphi_{P_\pfrak}(\alpha)
  \ \approx\ 
  \pfrakbar^{12}
  .
  \tag{37}
\]
}

\noindent
Let $\cfrak_\Rssf$ be
an integral
characteristic ideal of
$\Rssf$
prime to~$p$,
so that
$\cfrak_\Rssf
 \,
 \pfrak_\Rssf
$
becomes
a principal ideal
$\gamma\,\Rssf$,
and
let
\[
  C
  \,
  P_\pfrak
  \begin{pmatrix}
  \alpha_1
  \\
  \alpha_2
  \end{pmatrix}
\]
be a basis of
$\cfrak_\Rssf
 \,
 \pfrak_\Rssf
 \,
 \afrak_\Rssf
 =
 \gamma
 \,
 \afrak_\Rssf
$.
Then
$|C|
 =
 c
 =
 \Nm\,\cfrak
$
and
$\Nm\,\gamma
 =
 p
 \,
 c
$.
We assume
$\cfrak_\Rssf$
is without any rational part,
which is obviously
allowed;
so $C$ is primitive
and we have
\[
  \begin{aligned}
  \varphi_{CP_\pfrak}(\alpha)
  &
  \ =\ 
  c^{12}
  \,
  p^{12}
  \,
  \dfrac
  {\Delta
   \left(
   C
   \,
   P_\pfrak
   \begin{pmatrix}
   \alpha_1
   \\
   \alpha_2
   \end{pmatrix}
   \right)
  }
  {\Delta
   \begin{pmatrix}
   \alpha_1
   \\
   \alpha_2
   \end{pmatrix}
  }
  \ =\ 
  c^{12}
  \,
  p^{12}
  \,
  \dfrac
  {\Delta
   \left(
   \begin{pmatrix}
   \alpha_1
   \\
   \alpha_2
   \end{pmatrix}
   \gamma
   \right)
  }
  {\Delta
   \begin{pmatrix}
   \alpha_1
   \\
   \alpha_2
   \end{pmatrix}
  }
  \\[2ex]
  &
  \ =\ 
  c^{12}
  \,
  p^{12}
  \,
  \gamma^{-12}
  \ =\ 
  \gammabar^{12}
  \ \approx\ 
  \cfrakbar^{12}
  \,
  \pfrakbar^{12}
  .
  \end{aligned}
\]
On the other hand,
it follows directly
from the definition of
$\varphi_S(\omega)$
that 
one has the formula
\[
  \varphi_{CP_\pfrak}(\omega)
  \ =\ 
  \varphi_C(P_\pfrak(\omega))
  \,
  \varphi_{P_\pfrak}(\omega)
  ,
\]
whence
\[
  \varphi_{CP_\pfrak}(\alpha)
  \ =\ 
  \varphi_C(P_\pfrak(\alpha))
  \,
  \varphi_{P_\pfrak}(\alpha)
  ,
\]
and since
$\varphi_C(P_\pfrak(\alpha))$
as a divisor of
a power of $c$
is prime to~$p$,
it follows that
\[
  \varphi_{P_\pfrak}(\alpha)
  \ \approx\ 
  \pfrakbar^{12}
  .
\]
Naturally,
corresponding to
a matrix
$P_\pfrakbar$
for which
$P_\pfrakbar
 \begin{pmatrix}
 \alpha_1
 \\
 \alpha_2
 \end{pmatrix}
$
is a basis of
$\afrak_\Rssf
 \,
 \pfrakbar_\Rssf
$,
one has
\[
  \varphi_{P_\pfrakbar}(\alpha)
  \ \approx\ 
  \pfrak^{12}
  .
  \tag{38}
  \label{eqn:38}
\]

Let us now assume that
$\pfrak$ and $\pfrakbar$
are distinct.
Then
it follows from~\eqref{eqn:36}
that
\[
  \begin{array}{l}
  \textit{
  $\varphi_P(\alpha)$
  is a unit
  for all $P$ of determinant $p$
  }
  \\
  \qquad\qquad
  \textit{
  which is neither equivalent
  to $P_\pfrak$
  nor to $P_\pfrakbar$.
  }
  \end{array}
  \tag{39}
  \label{eqn:39}
\]

By means of
this divisor representation,
we prove the following
theorem:
\footnote{
\textit{Hasse},
loc. cit.
{\ \ \llap{${}^{2)}$}\ }
123 and 125.
}

\textit{
Let
$p$, $\pfrak$, $\pfrakbar$
and $\alpha$
have the same meaning
as above,
but suppose
$\pfrak \neq \pfrakbar$.
Let $f_P(\omega)$ be
an entire integral function
from $\Pssf_{\Mfrak_P}$,
and moreover
let the coefficients of
the $q$-expansion of
$f_{\left(\begin{smallmatrix}
        p & 0
        \\
        0 & 1
        \end{smallmatrix}\right)}
 (\omega)
$
be divisible by~$p$.
Then
the (complex) number
$f_{P_\pfrakbar}(\alpha)$
is algebraic
and divisible by $\pfrak$.
}

According to~\S\textbf{3}
one has
\[
  f_P(\omega)
  \,
  \Phi_p'(\varphi_P(\omega),j(\omega))
  \ =\ 
  a_0(j(\omega))
  \ +\ 
  a_1(j(\omega))
  \,
  \varphi_P(\omega)
  \ +\ 
  \cdots
  \ +\ 
  a_p(j(\omega))
  \,
  \varphi_P(\omega)^p
  ,
\]
with polynomials
$a_\nu(j(\omega))$
of $j(\omega)$
whose coefficients
are rational integers;
the coefficients
of $a_0(j(\omega))$
are in fact
divisible by~$p$.
Let us put
$\omega = \alpha$
here,
from which
it follows that
$f_P(\alpha)$
is algebraic,
and
if we take
$P = P_\pfrakbar$,
then
by~\eqref{eqn:38}
it follows that
\[
  f_{P_\pfrakbar}(\alpha)
  \,
  \Phi_p'(\varphi_{P_\pfrakbar}(\alpha),j(\alpha))
  \ \equiv\ 
  0
  \mod
  \pfrak
  ,
\]
but because
\[
  \begin{aligned}
  \Phi_p'(\varphi_{P_\pfrakbar}(\alpha),j(\alpha))
  &
  \ =\ 
  \prod_{P \not\sim P_\pfrakbar}
  \ 
  \bigl(\,
  \varphi_{P_\pfrakbar}(\alpha)
  -
  \varphi_P(\alpha)  
  \,\bigr)
  \\
  &
  \ \equiv\ 
  (-1)^p
  \,
  \prod_{P \not\sim P_\pfrakbar}
  \ 
  \varphi_{P_\pfrak}(\alpha)
  \mod
  \pfrak
  \end{aligned}
\]
is relatively prime to~$\pfrak$
by~\eqref{eqn:38} and~\eqref{eqn:39},
the assertion follows.


\paragraph{14. Proof of the fundamental congruence~\eqref{eqn:33}}

It is a matter of
showing the congruence
\[
  j
  (\,
  \kfrak_\Rssf
  \,
  \pfrak_\Rssf^{-1}
  \,)
  \ -\ 
  j(\kfrak_\Rssf)^p
  \ \equiv\ 
  0
  \mod
  \pfrak
\]
for degree-one prime divisors
$\pfrak$
of norm~$p$
not dividing
the conductor of $\Rssf$.
In the case of
$\pfrak = \pfrakbar$,
this follows
from~\eqref{eqn:30},
because
in this case
the two possibilities
undetermined in
\eqref{eqn:30}
obviously
say the same thing.
But in the case of
$\pfrak \neq \pfrakbar$,
we can apply
the theorem just proven,
because
the ideal
$\afrak_\Rssf$
with the basis
$\alpha_1,\alpha_2$
can be taken from
$\kfrak_\Rssf$,
thus
$\afrak_\Rssf\,\pfrakbar_\Rssf$
is an ideal from
$\kfrak_\Rssf\,\pfrak_\Rssf^{-1}$,
and with
\[
  f_P(\omega)
  \ =\ 
  j(P(\omega))
  -
  j(\omega)^p
\]
one now has
\[
  j
  (\,
  \kfrak_\Rssf
  \,
  \pfrak_\Rssf^{-1}
  \,)
  \ -\ 
  j(\kfrak_\Rssf)^p
  \ =\ 
  f_{P_\pfrak}(\alpha)
  .
\]
The conjugate functions
$f_P(\omega)$
obviously have
integer $q$-coefficients,
and
they are divisible by~$p$,
as seen from
$f_{\left(\begin{smallmatrix}
        p & 0
        \\
        0 & 1
        \end{smallmatrix}\right)}
 (\omega)
 =
 j(p\omega)
 -
 j(\omega)^p
$
and because
\[
  j(p\omega)
  \ \equiv\ 
  j(\omega)^p
  \mod
  p
  .
\]
%


\paragraph{15. The isomorphism $\sigma$
of $\Rfrak_\Rssf$ with $\Gfrak^{\Omega_\Rssf}_{\Sigma}$.}

Let $\pfrak$ be
a degree-one prime divisor
of $\Sigma$
not dividing
the conductor of $\Rssf$,
and which is
relatively prime to
all the differences
$j(\kfrak_\Rssf)
 -
 j(\kfrak_\Rssf')
$,
$\kfrak_\Rssf
 \neq
 \kfrak_\Rssf'
$,
or,
what amounts to the same thing,
to the discriminant
$D(H_\Rssf)$
of the class polynomial
$H_\Rssf(X)$.
Let $\Pfrak$ be
a prime factor of $\pfrak$
in $\Omega_\Rssf$
and
let $F_\Pfrak$ be
the \textit{Frobenius} automorphism
of $\Pfrak$
over $\Sigma$;
then
one has
\[
  j(\kfrak_\Rssf)^{\Nm\pfrak}
  \ \equiv\ 
  j(\kfrak_\Rssf)
  \,
  F_\Pfrak
  \mod
  \Pfrak
  .
  \tag{40}
\]
Hence
by~\eqref{eqn:33}
one has
\[
  j(\kfrak_\Rssf)
  \,
  F_\Pfrak
  \ \equiv\ 
  j(\kfrak_\Rssf\,\pfrak_\Rssf^{-1})
  \mod
  \Pfrak
  ,
\]
and because
$j(\kfrak_\Rssf)
 \,
 F_\Pfrak
$
is one of
$j(\kfrak_\Rssf')$
---
$H_\Rssf(x)$
having rational coefficients
---,
one has
\[
  j(\kfrak_\Rssf\,\pfrak_\Rssf^{-1})
  \ =\ 
  j(\kfrak_\Rssf)
  \,
  F_\Pfrak
  .
  \tag{41}
  \label{eqn:41}
\]
Now let
$\hfrak_\Rssf$
be a second class
of the order $\Rssf$
and
let
$\pfrak_\Rssf^{(1)}
 \pfrak_\Rssf^{(2)}
 \cdots
 \pfrak_\Rssf^{(m)}
$
be an integral ideal
from $\hfrak_\Rssf$
prime to
the conductor of $\Rssf$
and to
$D(H_\Rssf)$.
As degree-two
prime divisors
$\pfrak$
determine
principal ideals
$\pfrak_\Rssf$
of $\Rssf$,
the $\pfrak^{(\nu)}$
can be assumed
to be of degree-one;
then \eqref{eqn:41}
yields
\[
  j(\kfrak_\Rssf
    \hfrak_\Rssf^{-1})
  \ =\ 
  j(\kfrak_\Rssf)
  \,
  F_{\Pfrak^{(1)}}
  F_{\Pfrak^{(2)}}
  \cdots
  F_{\Pfrak^{(m)}}
  ,
\]
that is to say:
there exists
an automorphism
$\sigma(h_\Rssf)$
of $\Omega_\Rssf/\Sigma$
with
\[
  j(\kfrak_\Rssf)
  \,
  \sigma(h_\Rssf)
  \ =\ 
  j(\kfrak_\Rssf
    \hfrak_\Rssf^{-1})
  \tag{42}
  \label{eqn:42}
\]
for all classes
$\kfrak_\Rssf$
of $\Rssf$.
Also,
$\sigma(\hfrak_\Rssf)$
is uniquely determined
by~\eqref{eqn:42}.
On one hand,
it follows from~\eqref{eqn:42}
that
\[
  \sigma(\hfrak_\Rssf^{(1)})
  \,
  \sigma(\hfrak_\Rssf^{(1)})
  \ =\ 
  \sigma(
  \hfrak_\Rssf^{(1)}
  \hfrak_\Rssf^{(2)}
  )
  ,
\]
while
on the other hand,
$\sigma(\hfrak_\Rssf) = 1$
can hold
only for
$\hfrak_\Rssf = \efrak_\Rssf$:
thus
$\sigma$
maps $\Rfrak_\Rssf$
isomorphically \emph{into}
$\Gfrak^{\Omega_\Rssf}_{\Sigma}$.

In fact
one has
$\sigma(\Rfrak_\Rssf)
 =
 \Gfrak^{\Omega_\Rssf}_{\Sigma}
$,
thus
\textit{
$\sigma$ is
an isomorphism
of the class group of
$\Rssf$
with
the Galois group of
$\Omega_\Rssf/\Sigma$.
}

For that,
we first show that
for a single invariant
$j(\hfrak_\Rssf)$,
one already has
\[
  \Sigma(j(\hfrak_\Rssf))
  \ =\ 
  \Omega_\Rssf
  ,
\]
which also justifies
the definition of
$\Omega_\Rssf$
in~\S\textbf{11}
deviating from~\S\textbf{9}.
But that is to say:
any automorphism
$\lambda$
of $\Omega_\Rssf/\Sigma(j(\hfrak_\Rssf))$
is the identity,
or,
for any class
$\kfrak_\Rssf$,
one has
\[
  j(\kfrak_\Rssf)
  \,
  \lambda
  \ =\ 
  j(\kfrak_\Rssf)
  .
  \tag{43}
  \label{eqn:43}
\]

If~\eqref{eqn:43} holds
for a particular class
$\kfrak_\Rssf$,
then it also holds for
$\kfrak_\Rssf
 \pfrak_\Rssf^{-1}
$,
where
$\pfrak$ is understood
to be
a degree-one prime divisor
of $\Sigma$
prime to
$D(H_\Rssf)$
and to
the conductor of $\Rssf$:
one has
\[
  j(\kfrak_\Rssf\,\pfrak_\Rssf^{-1})
  \ \equiv\ 
  j(\kfrak_\Rssf)^{\Nm\pfrak}
  \mod
  \pfrak
  ,
\]
and so
since
$\pfrak\,\lambda
 =
 \pfrak
$,
one has
\[
  \begin{aligned}
  j(\kfrak_\Rssf\,\pfrak_\Rssf^{-1})
  \,
  \lambda
  &
  \ \equiv\ 
  \bigl(\,
  j(\kfrak_\Rssf)
  \,
  \lambda
  \,\bigr)^{\Nm\pfrak}
  \ =\ 
  j(\kfrak_\Rssf)^{\Nm\pfrak}
  \mod
  \pfrak
  ,
  \\
  j(\kfrak_\Rssf\,\pfrak_\Rssf^{-1})
  \,
  \lambda
  &
  \ \equiv\ 
  j(\kfrak_\Rssf
    \pfrak_\Rssf^{-1}
   )
  \mod
  \pfrak
  .
  \end{aligned}
\]
But the (algebraic) number
$j(\kfrak_\Rssf\,\pfrak_\Rssf^{-1})
 \,
 \lambda
$
is an invariant
$j(\kfrak_\Rssf')$,
because
$H_\Rssf(x)$
has rational coefficients,
so
it must be equal to
$j(\kfrak_\Rssf\,\pfrak_\Rssf^{-1})$.
From the obviously valid
equation
$j(\hfrak_\Rssf)
 \,
 \lambda
 =
 j(\hfrak_\Rssf)
$,
we can get to
$j(\kfrak_\Rssf)
 \,
 \lambda
 =
 j(\kfrak_\Rssf)
$
for an arbitrary
$\kfrak_\Rssf$
by multiplication with
a product of
degree-one prime divisors.

Since
$j(\hfrak_\Rssf)$
as a zero of
the rational integer polynomial
$H_\Rssf(X)$
has degree
at most~$h_\Rssf$
over $\Sigma$,
it follows that
$\Gfrak^{\Omega_\Rssf}_{\Sigma}
 =
 \Gfrak^{\Sigma(j(\hfrak_\Rssf))}_{\Sigma}
$
has order
at most~$h_\Rssf$,
and from that
it follows that
$\sigma(\Rfrak_\Rssf)
 =
 \Gfrak^{\Omega_\Rssf}_{\Sigma}
$.
At the same time
it is proven that
\textit{
the class polynomial
$H_\Rssf(X)$
is irreducible
over $\Sigma$.
}


\paragraph{16. The Galois group
$\Gfrak^{\Omega_\Rssf}_{\Pssf}$.}

Obviously
$\Omega_\Rssf$
is also normal
over $\Pssf$.
Let us compute
the Galois group
$\Gfrak^{\Omega_\Rssf}_{\Pssf}$.
Since one has
$[\,
 \Pssf(j(\hfrak_\Rssf))
 \,{:}\,
 \Pssf
 \,]
 =
 h_\Rssf
$,
it follows that
$[\,
 \Omega_\Rssf
 \,{:}\,
 \Pssf(j(\hfrak_\Rssf))
 \,]
 =
 2
$,
so
$\Gfrak^{\Omega_\Rssf}_{\Pssf(j(\hfrak_\Rssf))}$
is of order~2.
Let $\lambda(\hfrak_\Rssf)$ be
the automorphism of
$\Omega_\Rssf/\Pssf(j(\hfrak_\Rssf))$
distinct from
the identity.
For that
one has
\[
  j(\hfrak_\Rssf\kfrak_\Rssf)
  \,
  \lambda(\hfrak_\Rssf)
  \ =\ 
  j(\hfrak_\Rssf\kfrak_\Rssf^{-1})
  \qquad
  \text{for every $\kfrak_\Rssf$}
  .
  \tag{44}
  \label{eqn:44}
\]
For $\kfrak_\Rssf = \efrak_\Rssf$,
the equation
is valid.
Hence
it suffices to prove that,
assuming it is valid
for a particular class
$\kfrak_\Rssf$,
it is then so for
$\kfrak_\Rssf\,\pfrak_\Rssf^{-1}$,
where again
$\pfrak$
is of degree-one
and is prime to
$D(H_\Rssf)$
and to
the conductor of $\Rssf$.
Now
$\lambda(\hfrak_\Rssf)$
not the identity
on $\Sigma$,
whence
one has
$\pfrak
 \,
 \lambda(\hfrak_\Rssf)
 =
 \pfrakbar
$,
and from
\[
  j(\hfrak_\Rssf\kfrak_\Rssf\pfrak_\Rssf^{-1})
  \ \equiv\ 
  j(\hfrak_\Rssf\kfrak_\Rssf)^{\Nm\pfrak}
  \mod
  \pfrak  
\]
and~\eqref{eqn:44},
it follows that
\[
  j(\hfrak_\Rssf\kfrak_\Rssf\pfrak_\Rssf^{-1})
  \,
  \lambda(\hfrak_\Rssf)
  \ =\ 
  j(\hfrak_\Rssf\kfrak_\Rssf^{-1})^{\Nm\pfrak}
  \mod
  \pfrakbar
  ;
\]
but
on the other hand
one also has
\[
  j(\hfrak_\Rssf\kfrak_\Rssf^{-1})^{\Nm\pfrak}
  \ \equiv\ 
  j(\hfrak_\Rssf\kfrak_\Rssf^{-1}\pfrak_\Rssf)
  \mod
  \pfrakbar
  ,
\]
whence
\[
  j(\hfrak_\Rssf\kfrak_\Rssf\pfrak_\Rssf^{-1})
  \,
  \lambda(\hfrak_\Rssf)
  \ =\ 
  j(\hfrak_\Rssf
    \,
    (\kfrak_\Rssf\pfrak_\Rssf^{-1})^{-1}
   )
  .
\]
The Galois group
$\Gfrak^{\Omega_\Rssf}_{\Pssf}$
consists of all
$\sigma(\kfrak_\Rssf)$,
$\kfrak_\Rssf \in \Rfrak_\Rssf$,
and all
$\lambda(\hfrak_\Rssf)
 \sigma(\kfrak_\Rssf)
$,
$\kfrak_\Rssf \in \Rfrak_\Rssf$
for fixed
$\hfrak_\Rssf$.
The structure of
$\Gfrak^{\Omega_\Rssf}_{\Pssf}$
is known
if we known
the order of
$\lambda(\hfrak_\Rssf)$
and understand
to which
$\sigma(\kfrak_\Rssf')$
a $\sigma(\kfrak_\Rssf)$
is transformed by
$\lambda(\hfrak_\Rssf)$.
It follows directly
from~\eqref{eqn:44}
that
\[
  \lambda(\hfrak_\Rssf)^2
  \ =\ 
  1
  ,
  \tag{45}
\]
and furthermore
one has
\[
  j(\hfrak_\Rssf)
  \,
  \lambda(\hfrak_\Rssf)^{-1}
  \,
  \sigma(\kfrak_\Rssf)
  \,
  \lambda(\hfrak_\Rssf)
  \ =\ 
  j(\hfrak_\Rssf\kfrak_\Rssf^{-1})
  \,
  \lambda(\hfrak_\Rssf)
  \ =\ 
  j(\hfrak_\Rssf\kfrak_\Rssf)
  \ =\ 
  j(\hfrak_\Rssf)
  \,
  \sigma(\kfrak_\Rssf^{-1})
\]
and because
$\lambda(\hfrak_\Rssf)^{-1}
 \,
 \sigma(\kfrak_\Rssf)
 \,
 \lambda(\hfrak_\Rssf)
$
lies in
$\Gfrak^{\Omega_\Rssf}_{\Sigma}$,
it is one of
$\sigma(\kfrak_\Rssf')$,
but
on the other hand
it is uniquely determined
by its effect on
the primitive element
$j(\hfrak_\Rssf)$
of $\Omega_\Rssf/\Sigma$,
so one has
\[
  \lambda(\hfrak_\Rssf)^{-1}
  \,
  \sigma(\kfrak_\Rssf)
  \,
  \lambda(\hfrak_\Rssf)
  \ =\ 
  \sigma(\kfrak_\Rssf^{-1})
  .
  \tag{46}
\]

The automorphism
$\lambda(\hfrak_\Rssf)
 \,
 \sigma(\hfrak_\Rssf)^{-2}
$
is independent of $\hfrak_\Rssf$
because
$j(\kfrak_\Rssf)
 \cdot
 \lambda(\hfrak_\Rssf)
 \,
 \sigma(\hfrak_\Rssf)^{-2}
 =
 j(\kfrak_\Rssf^{-1})
$;
consequently
one has
\[
  \lambda(\hfrak_\Rssf')
  \ =\ 
  \lambda(\hfrak_\Rssf)
  \,
  \sigma({\hfrak_\Rssf'}^2 \hfrak_\Rssf^{-2})
  ,
  \qquad
  \text{in particular}
  \qquad
  \lambda(\hfrak_\Rssf)
  \ =\ 
  \lambda(\efrak_\Rssf)
  \,
  \sigma(\hfrak_\Rssf^2)
\]
and
\[
  \lambda(\hfrak_\Rssf)
  \,
  \lambda(\hfrak_\Rssf')
  \ =\ 
  \sigma(\hfrak_\Rssf^{-2} {\hfrak_\Rssf'}^2)
  .
\]
The automorphism
\textit{
$\lambda(\efrak_\Rssf)$
is the transition to
the complex conjugate,
}
because one has
$\alpha \cdot \lambda(\efrak_\Rssf)
 =
 \alphabar
$
for $\alpha \in \Sigma$
and
$j(\kfrak_\Rssf)
 \cdot
 \lambda(\efrak_\Rssf)
 =
 j(\kfrak_\Rssf^{-1})
$,
and
$j(\kfrak_\Rssf)$
and
$j(\kfrak_\Rssf^{-1})$
are complex conjugates
because,
taking
$\alpha_1,\alpha_2$
to be a basis of
an ideal $\afrak_\Rssf$
from $\kfrak_\Rssf$,
one has
$\alphabar_1,-\alphabar_2$
as a basis of
the ideal $\afrakbar_\Rssf$
from $\kfrak_\Rssf^{-1}$,
and
$j(\kfrak_\Rssf)
 =
 j
 \left(
 \dfrac
 {\alpha_1}
 {\alpha_2}
 \right)
$
and
$j(\kfrak_\Rssf^{-1})
 =
 j
 \left(
 \dfrac
 {\alphabar_1}
 {-\alphabar_2}
 \right)
$
are thus
values of the function
$j(\omega)$
at points
which are mirror-images in
the imaginary axis,
on which
$j(\omega)$
takes real values.

$\Pssf(j(\efrak_\Rssf))$
is the maximal real subfield
of $\Omega_\Rssf$.


\paragraph{17. The reciprocity law for
the ring class field
$\Omega_\Rssf/\Sigma$.}

As $\Omega_\Rssf$
is abelian
over $\Sigma$,
the \textit{Frobenius} automorphisms
$F_\Pfrak$
of all the
$\Omega_\Rssf$-prime factors
$\Pfrak$
of a prime divisor
$\pfrak$ of $\Sigma$
coincide
and henceforth
will be denoted by
$F_\pfrak$.
Let us now prove that
for the prime divisors
$\pfrak$ of $\Sigma$
unramified in $\Omega_\Rssf$,
the decomposition assertion of
the reciprocity law
holds:
\[
  F_\pfrak
  \ =\ 
  \sigma(\efrak_\Rssf\,\pfrak_\Rssf)
  ;
  \tag{47}
\] 
the proof here
accomplishes this
for almost all $\pfrak$.

If $\pfrak$ is
of degree-one
and is relatively prime to
the conductor of $\Rssf$
and to $D(H_\Rssf)$,
this follows directly from
the definition of
$\sigma(\kfrak_\Rssf)$.
For the $\pfrak$
of degree-two
unramified in $\Omega_\Rssf$,
we prove it
as follows:
since
one has
$\efrak_\Rssf
 \,
 \pfrak_\Rssf
 =
 \efrak_\Rssf
$
in this case,
it must be shown that
the decomposition group
$\Zfrak^{\Omega_\Rssf}_{\Sigma}(\Pfrak)$
of the $\Omega_\Rssf$-prime factor
$\Pfrak$ of $\pfrak$
over $\Sigma$
has order~1.
As $\Sigma$
does not lie in
the decomposition field of
$\Pfrak$
over $\Pssf$,
it follows
on one hand
that
\[
  [\,
  \Zfrak^{\Omega_\Rssf}_{\Pssf}(\Pfrak)
  \,{:}\,
  \Zfrak^{\Omega_\Rssf}_{\Sigma}(\Pfrak)
  \,]
  \ =\ 
  2
  ,
\]
and
on the other hand
that
a generating element of
the cyclic group
$\Zfrak^{\Omega_\Rssf}_{\Pssf}(\Pfrak)$
cannot lie in
$\Gfrak^{\Omega_\Rssf}_{\Sigma}$,
and so must be
one of
$\lambda(\hfrak_\Rssf)$.
But that means
\[
  [\,
  \Zfrak^{\Omega_\Rssf}_{\Pssf}(\Pfrak)
  \,{:}\,
  1
  \,]
  \ =\ 
  2
  ;
  \qquad
  [\,
  \Zfrak^{\Omega_\Rssf}_{\Sigma}(\Pfrak)
  \,{:}\,
  1
  \,]
  \ =\ 
  1
  .
\]


\paragraph{18. The genus field of the order $\Rssf$.}

The \emph{genus field $A_\Rssf$ of $\Rssf$}
can be defined as
the maximal subfield of $\Omega_\Rssf$
which is abelian over $\Pssf$.
It is the fixed field of
the commutator group of
$\Gfrak^{\Omega_\Rssf}_{\Pssf}$;
since
the commutator of
two $\sigma(\kfrak_\Rssf)$
equals~1,
the commutator of
$\lambda(\hfrak_\Rssf)
 \sigma(\kfrak_\Rssf)
$
and
$\sigma(\kfrak_\Rssf')$
equals
$\sigma({\kfrak_\Rssf'}^{-2})$
and
the commutator of
$\lambda(\hfrak_\Rssf)
 \sigma(\kfrak_\Rssf)
$
and
$\lambda(\hfrak_\Rssf')
 \sigma(\kfrak_\Rssf')
$
equals
$\sigma((\kfrak_\Rssf^{-1}{\kfrak_\Rssf'})^2)$,
it follows that
this commutator group
is equal to
$\sigma(\Rfrak_\Rssf^2)$,
where
$\Rfrak_\Rssf^2$
consists of
all the square classes
$\kfrak_\Rssf^2$
and stands for
the \emph{principal genus of $\Rssf$}.
The kernel of
the homomorphism
$\kfrak_\Rssf \to \kfrak_\Rssf^2$
from $\Rfrak_\Rssf$
onto $\Rfrak_\Rssf^2$
is the group
$\Rfrak_\Rssf^{(0)}$
of \emph{ambiguous classes of $\Rssf$}.
The index
$[\,
 \Rfrak_\Rssf
 \,{:}\,
 \Rfrak_\Rssf^{(0)}
 \,]
$
is consequently
equal to
the order of
$\Rfrak_\Rssf^2$,
and hence
the index
$[\,
 \Rfrak_\Rssf
 \,{:}\,
 \Rfrak_\Rssf^2
 \,]
$,
which is the
\emph{number of genus of $\Rssf$}
(residue classes of
$\Rfrak_\Rssf$
modulo the principal genus),
\textit{
is equal to
the number of
ambiguous classes
}
(order of $\Rfrak_\Rssf^{(0)}$).
The genus group
$\Rfrak_\Rssf/\Rfrak_\Rssf^2$
is (like $\Rfrak_\Rssf^{(0)}$)
of the type
$(
 \underset{1}
 2
 \,,\,
 \underset{2}
 2
 \,,\,
 \ldots
 \,,\,
 \underset{a}
 2
 )
$.
The factor group
$\Gfrak^{\Omega_\Rssf}_{\Pssf}
 /
 \sigma(\Rfrak_\Rssf^2)
$
is then
of the type
$(
 \underset{1}
 2
 \,,\,
 \underset{2}
 2
 \,,\,
 \ldots
 \,,\,
 \underset{a}
 2
 \,,\,
 \underset{a{+}1}
 2
 )
$,
which proves that
\textit{
the genus field
$A_\Rssf$
is the composite of
$a{+}1$
independent
quadratic fields
},
of which
one can be taken
to be equal to $\Sigma$.

By~\eqref{eqn:44}
one has
$\lambda(\kfrak_\Rssf')
 =
 \lambda(\kfrak_\Rssf)
$
if and only if
$\kfrak_\Rssf^2
 =
 {\kfrak_\Rssf'}^2
$;
thus
\textit{
one has
$\Pssf(j(\kfrak_\Rssf))
 =
 \Pssf(j(\kfrak_\Rssf'))
$
if and only if
the class
$\kfrak_\Rssf{\kfrak_\Rssf'}^{-1}$
is ambiguous.
}


\paragraph{19. The correspondence theorem for the ring class field.}

\textit{
Let $\Rssf$ and $\Rssf'$ be
orders of
an imaginary quadratic field
$\Sigma$,
and suppose
$\Rssf \supseteq \Rssf'$.
Then
one has
$\Omega_\Rssf
 \subseteq
 \Omega_{\Rssf'}
$,
and an automorphism
$\sigma(\hfrak_{\Rssf'})$
of $\Omega_{\Rssf'}/\Sigma$
induces on $\Omega_\Rssf$
the automorphism
$\sigma(\hfrak_\Rssf)$
when the class
$\hfrak_{\Rssf'}$ from $\Rssf'$
belongs to
the class
$\hfrak_\Rssf$ from $\Rssf$;
the Galois group
$\Gfrak^{\Omega_{\Rssf'}}_{\Omega_\Rssf}$
thus consists of
all $\sigma(\hfrak_{\Rssf'})$
for the $\hfrak_{\Rssf'}$
belonging to
the principal class of $\Rssf$.
}

Proof.
\footnote{
\textit{Hasse},
Monatsh. f. Math. u. Phys.
\textbf{38},
p.~325
(1931).
}
The conductor~$f'$
of $\Rssf'$
is divisible by
the conductor~$f$
of $\Rssf$;
let us say
$f'/f = s$.
Let a class
$\kfrak_\Rssf$
of $\Rssf$
be given,
let $\kfrak_{\Rssf'}$ be
a class of $\Rssf'$
contained in
$\kfrak_\Rssf$.
If
$j(\kfrak_{\Rssf'})
 =
 j(\alpha)
$,
then
by~\S\textbf{8},
exactly one zero
$j(S(\alpha))$
of
${\displaystyle
 J_S(X,j(\kfrak_{\Rssf'}))
 =
 \prod_S
 (X - j(S(\alpha)))
}$
is a class invariant of $\Rssf$,
and in fact
equal to
$j(\kfrak_\Rssf)$.
Thus
$X-j(\kfrak_\Rssf)$
is the greatest common divisor of
$J_S(X,j(\kfrak_{\Rssf'}))$
and
$H_\Rssf(X)$,
which proves that
$j(\kfrak_\Rssf)$
lies in $\Omega_{\Rssf'}$,
and thus
one has
$\Omega_\Rssf
 \subseteq
 \Omega_{\Rssf'}
$.
Then
$\sigma(\hfrak_{\Rssf'})$
sends
$X-j(\kfrak_\Rssf)$
to the greatest common divisor of
$J_S(X,j(\kfrak_{\Rssf'})\sigma(\hfrak_{\Rssf'}))
 =
 J_S(X,j(\kfrak_{\Rssf'}\hfrak_{\Rssf'}^{-1}))
$
and $H(X)$,
that is to say
to
$X-j(\kfrak_\Rssf\hfrak_{\Rssf'}^{-1})$:
\[
  j(\kfrak_\Rssf)
  \,
  \sigma(\hfrak_{\Rssf'})
  \ =\ 
  j(\kfrak_\Rssf\hfrak_{\Rssf'}^{-1})
  \ =\ 
  j(\kfrak_\Rssf\hfrak_\Rssf^{-1})
  \ =\ 
  j(\kfrak_\Rssf)
  \,
  \sigma(\hfrak_\Rssf)
  ,
\]
q.e.d.


\paragraph{20. The singular values of  functions
from $\Pssf_{\Mfrak_S}$.}

The fact that
a function
$f_S(\omega)$
from $\Pssf_{\Mfrak_S}$
is a rational function of
$j(\omega)$ and $j(S(\omega))$
with rational coefficients
does not allow one
to conclude that
the function value
$f_S(\alpha)$
for a (complex) number $\alpha$
lies in the field
$\Pssf(j(\alpha),j(S(\alpha)))$.
For an entire function
$f_S(\omega)$,
the product
$f_S(\omega)
 \,
 J_S'(j(S(\omega)),j(\omega))
$
is a rational integer polynomial of
$j(\omega)$ and $j(S(\omega))$,
and from that
it follows that
$f_S(\alpha)
 \in
 \Pssf(j(\alpha),j(S(\alpha)))
$
if
$j(S(\alpha))$
is a \emph{simple} zero of
$J_S(X,j(\alpha))$.
But
for imaginary quadratic $\alpha$
that is in general
not the case.

Let us assume that
among the conjugates
\[
  j(S_\nu(\omega))
  ,
  \quad
  \nu = 1,2,\ldots,\psi(s)
  ,
\]
of $j(S(\omega))$
over $\Pssf_{\Mfrak_S}$,
exactly~$r$ of them,
say
$j(S_1(\omega)),
 \ldots,
 j(S_r(\omega))
$,
take on the value
$j(S(\alpha))$
for $\omega = \alpha$,
so that
$S$ is equivalent to
one of the matrices
$S_1,\ldots,S_r$.
As we are going to show,
\textit{
$f_S(\alpha)$
is then algebraic over
$\Pssf(j(\alpha),j(S(\alpha)))$
and
the conjugates of
$f_S(\alpha)$
over this field
are among
$f_{S_1}(\alpha),
 \ldots,
 f_{S_r}(\alpha)
$.
}
For the proof,
let us denote by
$\Gfrak$
the Galois group of
the normal closure
$K
 =
 \Pssf(
 j(\omega),
 j(S_1(\omega)),
 \ldots,
 j(S_{\psi(s)}(\omega))
 )
$
of $\Pssf_{\Mfrak_S}$
over $\Pssf_{\Mfrak}$,
and by
$\Hfrak$
the subgroup of $\Gfrak$
of the automorphisms
which permute the functions
\[
  j(S_1(\omega))
  \,,\,
  \ldots
  \,,\,
  j(S_r(\omega))
\]
among one another,
and by
$K_0$
the fixed field
of $\Hfrak$
in $K$.
The polynomial
\[
  F(X,\omega)
  \ =\ 
  \prod_{\nu=1}^{r}
  (X - j(S_\nu(\omega)))
\]
is sent to itself
by every
$\sigma \in \Hfrak$;
in contrast,
for
$\sigma \in \Gfrak$,
$\sigma \notin \Hfrak$,
not only is
$F(X,\omega)\,\sigma
 \neq
 F(X,\omega)
$,
but in fact
\[
  F(X,\alpha)\,\sigma
  \ \neq\ 
  F(X,\alpha)
  ,
\]
whence
for suitable rational $a$,
\[
  \ell(\omega)
  \ =\ 
  F(a,\omega)
\]
is a primitive element of
$K_0/\Pssf_\Mfrak$
with the following property:
writing
$L(X,j(\omega))$
for the irreducible polynomial
over $\Pssf_\Mfrak$
with
$L(\ell(\omega),j(\omega)) = 0$,
one has
$\ell(\alpha)$
as a simple zero of
$L(X,j(\alpha))$.
For an entire function
$g(\omega)$
from $K_0$,
the product
$g(\omega)
 \,
 L'(\ell(\omega),j(\omega))
$
is a rational integer polynomial of
$j(\omega)$ and $\ell(\omega)$,
whence
$g(\alpha)$
is an element of
$\Pssf(j(\alpha),\ell(\alpha))$,
and thus of
$\Pssf(j(\alpha),j(S(\alpha)))$
because
$\ell(\alpha)
 =
 (a - j(S(\alpha)))^r
$.

For an entire function
$f_S(\omega)$
from $\Pssf_{\Mfrak_S}$,
the polynomial
${\displaystyle
 \prod_{\nu=1}^{r}
 (X-f_{S_\nu}(\omega))
}$
has entire functions
from $K_0$
as coefficients;
thus
${\displaystyle
 \prod_{\nu=1}^{r}
 (X-f_{S_\nu}(\alpha))
}$
has coefficients in
$\Pssf(j(\alpha),j(S(\alpha)))$,
from which
in fact
it follows that
$f_S(\alpha)$
is algebraic over
$\Pssf(j(\alpha),j(S(\alpha)))$
and that
its conjugates over this field
are among
$f_{S_1}(\alpha),
 \ldots,
 f_{S_r}(\alpha)
$.

In a special case,
it can also
be deduced that
$f_S(\alpha)
 \in
 \Pssf(j(\alpha),j(S(\alpha)))
$
when $r>1$:
\footnote{
\textit{Hasse},
J. f. Math.
\textbf{165},
82 to 83 (1931).
There,
a special function
$f_S(\alpha)$
is investigated
and
$\Rssf$
is the principal order.
}

\textit{
Let $\alpha_1,\alpha_2$ be
a basis of
a characteristic ideal
$\afrak_\Rssf$
of the order $\Rssf$
of $\Sigma$
and let
$S
 \begin{pmatrix}
 \alpha_1
 \\
 \alpha_2
 \end{pmatrix}
$
be likewise
a basis of
a characteristic ideal
$\afrak_\Rssf \bfrak_\Rssf$
of $\Rssf$.
Set
$\dfrac{\alpha_1}{\alpha_2}
 =
 \alpha
$.
The entire function
$f_S(\omega)$
from $\Pssf_{\Mfrak_S}$
has the property that
a suitable power
$f_S(\omega)^n$,
$n = 1,2,\ldots$
of it
is represented in the form
\[
  f_S(\omega)^n
  \ =\ 
  \dfrac
  {H_1
   \left(\,
   S
   \begin{pmatrix}
     \omega_1
     \\
     \omega_2
   \end{pmatrix}
   \,\right)
  }
  {H_2
   \begin{pmatrix}
     \omega_1
     \\
     \omega_2
   \end{pmatrix}
  }
\]
with two modular forms
$H_1,H_2$
of the same weight~$t$
distinct from~0,
and furthermore
$H_2(\alpha) \neq 0$.
Thus
$f_S(\alpha)$
lies in the field
$\Omega_\Rssf
 =
 \Sigma(j(\alpha))
$.
}

Proof.
From
$j(S_\nu(\alpha))
 =
 j(S(\alpha))
$,
$\nu = 1,2,\ldots,r$,
it follows that
there exist
matrices
$M_\nu \in \Mfrak$
with
\[
  S(\alpha)
  \ =\ 
  M_1\,S_1(\alpha)
  \ =\ 
  \cdots
  \ =\ 
  M_r\,S_r(\alpha)
  .
\]
More precisely,
one has
\[
  S
  \begin{pmatrix}
    \alpha_1
    \\
    \alpha_2
  \end{pmatrix}
  \,
  \xi_\nu
  \ =\ 
  M_\nu
  \,
  S_\nu
  \begin{pmatrix}
    \alpha_1
    \\
    \alpha_2
  \end{pmatrix}
  ,
  \quad
  \nu = 1,2,\ldots,r
\]
with elements
$\xi_\nu \neq 0$
from $\Sigma$.
These $\xi_\nu$
are congruent to
rational integers
modulo the conductor~$f$
of $\Rssf$.
To see that,
let us take
an integral characteristic ideal
$\cfrak_\Rssf$ of $\Rssf$
so that
$\afrak_\Rssf \cfrak_\Rssf$
becomes
a principal ideal
$\afrak_\Rssf \cfrak_\Rssf
 =
 \gamma
 \,
 \Rssf
$.
Let $\varrho,1$ be
a basis of
the principal order
of $\Sigma$;
thus
$f\varrho,1$
is a basis of $\Rssf$.
Hence
one has
\[
  \begin{pmatrix}
  \varrho
  \,
  f
  \,
  \gamma
  \\
  \gamma
  \end{pmatrix}
  \ =\ 
  C
  \begin{pmatrix}
  \alpha_1
  \\
  \alpha_2
  \end{pmatrix}
\]
with
a rational integer matrix $C$,
and it follows that
\[
  \begin{pmatrix}
  \varrho
  \,
  f
  \,
  \xi_\nu
  \\
  \xi_\nu
  \end{pmatrix}
  \ =\ 
  C
  \,
  S^{-1}
  \,
  M_\nu
  \,
  S_\nu
  \,
  C^{-1}
  \begin{pmatrix}
  \varrho\,f
  \\
  1
  \end{pmatrix}
  .
\]
As
$|S|
 =
 \Nm\,\bfrak
$
and
$|C|
 =
 \Nm\,\cfrak
$
are prime to~$f$,
the entries of
the matrix
$C\,S^{-1}\,M_\nu\,S_\nu\,C^{-1}$
are (rational) numbers
with denominators
prime to~$f$,
so that
the second row of
the last equation
reads as
\[
  \xi_\nu
  \ =\ 
  A
  \,
  f
  \,
  \varrho
  \ +\ 
  B
\]
with $f$-integral rationals
$A,B$;
consequently
one has
\[
  \xi_\nu
  \ \equiv\ 
  B
  \mod
  f
  .
\]
The matrices
$S_1,\ldots,S_r$
equivalent to $S$
are thereby
characterized by
$\xi_\nu$ being
a unit.
Since
\[
  S^{-1}
  \,
  M_\nu
  \,
  S_\nu
  \begin{pmatrix}
  \alpha_1
  \\
  \alpha_2
  \end{pmatrix}
  \ =\ 
  \begin{pmatrix}
  \alpha_1
  \\
  \alpha_2
  \end{pmatrix}
  \,
  \xi_\nu
  ,
\]
$S^{-1}
 \,
 M_\nu
 \,
 S_\nu
$
is the representing matrix
of $\xi_\nu$
with respect to
the basis
$\alpha_1,\alpha_2$,
and hence
it follows that
one has
$S^{-1}
 \,
 M_\nu
 \,
 S_\nu
 \in
 \Mfrak
$,
i.e.~$S$ is equivalent to
$S_\nu$
if and only if
$\xi_\nu$
is a unit in $\Rssf$;
but when
$\xi_\nu$
is a unit,
it certainly lies in $\Rssf$,
because
it is congruent to
a rational number
modulo~$f$.

Let $\sigma$ be
an automorphism of
the normal field of
$f_S(\alpha)$
over
\[
  \Pssf(j(\alpha),j(S(\alpha)))
  \ =\ 
  \Omega_\Rssf
  ;
\]
thus
$f_S(\alpha)
 \,
 \sigma
 =
 f_{S_\nu}(\alpha)
$
for a suitable
$\nu = 1,2,\ldots,r$.
One has
\[
  f_S(\alpha)^n
  \,
  \sigma
  \ =\ 
  f_{S_\nu}(\alpha)^n
  \ =\ 
  \dfrac
  {H_1
   \left(\,
   S_\nu
   \begin{pmatrix}
     \alpha_1
     \\
     \alpha_2
   \end{pmatrix}
   \,\right)
  }
  {H_2
   \begin{pmatrix}
     \alpha_1
     \\
     \alpha_2
   \end{pmatrix}
  }
  \ =\ 
  \dfrac
  {H_1
   \left(\,
   M_\nu^{-1}
   \,
   S
   \begin{pmatrix}
     \alpha_1
     \\
     \alpha_2
   \end{pmatrix}
   \,
   \xi_\nu
   \,\right)
  }
  {H_2
   \begin{pmatrix}
     \alpha_1
     \\
     \alpha_2
   \end{pmatrix}
  }
  \ =\ 
  \xi_\nu^t
  \,
  f_S(\alpha)^n
  .
\]
If $\sigma$ has order~$N$,
then one has
$\xi_\nu^{tN} = 1$,
and as
$t \neq 0$,
it follows that
$\xi_\nu$
is a root of unity,
whence
$S_\nu$ is equivalent to $S$,
so
$f_S(\alpha)
 \,
 \sigma
 =
 f_{S_\nu}(\alpha)
 =
 f_S(\alpha)
$,
but that means that
$f_S(\alpha)$
lies in $\Omega_\Rssf$.


\paragraph{21. The principal ideal theorem
for imaginary quadratic number fields.}

We will need
the functions
\[
  \varphi_S(\omega)
  \ =\ 
  |S|^{12}
  \,
  \dfrac
  {\Delta
   \left(\,
   S
   \begin{pmatrix}
     \omega_1
     \\
     \omega_2
   \end{pmatrix}
   \,\right)
  }
  {\Delta
   \begin{pmatrix}
     \omega_1
     \\
     \omega_2
   \end{pmatrix}
  }
\]
for arbitrary rational matrices
with $|S|>0$.
That is not something
essentially new,
for $S$ is uniquely
representable as
$S = r\,S_0$
with primitive $S_0$
and
a positive rational number~$r$,
and one has
\[
  \varphi_S(\omega)
  \ =\ 
  r^{12}
  \,
  \varphi_{S_0}(\omega)
  .
\]
Accordingly,
we set
\[
  \sqrt[24]{\varphi_S(\omega)}
  \ =\ 
  \sqrt{r}
  \,
  \sqrt[24]{\varphi_{S_0}(\omega)}
  \qquad
  \text{with}
  \quad
  \sqrt{r} > 0
  .
\]

Let $\alpha_1,\alpha_2$ be
a basis of
a characteristic ideal
$\afrak_\Rssf$
of the order
$\Rssf$ of $\Sigma$
and let
$S
 \begin{pmatrix}
 \alpha_1
 \\
 \alpha_2
 \end{pmatrix}
$
be a basis of
a characteristic ideal
$\afrak_\Rssf \sfrak_\Rssf$
of $\Rssf$.
As in~\S\textbf{13}
it will be proven that
the divisor of $\varphi_S(\alpha)$
is equal to
$\sfrakbar^{12}$.

From the theorem
proven in~\S\textbf{20},
it follows that
$\varphi_S(\alpha)$
belongs to
the field $\Omega_\Rssf$,
and that
by means of~\S\textbf{4},
the equation for
$\sqrt[24]{\varphi_S(\alpha)}$
holds
if additionally
it is assumed that
$|S|$ is the square of
a rational number
prime to~6.
From that
we conclude that
\footnote{
cf.~\textit{Fricke},
p.~361--362,
also
\textit{Hasse},
Monatsh. f. Math. u. Phys.
\textbf{38},
315--322
(1931).
}

\textit{
Let $\bfrak$ be
a given divisor of $\Sigma$
prime to
the conductor~$f$ of $\Rssf$
and to~6.
Let $\alpha_1,\alpha_2$ be
a basis of
a characteristic ideal
$\afrak_\Rssf$ of $\Rssf$
and let
the rational-entry matrix $B$
transform
$\alpha_1,\alpha_2$
to a basis of
$\afrak_\Rssf \bfrakbar_\Rssf^{-2}$.
Set
$\alpha
 =
 \dfrac{\alpha_1}{\alpha_2}
$;
then
$\sqrt[24]{\varphi_B(\alpha)}$
is an element of
the field $\Omega_\Rssf$
with the divisor $\bfrak$.
}

If the principal order of $\Sigma$
is taken for $\Rssf$
($f=1$),
this is
an explicit version of
the principal ideal theorem
for the field $\Sigma$,
for then
$\Omega_\Rssf$
is the absolute class field
of $\Sigma$.


\paragraph{22. The singular values of
the functions $\varphi_S(\omega)$.
Second part.}

The divisors of
the singular values
$\varphi_S(\alpha)$
can be completely determined.
First
a preliminary remark.
\textit{
Let $\alpha_1,\alpha_2$ be
a basis of
a characteristic ideal
$\afrak_{\Rssf_f}$
of the order $\Rssf_f$
with conductor~$f$;
let $S$ and $S'$ be
two (primitive) matrices
of determinant~$s$
such that
$S
 \begin{pmatrix}
 \alpha_1
 \\
 \alpha_2
 \end{pmatrix}
$
and
$S'
 \begin{pmatrix}
 \alpha_1
 \\
 \alpha_2
 \end{pmatrix}
$
are bases of
characteristic ideals
$\bfrak_{\Rssf_{fs}}$
and
$\bfrak'_{\Rssf_{fs}}$
of the order
$\Rssf_{fs}$
with conductor~$fs$.
Then
the two (algebraic) numbers
$\varphi_S(\alpha)$
and
$\varphi_{S'}(\alpha)$
are associated
(where
$\alpha
 =
 \dfrac{\alpha_1}{\alpha_2}
$).
}

Proof.
In the class of
the ideal
$\bfrak'_{\Rssf_{fs}}
 \bfrak_{\Rssf_{fs}}^{-1}
$,
let $\cfrak_{\Rssf_{fs}}$ be
an integral ideal
which is prime to~$fs$
and to a given
prime number~$p$;
then
$\cfrak_{\Rssf_{fs}}$
is canonically associated to
a divisor
$\cfrak$
prime to~$fsp$.
One has
\[
  \bfrak'_{\Rssf_{fs}}
  \,
  \gamma
  \ =\ 
  \bfrak_{\Rssf_{fs}}
  \,
  \cfrak_{\Rssf_{fs}}
  .
  \tag{48}
  \label{eqn:48}
\]

By~\S\textbf{8}
it follows
by multiplication with
$\Rssf_f$
that
\[
  \begin{aligned}
  \afrak_{\Rssf_f}
  \,
  \gamma
  &
  \ =\ 
  \afrak_{\Rssf_f}
  \,
  \cfrak_{\Rssf_f}
  ,
  \\
  \gamma
  \,
  \Rssf_f
  &
  \ =\ 
  \cfrak_{\Rssf_f}
  ,
  \end{aligned}
\]
and by multiplication with
$\Rssf_1$
that
\[
  \gamma
  \,
  \Rssf_1
  \ =\ 
  \cfrak_{\Rssf_1}
  ;
\]
the divisor of $\gamma$
is thus
$\cfrak$.
From~\eqref{eqn:48}
it further follows that
\[
  S'
  \begin{pmatrix}
    \alpha_1
    \\
    \alpha_2
  \end{pmatrix}
  \,
  \gamma
  \ =\ 
  C
  \,
  S
  \begin{pmatrix}
    \alpha_1
    \\
    \alpha_2
  \end{pmatrix}
\]
for a integer-entry matrix $C$
whose determinant
\[
  |C|
  \ =\ 
  [\,
  \bfrak_{\Rssf_{fs}}
  \,{:}\,
  \bfrak_{\Rssf_{fs}}
  \,
  \cfrak_{\Rssf_{fs}}
  \,]
\]
is equal to
the norm
$\Nm\,\cfrak
 =
 \Nm\,\gamma
$.
Hence
one has
\[
  \begin{aligned}
  \dfrac
  {\varphi_{S'}(\alpha)}
  {\varphi_S(\alpha)}
  &
  \ =\ 
  \dfrac
  {\Delta
   \left(\,
   S'
   \begin{pmatrix}
     \alpha_1
     \\
     \alpha_2
   \end{pmatrix}
   \,\right)
  }
  {\Delta
   \left(\,
   S
   \begin{pmatrix}
     \alpha_1
     \\
     \alpha_2
   \end{pmatrix}
   \,\right)
  }
  \ =\ 
  \dfrac
  {\Delta
   \left(\,
   C
   \,
   S
   \begin{pmatrix}
     \alpha_1
     \\
     \alpha_2
   \end{pmatrix}
   \,
   \gamma^{-1}
   \,\right)
  }
  {\Delta
   \left(\,
   S
   \begin{pmatrix}
     \alpha_1
     \\
     \alpha_2
   \end{pmatrix}
   \,\right)
  }
  \ =\ 
  \gamma^{12}
  \,
  \dfrac
  {\Delta
   \left(\,
   C
   \,
   S
   \begin{pmatrix}
     \alpha_1
     \\
     \alpha_2
   \end{pmatrix}
   \,\right)
  }
  {\Delta
   \left(\,
   S
   \begin{pmatrix}
     \alpha_1
     \\
     \alpha_2
   \end{pmatrix}
   \,\right)
  }
  \\[2ex]
  &
  \ =\ 
  \gamma^{12}
  \,
  \Nm\,\cfrak^{-12}
  \,
  \varphi_C(S(\alpha))
  \ =\ 
  \gammabar^{-12}
  \,
  \varphi_C(S(\alpha))
  ,
  \end{aligned}
\]
and as
$\varphi_C(S(\alpha))$
is a factor of
a power of
$\Nm\,\cfrak = \Nm\,\gamma$
by~\eqref{eqn:17},
it follows that
$\dfrac
 {\varphi_{S'}(\alpha)}
 {\varphi_S(\alpha)}
$
is prime to
the arbitrarily specified
prime number~$p$,
whence
a unit.

Since
a primitive matrix $S$
of determinant~$|S|$
can be put in the form
\[
  S
  \ =\ 
  P_1
  P_2
  \cdot
  \cdots
  \cdot
  P_r
  ,
\]
where each
$|P_\nu|$ is a prime number,
by the formula
$\varphi_{AB}(\alpha)
 =
 \varphi_A(B(\alpha))
 \varphi_B(\alpha)
$,
in order to compute
the divisor of
a singular value
$\varphi_S(\alpha)$,
it suffices to know
the divisors of the
$\varphi_P(\alpha)$
where
$|P|$ is a prime number.

\textit{
Now let $p$ be
a prime number,
let $\Rssf_f$ be
the order
with conductor~$f$,
and let $p^t$ be
the power of~$p$
contained in $f$.
Let $P$ denote
the primitive matrix
of determinant~$p$.
Let $\alpha_1,\alpha_2$ be
a basis of
a characteristic ideal
$\afrak_{\Rssf_f}$
of $\Rssf_f$.
Then
the divisors of
the singular value
$\varphi_P(\alpha)$
are given as follows:
\footnote{
Essentially in
\textit{Hasse},
Monatsh. f. Math.
\textbf{38},
p.~331
(1931).
}
\begin{itemize}
\item[1.]
If $p$ is
the product of
two distinct prime divisors
$\pfrak$ and $\pfrakbar$
in $\Sigma$:
\begin{itemize}
\item[1.1.]
$\varphi_P(\alpha)$
is a unit
when
$P
 \begin{pmatrix}
 \alpha_1
 \\
 \alpha_2
 \end{pmatrix}
$
is a basis of
a characteristic ideal of
$\Rssf_{fp}$,
\item[1.2.]
$\varphi_P(\alpha)
 \approx
 p^{12}
$
when
$P
 \begin{pmatrix}
 \alpha_1
 \\
 \alpha_2
 \end{pmatrix}
$
is a basis of
a characteristic ideal of
$\Rssf_{fp^{-1}}$,
\item[1.3.]
$\varphi_{P_\pfrak}(\alpha)
 \approx
 \pfrakbar^{12}
$
and
$\varphi_{P_\pfrakbar}(\alpha)
 \approx
 \pfrak^{12}
$
when
$p \nmid f$
($t=0$)
and
$P_\pfrak
 \begin{pmatrix}
 \alpha_1
 \\
 \alpha_2
 \end{pmatrix}
$
is a basis of
$\afrak_{\Rssf_f}
 \,
 \pfrak_{\Rssf_f}
$
and
$P_\pfrakbar
 \begin{pmatrix}
 \alpha_1
 \\
 \alpha_2
 \end{pmatrix}
$
is a basis of
$\afrak_{\Rssf_f}
 \,
 \pfrakbar_{\Rssf_f}
$.
\end{itemize}
\item[2.]
If $p$ is
the square of
a prime divisor
$\pfrak$
in $\Sigma$:
\begin{itemize}
\item[2.1.]
$\varphi_P(\alpha)
 \approx
 p^{\textstyle \frac{6}{p^{t+1}}}
$
when
$P
 \begin{pmatrix}
 \alpha_1
 \\
 \alpha_2
 \end{pmatrix}
$
is a basis of
a characteristic ideal of
$\Rssf_{fp}$.
\item[2.2.]
$\varphi_P(\alpha)
 \approx
 p^{12
    -
    \textstyle \frac{6}{p^t}}
$
when
$P
 \begin{pmatrix}
 \alpha_1
 \\
 \alpha_2
 \end{pmatrix}
$
is a basis of
a characteristic ideal of
$\Rssf_{fp^{-1}}$.
\item[2.3.]
$\varphi_{P_\pfrak}(\alpha)
 \approx
 p^6
$
when
$P_\pfrak
 \begin{pmatrix}
 \alpha_1
 \\
 \alpha_2
 \end{pmatrix}
$
is a basis of
$\afrak_{\Rssf_f}
 \,
 \pfrak_{\Rssf_f}
$.
\end{itemize}
\item[3.]
If $p$ is prime
in $\Sigma$:
\begin{itemize}
\item[3.1.]
$\varphi_P(\alpha)
 \approx
 p^{\textstyle \frac{12}{p^t(p+1)}}
$
when
$P
 \begin{pmatrix}
 \alpha_1
 \\
 \alpha_2
 \end{pmatrix}
$
is a basis of
a characteristic ideal of
$\Rssf_{fp}$.
\item[3.2.]
$\varphi_P(\alpha)
 \approx
 p^{12
    \textstyle 
    \left[
    1
    -
    \frac{1}{p^{t-1}(p+1)}
    \right]}
$
when
$P
 \begin{pmatrix}
 \alpha_1
 \\
 \alpha_2
 \end{pmatrix}
$
is a basis of
a characteristic ideal of
$\Rssf_{fp^{-1}}$.
\end{itemize}
\end{itemize}
}

We prove this
by induction on~$t$.
For $t=0$
the assertions of~1.3 and~2.3
are proven in~\S\textbf{13};
for~3.1
the assertion follows
immediately
from the preliminary remark,
for then
\emph{all}
$\begin{pmatrix}
 \alpha_1
 \\
 \alpha_2
 \end{pmatrix}
$
are bases of
characteristic ideals in
$\Rssf_{fp}$,
thus
letting
$P_1,\ldots,P_{p+1}$ be
a system of representatives of
the equivalence classes of
matrices
of determinant~$p$,
one has
\[
  \varphi_{P_\nu}(\alpha)^{p+1}
  \ \approx\ 
  \prod_{\nu=1}^{p+1}
  \varphi_{P_\nu}(\alpha)
  \ =\ 
  p^{12}
  ,
\]
and we conclude
just as
in the cases of~1.1 and~2.1.

Now let $t>0$.
Then by~\S\textbf{8}
there exists
exactly one of
the representatives~$P_\nu$,
say $P_1$,
for which
$P_1
 \begin{pmatrix}
 \alpha_1
 \\
 \alpha_2
 \end{pmatrix}
$
is a basis of
a characteristic ideal of
$\Rssf_{fp^{-1}}$.
Since
$P_\nu
 \begin{pmatrix}
 \alpha_1
 \\
 \alpha_2
 \end{pmatrix}
$
for the remaining
$P_\nu$
are bases of
characteristic ideals of
$\Rssf_{fp}$,
it follows that
one has
\[
  p^{12}
  \ =\ 
  \prod_{\nu=1}^{p+1}
  \varphi_{P_\nu}(\alpha)
  \ \approx\ 
  \varphi_{P_1}(\alpha)
  \,
  \varphi_{P_\nu}(\alpha)^p
  ,
  \qquad
  \nu > 1,
\]
so it suffices to determine
the divisor of
$\varphi_{P_1}(\alpha)$.
Now one has
\[
  \varphi_{p\,P_1^{-1}}(P_1(\alpha))
  \,
  \varphi_{P_1}(\alpha)
  \ =\ 
  \varphi_{p\,E}(\alpha)
  \ =\ 
  p^{12}
  ,
\]
and
$P_1
 \begin{pmatrix}
 \alpha_1
 \\
 \alpha_2
 \end{pmatrix}
$
is a basis of
a characteristic ideal of
$\Rssf_{fp^{-1}}$,
which is transformed
by the matrix
$p\,P_1^{-1}$
of determinant~$p$
to the basis
$\alpha_1,\alpha_2$
of the characteristic ideal
$\afrak_{\Rssf_f}$
of $\Rssf_f$,
and hence
by the induction hypothesis
one has
$\varphi_{p\,P_1^{-1}}(\alpha)
 \approx
 p^{\textstyle \frac{12}{p^{t-1}(p+1)}}
$
and
$\varphi_{P_1}(\alpha)
 \approx
 p^{12
    \textstyle 
    \left[
    1
    -
    \frac{1}{p^{t-1}(p+1)}
    \right]}
$,
q.e.d.


\paragraph{23. Congruences for
the singular values of functions from $\Pssf_{\Mfrak_S}$.}

The theorem proven in~\S\textbf{12}
on the divisibility of
the singular value of
functions from
$\Pssf_{\Mfrak_P}$,
$|P|=p$ prime,
by a $\Sigma$-prime factor
$\pfrak$ of $p$,
is for the case
$p = \pfrak\,\pfrakbar$,
$\pfrak \neq \pfrakbar$.
Let us now
establish
a similar theorem
for the cases of
$p = \pfrak$
and
$p = \pfrak^2$
put aside.

Let $p$ be
a prime number,
and let
$P_1,\ldots,P_{p+1}$ be
a system of representatives of
the equivalence classes of
matrices
of determinant~$p$;
let $P$ stand for
an arbitrary matrix
of determinant~$p$.
Let $f_P(\omega)$ be
an entire integral function
from $\Pssf_{\Mfrak_P}$.
By~\S\textbf{3}
the principal polynomial
\[
  F(X,j(\omega))
  \ =\ 
  \prod_{\nu=1}^{p+1}
  (X-f_{P_\nu}(\omega))
  \ =\ 
  \sum_{\mu=0}^{p}
  Q_\mu(j(\omega))
  \,
  X^\mu
\]
of $f_P(\omega)$
over $\Pssf_{\Mfrak}$
has rational integer coefficients
as a polynomial
in $X$ and $j(\omega)$;
thus
$f_P(\alpha)$
is an integer
for imaginary quadratic $\alpha$.

\textit{
The entire integral function
$f_P(\omega)$
from $\Pssf_{\Mfrak_P}$
has the property that
the $q$-coefficients of
the conjugate
$f_{\left(\begin{smallmatrix}
          p & 0
          \\
          0 & 1
          \end{smallmatrix}\right)}
 (\omega)
$
are divisible by~$p$.
Let $\alpha$ be
imaginary quadratic,
$\Pssf(\alpha) = \Sigma$,
and let
$\alpha
 =
 \dfrac{\alpha_1}{\alpha_2}
$
where
$\alpha_1,\alpha_2$ is
a basis of
a characteristic ideal
$\afrak_{\Rssf_f}$
of the order
$\Rssf_f$
with conductor~$f$
in $\Sigma$.
Let $p$ be
prime or ramified
in $\Sigma$;
thus
$p = \pfrak$
or
$p = \pfrak^2$
where $\pfrak$ is
a prime divisor of $\Sigma$.
Let $\Pfrak$ be
a prime factor of $p$
in the number field
$\Lambda_P
 =
 \Pssf(j(\alpha),f_P(\alpha))
$.
Then
in the following cases,
$f_P(\alpha)$
is divisible by $\Pfrak$:
\begin{itemize}
\item[1.]
$p = \pfrak$,
$f \not\equiv 0
 \mod
 p
$,
$P$ arbitrary.
\item[2.]
$p = \pfrak^2$,
$f \not\equiv 0
 \mod
 p
$,
$P
 \begin{pmatrix}
 \alpha_1
 \\
 \alpha_2
 \end{pmatrix}
$
is a basis of
$\afrak_{\Rssf_f}
 \pfrak_{\Rssf_f}
$.
\item[3.]
$f \equiv 0
 \mod
 p
$,
$P
 \begin{pmatrix}
 \alpha_1
 \\
 \alpha_2
 \end{pmatrix}
$
is a basis of
a characteristic ideal of
$\Rssf_{fp^{-1}}$.
\end{itemize}
But we also make
these two assumptions
on $p$:
\ 
$p{\,>\,}12$
and
$\pfrak$ is unramified
in $\Omega_{\Rssf_{f_0}}$
where
$f = f_0\,p^t$
with
$f_0 \not\equiv 0
 \mod
 p
$.
}

Proof.
By~\S\textbf{3}
one has
the congruence
in the $q$-coefficients
\[
  F(X,j(\omega))
  \ \equiv\ 
  \bigl(\,
  X
  -
  f_{\left(\begin{smallmatrix}
          p & 0
          \\
          0 & 1
          \end{smallmatrix}\right)}
  (\omega)
  \,\bigr)
  \,
  \bigl(\,
  X^p
  -
  f_{\left(\begin{smallmatrix}
          1 & 0
          \\
          0 & p
          \end{smallmatrix}\right)}
  (\omega)^p
  \,\bigr)
  \mod
  (1-\zeta_p)
  ;
\]
thus
by the assumption on
$f_{\left(\begin{smallmatrix}
          p & 0
          \\
          0 & 1
          \end{smallmatrix}\right)}
 (\omega)
$,
one even has
\[
  F(X,j(\omega))
  \ \equiv\ 
  X
  \,
  \bigl(\,
  X^p
  -
  f_{\left(\begin{smallmatrix}
          1 & 0
          \\
          0 & p
          \end{smallmatrix}\right)}
  (\omega)^p
  \,\bigr)
  \mod
  (1-\zeta)
  ,
\]
and
from that,
as in~\S\textbf{3},
it follows that
one has
the congruence
\[
  F(X,j(\omega))
  \ \equiv\ 
  X
  \,
  (X^p-Q_1(j(\omega)))
  \mod
  p
\]
for $F$ as a polynomial
in $X$ and $j(\omega)$,
and thus also
\[
  F(X,j(\alpha))
  \ \equiv\ 
  X
  \,
  (X^p-Q_1(j(\alpha)))
  \mod
  p
  .
  \tag{49}
  \label{eqn:49}
\]
One has
$\Pssf(j(\alpha))
 = 
 \Omega_{\Rssf_f}
$,
$\Lambda_P
 =
 \Omega_{\Rssf_f}(f_P(\alpha))
$
and
$\Lambda^*
 =
 \Omega_{\Rssf_f}
 (f_{P_1}(\alpha),
  \ldots,
  f_{P_{p+1}}(\alpha)
 )
$
is an extension field of $\Lambda$
normal over
$\Omega_{\Rssf_f}$.
Let
$\Pfrak^*$
be a prime factor of
$\Pfrak$
in $\Lambda^*$.
From~\eqref{eqn:49}
it follows that
at least one of the
$p{+}1$~elements
$f_{P_\nu}(\alpha)$
is divisible by
$\Pfrak^*$.
Let us say
$f_{P_1}(\alpha)
 \equiv
 0
 \mod
 \Pfrak^*
$;
then
one has
$f_{P_1}(\alpha)^p
 \equiv
 Q_1(j(\alpha))
 \mod
 \Pfrak^*
$
for $\nu \geqslant 2$.

If
$f_{P_1}(\alpha)$
is a multiple zero of
$F(X,j(\alpha))$,
then one has
\[
  Q_1(j(\alpha))
  \ \equiv\ 
  0
  \mod
  \Pfrak^*
  ,
  \qquad
  f_{P_\nu}(\alpha)
  \ \equiv\ 
  0
  \mod
  \Pfrak^*
\]
for all $\nu$
and in particular
$f_P(\alpha)
 \equiv
 0
 \mod
 \Pfrak
$.

Thus let us assume that
$f_{P_1}(\alpha)$
is a simple zero of
$F(X,j(\alpha))$.
Then
$f_{P_1}(\omega)$
is all the more
a simple zero of
$F(X,j(\omega))$
and consequently
$f_P(\omega)$
is a generating element
of $\Pssf_{\Mfrak_P}$
over $\Pssf_{\Mfrak}$.
From this
we conclude
as in~\S\textbf{3}
that
\[
  F'
  \bigl(\,
  f_P(\omega)
  \,,\,
  j(\omega)
  \,\bigr)
  \,
  \varphi_P(\omega)
  \ =\ 
  \sum_{\mu=0}^{p}
  b_\mu(j(\omega))
  \,
  f_P(\omega)^\mu
\]
with polynomials
$b_\mu$,
whose coefficients are
rational numbers
(even rational integers).
As one has
$F'(f_{P_1}(\alpha),j(\alpha))
 \neq
 0
$,
the (algebraic) number
$\varphi_{P_1}(\alpha)$
lies in the field
$\Lambda_{P_1}$.

Let us first consider
the case
$p = \pfrak$,
$f \not\equiv 0
 \mod
 p
$.
Then one has
$\varphi_P(\alpha)
 \approx
 \pfrak^{\textstyle \frac{12}{p+1}}
$
and because
$\pfrak$ is unramified
in $\Omega_{\Rssf_f}$
and $p > 12$,
it follows that
$\Pfrak$ is ramified
over $\Omega_{\Rssf_f}$.
The inertia group of
$\Pfrak^*$ over $\Omega_{\Rssf_f}$
is consequently
not contained in
the Galois group of
$\Lambda^*/\Lambda_{P_1}$,
and so
there exists
an inertia automorphism
$\lambda$
of $\Pfrak^*/\Omega_{\Rssf_f}$
with
$f_{P_1}(\alpha)
 \,
 \lambda
 \neq
 f_{P_1}(\alpha)
$.

Thus
$f_{P_1}(\alpha)
 \,
 \lambda
 =
 f_{P_\nu}(\alpha)
$
for a suitable
$\nu \geqslant 2$
is divisible by $\Pfrak^*$,
\[
  Q_1(j(\alpha))
  \ \equiv\ 
  0
  \mod
  \Pfrak^*
  ,
  \qquad
  f_{P_\nu}(\alpha)
  \ \equiv\ 
  0
  \mod
  \Pfrak^*
\]
for all $\nu$
and in particular
$f_P(\alpha)
 \equiv
 0
 \mod
 \Pfrak
$.

In the two other cases,
$p = \pfrak$
but
$f \equiv 0
 \mod
 p
$
and
$p = \pfrak^2$,
the matrix $P$
to be examined
is equivalent to
a certain $P_\nu$.
If $P \sim P_1$,
then
the assertion holds.
If $P$ is
not equivalent to $P_1$,
then
$P_1
 \begin{pmatrix}
 \alpha_1
 \\
 \alpha_2
 \end{pmatrix}
$
is a basis of
a characteristic ideal of
$\Rssf_{fp}$
and hence
$\varphi_{P_1}(\alpha)
 \approx
 p^{\textstyle \frac{12}{p^t(p+1)}}
$
in the case
$p = \pfrak$
and
$\varphi_{P_1}(\alpha)
 \approx
 p^{\textstyle \frac{12}{p^{t+1}}}
$
in the case
$p = \pfrak^2$.
But the ramification index
of $\pfrak$ in $\Omega_{\Rssf_f}$
is a factor of
the degree
\[
  [\,
  \Omega_{\Rssf_f}
  \,{:}\,
  \Omega_{\Rssf_{f_0}}
  \,]
  \ =\ 
  h_{\Rssf_f}
  \,/\,
  h_{\Rssf_{f_0}}
  \ =\ 
  \dfrac
  {e_{\Rssf_f}}
  {e_{\Rssf_{f_0}}}
  \,
  \begin{cases}
  p^{t-1}(p{+}1)
  ,
  &
  \quad
  p = \pfrak
  \\
  p^t
  ,
  &
  \quad
  p = \pfrak^2
  \end{cases}
\]
(because
$\pfrak$ is unramified
in $\Omega_{\Rssf_{f_0}}$).
In the course of this
we have used
the correspondence theorem
from~\S\textbf{19}.
We have to assume
the formula for
$h_{\Rssf_f}
 /
 h_{\Rssf_{f_0}}
$
as well-known.
Since $p{\,>\,}12$,
it follows that
$\Pfrak$
is ramified over
$\Omega_{\Rssf_f}$
and as above
we conclude that
$f_{P_\nu}(\alpha)
 \equiv
 0
 \mod
 \Pfrak^*
$
for every~$\nu$,
in particular
$f_P(\alpha)
 \equiv
 0
 \mod
 \Pfrak
$.


\section*{D. The second main theorem}
\addcontentsline{toc}{section}{\hfill
D. The second main theorem
\hfill}


\paragraph{24. The singular elliptic functions.}

Let $K(\wfrak)$ denote
the field of
elliptic functions
with the period lattice
$\wfrak$.
The field
$K(\wfrak)$
is called
\emph{singular}
if $\wfrak = \afrak_\Rssf$
is an ideal
in the order
$\Rssf$
of an imaginary quadratic
number field
$\Sigma$.

Every unit
$\varepsilon$
of $\Rssf$
generates
an automorphism
$\varrho_\varepsilon$
of the singular field
$K(\afrak_\Rssf)$:
For a function
$f(z)$
from $K(\afrak_\Rssf)$
one sets
$f(z)
 \,
 \varrho_\varepsilon
 =
 f(\varepsilon\,z)
$.
The map
$\varepsilon
 \to
 \varrho_\varepsilon
$
is an isomorphism
from the unit group of $\Rssf$
onto
an automorphism group
of order~$e$
(number of units of $\Rssf$)
of $K(\afrak_\Rssf)$,
because
one has
$\varrho_{\varepsilon_1}
 \cdot
 \varrho_{\varepsilon_2}
 =
 \varrho_{\varepsilon_1 \varepsilon_2}
$
and
since
\[
  \wp(z,\afrak_\Rssf)
  \,
  \wp'(z,\afrak_\Rssf)^{-1}
  \cdot
  \varrho_\varepsilon
  \ =\ 
  \varepsilon
  \ 
  \wp(z,\afrak_\Rssf)
  \,
  \wp'(z,\afrak_\Rssf)^{-1}
\]
one has
$\varrho_\varepsilon = 1$
only for
$\varepsilon = 1$.

The
\emph{singular Weber function}
$\tau_\Rssf(z,\afrak_\Rssf)$
is invariant under the
$\varrho_\varepsilon$,
thus also under
the transformations
\[
  z
  \to
  \varepsilon\,z + \alpha
  ,
  \qquad
  \text{$\varepsilon$ unit in $\Rssf$}
  ,
  \ 
  \text{$\alpha$ in $\afrak_\Rssf$}
  .
\]
In the group of
these transformations,
the group of
the translations
\[
  z
  \to
  z + \alpha
  ,
  \qquad
  \text{$\alpha$ in $\afrak_\Rssf$}
\]
has index~$e$,
and since
$\tau_\Rssf(z,\afrak_\Rssf)$
as an elliptic function
has order~$e$,
\textit{
one has
$\tau_\Rssf(z',\afrak_\Rssf)
 =
 \tau_\Rssf(z,\afrak_\Rssf)
$
if and only if
one has
$z'
 \equiv
 \varepsilon
 \,
 z
 \mod
 \afrak_\Rssf
$,
where
$\varepsilon$
is a unit in $\Rssf$.
}
The field
$K(\afrak_\Rssf)$
has degree~$e$
over the field of
rational functions in
$\tau_\Rssf(z,\afrak_\Rssf)$,
which we will denote by
$K_0(\afrak_\Rssf)$,
and hence
\textit{
$K_0(\afrak_\Rssf)$
is the fixed field
of the
$\varrho_\varepsilon$.
}

In the power series
\[
  \wp(z,\wfrak)
  \ =\ 
  \sum_{\nu=-2}^{\infty}
  a_\nu(\wfrak)
  \,
  z^\nu
\]
the coefficient
$a_\nu(\wfrak)$
is obviously
an entire modular form
of weight~$-\nu$,
and as
$a_\nu(\omega)$
is known to be
a polynomial
in $g_2$ and $g_3$
with rational coefficients,
the Fourier coefficients
of $a_\nu(\wfrak)$
are rational.
It follows that
\[
  \tau_\Rssf(z,\wfrak)
  \ =\ 
  \sum_{\nu=-e}^{\infty}
  b_\nu(\wfrak)
  \,
  z^\nu
  ,
\]
where
$b_\nu(\wfrak)$
is an entire modular form
of weight~$-\nu$
with rational Fourier coefficients.
The \emph{singular} function
$\tau_\Rssf(z,\afrak_\Rssf)$
depends only on
$z^e$,
i.e.~one has
\[
  \tau_\Rssf(z,\afrak_\Rssf)
  \ =\ 
  \sum_{\nu=-1}^{\infty}
  b_{e\nu}(\afrak_\Rssf)
  \,
  z^{e\nu}
  .
\]
We introduce
$u = z^e\,g^{(e)}(\afrak_\Rssf)^{-1}$
as a new variable:
\[
  \tau_\Rssf(z,\afrak_\Rssf)
  \ =\ 
  \sum_{\nu=-1}^{\infty}
  b_{e\nu}(\afrak_\Rssf)
  \,
  g^{(e)}(\afrak_\Rssf)^\nu
  \,
  u^\nu
  .
\]
Here
$b_{e\nu}(\wfrak)
 \,
 g^{(e)}(\wfrak)^\nu
$
is an entire modular \emph{function}
with rational Fourier coefficients,
thus
a polynomial in $j(\wfrak)$
with rational coefficients,
so
$b_{e\nu}(\afrak_\Rssf)
 \,
 g^{(e)}(\afrak_\Rssf)^\nu
$
lies in $\Omega_\Rssf$.
\textit{
We thus have
\[
  \tau_\Rssf(z,\afrak_\Rssf)
  \ =\ 
  \sum_{\nu=-1}^{\infty}
  \eta_\nu
  \,
  u^\nu
  ,
  \qquad
  \text{$\eta_\nu$ in $\Omega_\Rssf$}
  ,
  \quad
  u = z^e\,g^{(e)}(\afrak_\Rssf)^{-1}
  .
  \tag{50}
\]
}


\paragraph{25. The singular values of the Weber function.
Ray class invariants.}

In the following,
let $\Sigma$ be
a fixed imaginary quadratic
number field.
By a \emph{singular value of
the Weber function
$\tau_\Rssf(z,\wfrak)$
associated to
an order $\Rssf$ of $\Sigma$},
we mean
a function value
$\tau_\Rssf(\gamma,\afrak_\Rssf)$,
where $\gamma$ is
an element of $\Sigma$
and $\afrak_\Rssf$
is a characteristic ideal of $\Rssf$.
For this,
we assume that
$\gamma
 \not\equiv
 0
 \mod
 \afrak_\Rssf
$,
so that
one has
$\tau_\Rssf(\gamma,\afrak_\Rssf)
 \neq
 \infty
$.
\textit{
The singular values
$\tau_\Rssf(\gamma,\afrak_\Rssf)$
are algebraic numbers,
}
because
taking
$\alpha_1,\alpha_2$
to be a basis of $\afrak_\Rssf$,
one has
$\gamma
 =
 N^{-1}(x_1\alpha_1 + x_2\alpha_2)
$
with rational integers
$N,x_1,x_2$,
$(N,x_1,x_2)=1$,
and hence
$\tau_\Rssf(\gamma,\afrak_\Rssf)$
is a zero of
the division polynomial
$T_N(X,j(\alpha))$,
whose coefficients
lie in $\Omega_\Rssf$.
The denominator (-divisor) of
$\tau_\Rssf(1,\afrak_\Rssf)$
contains at most
prime factors of $N$
and
$\tau_\Rssf(1,\afrak_\Rssf)$
is even integral
when $N$ is not
a prime power.

The main result on
the singular values of
$\tau_\Rssf(z,\wfrak)$
is:
\textit{
the number field
$\Omega_\Rssf(\tau_\Rssf(\gamma,\afrak_\Rssf))$
is the class field of
$\Sigma$
for a divisor class group
determined
in a simple way
by $\afrak_\Rssf$
and $\gamma$.
}
We will prove this
only for the case when
$\Rssf$ is
the principal order of
$\Sigma$
(second main theorem);
the general case
can be treated
quite similarly,
but it requires
some complicated
considerations
from the ideal theory of
non-principal orders,
which unnecessarily
encumbers
the simple train of thought.

Thus let $\Rssf$
from now on
be the principal order
of $\Sigma$.
Aside from that,
only the orders
$\Rssf_p$
with prime conductor~$p$
play a role.

Let $\mfrak$ be
an integral divisor of $\Sigma$
distinct from~$1$.
The singular value
$\tau_\Rssf(\gamma,\afrak_\Rssf)$
is called
an \emph{$\mfrak$-th division value
of $\afrak$}
if one has
$\mfrak\,\gamma
 \equiv
 0
 \mod
 \afrak
$,
i.e.~if
$\gamma
 \equiv
 0
 \mod
 \afrak\,\mfrak^{-1}
$.
The value
$\tau_\Rssf(\gamma,\afrak_\Rssf)$
is called
a \emph{characteristic
$\mfrak$-th division value
of $\afrak$}
or
a \emph{division value
of order $\mfrak$
of $\afrak$}
if $\mfrak$ is
the smallest
(in the sense of
divisibility)
integral divisor
with
$\mfrak\,\gamma
 \equiv
 0
 \mod
 \afrak
$;
in other words:
if
$(\gamma)
 =
 \afrak
 \,
 \rfrak
 \,
 \mfrak^{-1}
$
with integral
$\rfrak$
prime to $\mfrak$,
or also:
if $\mfrak$ is
the denominator of
the divisors
$\gamma\,\afrak^{-1}$.
The number $N$
in the representation
\[
  x
  \ =\ 
  N^{-1}
  (x_1\alpha_1 + x_2\alpha_2)
\]
with relatively prime
rational integers
$x_1,x_2,N$,
$N>0$,
is the smallest
natural number
divisible by
the order $\mfrak$
of $\tau_\Rssf(\gamma,\afrak_\Rssf)$;
since
$\Nm\,\mfrak \mid N^2$
and
$N \mid \Nm\,\mfrak$,
it follows that
the denominator of
$\tau_\Rssf(\gamma,\afrak_\Rssf)$
is at most divisible by
the prime factors of
$\Nm\,\mfrak$,
and
$\tau_\Rssf(\gamma,\afrak_\Rssf)$
is even integral
when
$\Nm\,\mfrak$
is not a prime power.

The value
$\tau_\Rssf(\gamma,\afrak_\Rssf)$
depends only on
$\gamma$ modulo $\afrak_\Rssf$,
and one has
\[
  \tau_\Rssf(\gamma,\afrak_\Rssf)
  \ =\ 
  \tau_\Rssf(\gamma',\afrak_\Rssf)
\]
if and only if
one has
$\gamma'
 \equiv
 \gamma
 \,
 \varepsilon
 \mod
 \afrak_\Rssf
$
with a unit
$\varepsilon$
of $\Rssf$.

If
$\tau_\Rssf(\gamma,\afrak_\Rssf)$
is the $\mfrak$-th
division value of
$\afrak_\Rssf$,
then
one has
$(\gamma)
 =
 \afrak
 \,
 \mfrak^{-1}
 \,
 \rfrak
$
with integral
$\rfrak$,
and
$\rfrak^{-1}$
lies in the (absolute)
class of the divisor
$\afrak\,\mfrak^{-1}$.
Then
one has
\[
  \tau_\Rssf(\gamma,\afrak_\Rssf)
  \ =\ 
  \tau_\Rssf(1,\mfrak_\Rssf\rfrak_\Rssf^{-1})
  .
\]
Conversely,
if $\rfrak$ is
an integral divisor,
with
$\rfrak^{-1}$
in the class of
$\afrak\,\mfrak^{-1}$,
then
one has
\[
  \tau_\Rssf(1,\mfrak_\Rssf\rfrak_\Rssf^{-1})
  \ =\ 
  \tau_\Rssf(\gamma,\afrak_\Rssf)
\]
with
$(\gamma)
 =
 \afrak
 \,
 \mfrak^{-1}
 \,
 \rfrak
$,
so
$\tau_\Rssf(1,\mfrak_\Rssf\,\rfrak_\Rssf^{-1})$
is an $\mfrak$-th division value
of $\afrak_\Rssf$.
\textit{
We will henceforth
regard
$\tau_\Rssf(1,\mfrak_\Rssf\rfrak_\Rssf^{-1})$
also as
an $\mfrak$-th division value
of the divisor class $\kfrak$
of $\afrak$,
i.e.~of $\mfrak\,\rfrak^{-1}$.
}

The value
$\tau_\Rssf(1,\mfrak_\Rssf\,\rfrak_\Rssf^{-1})$
is a \emph{characteristic}
$\mfrak$-th division value
of the class $\kfrak$
of $\mfrak\,\rfrak^{-1}$
(i.e.~of $\mfrak\,\rfrak^{-1}$)
if and only if
$\rfrak$ is
prime to $\mfrak$.

\textit{
The value
$\tau_\Rssf(1,\mfrak_\Rssf\,\rfrak_\Rssf^{-1})$
depends only on
the ray class
$\kfrak_\mfrak$
to which
$\rfrak^{-1}$
belongs.
}
(The ray classes
of $\Sigma$
modulo $\mfrak$
will in general
be denoted by
$\kfrak_\mfrak,
 \kfrak_\mfrak',
 \kfrak_\mfrak^{(\nu)},
 \hfrak_\mfrak,
 \ldots
$.)
This is because
if
$\rfrak$ and $\rfrak'$
are equivalent
modulo~$\mfrak$,
then
$\rfrak'
 =
 \lambda
 \,
 \rfrak
$
with
$\lambda
 \equiv
 1
 \ \textup{mod}^*\ 
 \mfrak
$
\footnote{
$\alpha \equiv \beta
 \ \textup{mod}^*\ 
 \mfrak
$
means
multiplicative congruence.
},
and since
the denominator
of $\lambda$
divides $\rfrak$,
it follows that
$\lambda
 \equiv
 1
 \mod
 \mfrak\,\rfrak^{-1}
$,
thus
\[
  \tau_\Rssf(1,\mfrak_\Rssf\,\rfrak_\Rssf^{-1})
  \ =\ 
  \tau_\Rssf(\lambda,\mfrak_\Rssf\,\rfrak_\Rssf^{-1})
  \ =\ 
  \tau_\Rssf(1,(\mfrak_\Rssf\,\rfrak_\Rssf^{-1})\,\lambda^{-1})
  \ =\ 
  \tau_\Rssf(1,\mfrak_\Rssf\,{\rfrak_\Rssf'}^{-1})
  .
\]
Henceforth,
when
$\rfrak^{-1}$
lies in the ray class
$\kfrak_\mfrak$
modulo $\mfrak$,
we will call
$\tau_\Rssf(1,\mfrak_\Rssf\,\rfrak_\Rssf^{-1})$
\emph{the invariant of
the ray class
$\kfrak_\mfrak$}
and denote it by
$\tau_\Rssf(\kfrak_\mfrak)$,
analogous to
the class invariant
$j(\kfrak_\Rssf)$.

But unlike
the $j(\kfrak_\Rssf)$,
the invariant
$\tau_\Rssf(\kfrak_\mfrak)$
does not necessarily
determine
the ray class
$\kfrak_\mfrak$;
rather,
one only has
the much weaker statement:
\textit{
$\tau_\Rssf(\kfrak_\mfrak)
 =
 \tau_\Rssf(\kfrak_\mfrak')
$
implies
$\kfrak_\mfrak
 =
 \kfrak_\mfrak'
$
if
$\kfrak_\mfrak$
and
$\kfrak_\mfrak'$
belong to the same
absolute class.
}

For then
one has
$\tau_\Rssf(\kfrak_\mfrak)
 =
 \tau_\Rssf(1,\mfrak_\Rssf\,\rfrak_\Rssf^{-1})
$
and
$\tau_\Rssf(\kfrak_\mfrak')
 =
 \tau_\Rssf(1,\mfrak_\Rssf\,{\rfrak_\Rssf'}^{-1})
$
with integral
$\rfrak,\rfrak'$
prime to $\mfrak$
and
$\rfrak'
 =
 \lambda\,\rfrak
$
with
$\lambda \in \Sigma$;
thus
\[
  \tau_\Rssf(\kfrak_\mfrak')
  \ =\ 
  \tau_\Rssf(1,\mfrak_\Rssf\,\rfrak_\Rssf^{-1}\,\lambda^{-1})
  \ =\ 
  \tau(\lambda,\mfrak_\Rssf\,\rfrak_\Rssf^{-1})
  ,
\]
consequently
one has
$\lambda
 \equiv
 \varepsilon
 \mod
 \mfrak\,\rfrak^{-1}
$
with a unit
$\varepsilon$ of $\Rssf$,
and since
$\rfrak$ is prime to $\mfrak$,
it follows that
$\lambda
 \equiv
 \varepsilon
 \ \textup{mod}^*\ 
 \mfrak
$,
thus
$\rfrak'$ and $\rfrak$
are equivalent
modulo~$\mfrak$,
which is to say
$\kfrak_\mfrak
 =
 \kfrak_\mfrak'
$.


\paragraph{26. The second main theorem.}

\textit{
The field
$\Omega_\Rssf(\tau_\Rssf(\kfrak_\mfrak))$
is the ray class field
of $\Sigma$
modulo~$\mfrak$.
}

Again
we first give a proof
by means of
the general theory of
abelian number fields
and then
a proof which,
except for
the elementary parts of
number theory,
relies only on
the theory of
modular functions
and
elliptic functions.

The first proof,
quite similar to that
in \S\textbf{III.2},
is based on
a congruence for
the ray class invariants,
namely:
\footnote{
\textit{Hasse},
loc. cit.${}^{2)}$
p.~134.
}

\textit{
Let $\pfrak$ be
a degree-one prime divisor
of $\Sigma$
not dividing
$\Nm\,\mfrak$,
which is distinct from
its conjugate
$\pfrakbar$,
and let $p$ be
the prime number
divisible by $\pfrak$.
Then
one has
\[
  \tau_\Rssf(\kfrak_\mfrak\,\pfrak^{-1})
  \ \equiv\ 
  \tau_\Rssf(\kfrak_\mfrak)^p
  \mod
  \Pfrak
  \tag{51}
  \label{eqn:51}
\]
for every prime factor
$\Pfrak$ of $\pfrak$
in
$\Omega_\Rssf(
 \tau_\Rssf(\kfrak_\mfrak)
 ,
 \tau_\Rssf(\kfrak_\mfrak\pfrak^{-1})
 )
$.
}

Proof.
We write
$\tau_\Rssf(\kfrak_\mfrak)$
in the form
\[
  \tau_\Rssf(\kfrak_\mfrak)
  \ =\ 
  \tau_\Rssf(1,\mfrak_\Rssf\,\rfrak_\Rssf^{-1})
  \ =\ 
  \tau_\Rssf
  \Bigl(\ 
  N^{-1}
  (x_1,x_2)
  \begin{pmatrix}
    \alpha_1
    \\
    \alpha_2
  \end{pmatrix}
  \,,\,
  \begin{pmatrix}
    \alpha_1
    \\
    \alpha_2
  \end{pmatrix}
  \ \Bigr)
  ,
\]
in which
$\alpha_1,\alpha_2$
is a basis of
$\afrak_\Rssf
 =
 \mfrak_\Rssf\,\rfrak_\Rssf^{-1}
$
and
$1
 =
 N^{-1}(x_1\alpha_1 + x_2\alpha_2)
$
with
relatively prime
rational integers
$x_1,x_2,N$;
here
$N$ contains
the same prime factors as
$\Nm\,\mfrak$
and is hence
not divisible by~$p$.
The matrix
$P_\pfrakbar$
of determinant~$p$
transforms
$\alpha_1,\alpha_2$
to a basis
$P_\pfrakbar
 \begin{pmatrix}
 \alpha_1
 \\
 \alpha_2
 \end{pmatrix}
$
of $\afrak_\Rssf\,\pfrak_\Rssf$.
Then
one has
\[
  \tau_\Rssf(\kfrak_\mfrak\,\pfrak^{-1})
  \ =\ 
  \tau_\Rssf(1,\mfrak_\Rssf\,\rfrak_\Rssf\,\pfrak_\Rssf^{-1})
  \ =\ 
  \tau_\Rssf(p,\mfrak_\Rssf\,\rfrak_\Rssf^{-1}\,\pfrakbar_\Rssf)
  \ =\ 
  \tau_\Rssf
  \Bigl(\ 
  p
  \,
  N^{-1}
  (x_1,x_2)
  \,,\,
  P_\pfrakbar
  \begin{pmatrix}
    \alpha_1
    \\
    \alpha_2
  \end{pmatrix}
  \ \Bigr)
\]
and thus
\[
  \tau_\Rssf(\kfrak_\mfrak)^p
  -
  \tau_\Rssf(\kfrak_\mfrak\,\pfrak^{-1})
  \ =\ 
  \delta_{P_\pfrak}
  (\,
  (x_1,x_2)
  \,{;}\,
  \alpha
  \,)
  ,
  \qquad
  \left(
  \alpha
  =
  \dfrac{\alpha_1}{\alpha_2}
  \right)
  .
\]
Here
$\delta_{P_\pfrak}
 (\,
 (x_1,x_2)
 \,{;}\,
 \alpha
 \,)
$
is a zero of
the polynomial
\[
  S_{P_\pfrakbar}(X,\alpha)
  \ =\ 
  \sum_\nu
  D_{P_\pfrakbar}^{(\nu)}(\alpha)
  \,
  X^\nu
  ,
\]
whose coefficients
$D_{P_\pfrakbar}^{(\nu)}(\alpha)$
are divisible by $\pfrak$
by \S\textbf{I.7} and \S\textbf{III.5}.
From this
the assertion
follows.

The proof of
the second main theorem
now goes
quite like
the proof of
the first main theorem
in~\S\textbf{III.2}.
\footnote{
\textit{Hasse},
loc. cit.${}^{2)}$
p.~137 to~138.
}

If $\pfrak$ is
a degree-one prime divisor
of $\Sigma$
in the ray modulo~$\mfrak$
which does not divide
the discriminant of $\Sigma$,
then by~\eqref{eqn:51},
for every ray class
$\kfrak_\mfrak$
and
the absolute divisor class
$\kfrak_\Rssf$
belonging to
$\kfrak_\mfrak$,
one has
\[
  \tau_\Rssf(\kfrak_\mfrak)^p
  \ \equiv\ 
  \tau_\Rssf(\kfrak_\mfrak)
  ,
  \qquad
  j(\kfrak_\Rssf)^p
  \ \equiv\ 
  j(\kfrak_\Rssf)
  \qquad
  \mod
  \Pfrak
\]
for every prime factor
$\Pfrak$ of $\pfrak$
in $\Omega_\Rssf(\tau_\Rssf(\kfrak_\mfrak))$.
If now
it is also assumed that
$\pfrak$
is prime to
the discriminant of
$j(\kfrak_\Rssf)$
and
$\tau_\Rssf(\kfrak_\mfrak)$
over $\Sigma$,
it follows that
for all
$\pfrak$-integral $\alpha$
from
$\Omega_\Rssf(\tau_\Rssf(\kfrak_\mfrak))$,
one has
the congruence
\[
  \alpha^p
  \ \equiv\ 
  \alpha
  \mod
  \Pfrak
  .
\]
Thus
$\Pfrak$
is of degree-one,
and
$\pfrak$
splits completely in
$\Omega_\Rssf(\tau_\Rssf(\kfrak_\mfrak))$.
Conversely
we suppose that
$\pfrak$
splits into
prime factors
$\Pfrak$
of degree-one
in $\Omega_\Rssf(\tau_\Rssf(\kfrak_\mfrak))$.
Then
one has
\[
  j(\kfrak_\Rssf)^p
  \ \equiv\ 
  j(\kfrak_\Rssf)
  \mod
  \Pfrak
  ;
\]
since
\[
  j(\kfrak_\Rssf)^p
  \ \equiv\ 
  j(\kfrak_\Rssf\,\pfrak^{-1})
  \mod
  \Pfrak
  ,
\]
thus
\[
  j(\kfrak_\Rssf\,\pfrak^{-1})
  \ \equiv\ 
  j(\kfrak_\Rssf)
  \mod
  \Pfrak
  ,
\]
and if we assume
$\pfrak$
is prime to
the differences
$j(\kfrak_\Rssf)
 -
 j(\kfrak_\Rssf')
$,
$\kfrak_\Rssf
 \neq
 \kfrak_\Rssf'
$,
then
\[
  j(\kfrak_\Rssf\,\pfrak_\Rssf)
  \ =\ 
  j(\kfrak_\Rssf)
  ,
  \qquad
  \text{so}
  \quad
  \kfrak_\Rssf\,\pfrak_\Rssf
  \ =\ 
  \kfrak_\Rssf
  ,
\]
and
$\pfrak$
is a principal divisor.
Furthermore
one has
\[
  \tau_\Rssf(\kfrak_\mfrak)^p
  \ \equiv\ 
  \tau_\Rssf(\kfrak_\mfrak)
  \mod
  \Pfrak
  ;
\]
since
\[
  \tau_\Rssf(\kfrak_\mfrak)^p
  \ \equiv\ 
  \tau_\Rssf(\kfrak_\mfrak\,\pfrak^{-1})
  \mod
  \Pfrak
  ,
\]
thus
\[
  \tau_\Rssf(\kfrak_\mfrak\,\pfrak^{-1})
  \ \equiv\ 
  \tau_\Rssf(\kfrak_\mfrak)
  \mod
  \Pfrak
  .
\]
As
$\kfrak_\mfrak$
and
$\kfrak_\mfrak\,\pfrak^{-1}$
belong to
the same absolute class
$\kfrak_\Rssf$,
it follows that
if we still assume
$\pfrak$
to be prime to
the (non-0) difference
\[
  \tau_\Rssf(\kfrak_\mfrak)
  -
  \tau_\Rssf(\kfrak_\mfrak')
  ,
  \qquad
  \kfrak_\mfrak
  \neq
  \kfrak_\mfrak'
  ,
  \qquad
  \text{
  $\kfrak_\mfrak$
  and
  $\kfrak_\mfrak'$
  in the same
  absolute class
  }
  ,
\]
then
\[
  \tau_\Rssf(\kfrak_\mfrak\,\pfrak^{-1})
  \ =\ 
  \tau_\Rssf(\kfrak_\mfrak)
  ,
  \qquad
  \text{so}
  \quad
  \kfrak_\mfrak\,\pfrak
  \ =\ 
  \kfrak_\mfrak
  ,
\]
and
$\pfrak$
lies in the ray
modulo~$\mfrak$.


\paragraph{27. Proof of the second main theorem
not relying on
general class field theory.}

The proof of the second main theorem
not relying on
general class field theory
also runs parallel to
the corresponding proof of
the first main theorem.
We set
\[
  \Omega_\mfrak
  \ =\ 
  \Omega_\Rssf
  (\,
  \tau_\Rssf(\kfrak_\mfrak^{(1)})
  \,,\,
  \ldots
  \,,\,
  \tau_\Rssf(\kfrak_\mfrak^{(h_\mfrak)})
  \,)
\]
and prove that
$\Omega_\mfrak$
is the ray class field
of $\Sigma$
modulo~$\mfrak$,
by showing that
the Artin reciprocity law
holds for
$\Omega_\mfrak$
(as class field of $\Sigma$
for the ray
modulo~$\mfrak$),
thereby also getting
\[
 \Omega_\mfrak
 =
 \Omega_\Rssf(\tau_\Rssf(\kfrak_\mfrak))
 =
 \Sigma(j(\kfrak_\Rssf),\tau_\Rssf(\kfrak_\mfrak))
 .
\]

The polynomial
with roots
$\tau_\Rssf(\kfrak_\mfrak)$,
thus
\[
  S_\mfrak(X)
  \ =\ 
  \prod_{\kfrak_\mfrak}
  (X-\tau_\Rssf(\kfrak_\mfrak))
  ,
\]
is called
the \emph{ray class polynomial
of $\Sigma$
modulo~$\mfrak$}.
If $\kfrak_\Rssf$ is
a divisor class of $\Sigma$,
then
\[
  S_{\mfrak,\kfrak_\Rssf}(X)
  \ =\ 
  \prod_{
    \text{$\kfrak_\mfrak$ in $\kfrak_\Rssf$}
  }
  (X-\tau_\Rssf(\kfrak_\mfrak))
\]
shall be called
the \emph{ray class polynomial
of $\Sigma$
modulo~$\mfrak$
over the class
$\kfrak_\Rssf$};
one has
\[
  S_\mfrak(X)
  \ =\ 
  \prod_{\kfrak_\Rssf}
  S_{\mfrak,\kfrak_\Rssf}(X)
  .
\]
\textit{
The coefficients of
$S_{\mfrak,\kfrak_\Rssf}(X)$
lie in $\Omega_\Rssf$.
}
This fact
plays the same role
in the proof of
the second main theorem
as the fact that
the class polynomial
$H_\Rssf(X)$
has rational coefficients
plays in the proof of
the first main theorem;
in order to
not interrupt
the train of thoughts,
we postpone
the proof of this fact
to the end of
our deliberations.
%


\paragraph{28. The fundamental congruence.}

While we got by
for the proof of
the first main theorem
with the congruence
\[
  j(\kfrak_\Rssf\,\pfrak^{-1})
  \ \equiv\ 
  j(\kfrak_\Rssf)^p
  \mod
  \pfrak
\]
which is
valid for
prime divisors $\pfrak$
of degree~\emph{one},
the corresponding congruence
\[
  \tau_\Rssf(\kfrak_\mfrak\,\pfrak^{-1})
  \ \equiv\ 
  \tau_\Rssf(\kfrak_\mfrak)^p
  \mod
  \pfrak
\]
is sufficient
for the proof of
the second main theorem
only if
we want to make use of
the fact that
in each ray class
modulo~$\mfrak$
there exist divisors
which are composed of
only prime factors
of degree~one,
which is a fact that
can only be seen
by means of
analytic methods
or by the use of
essential tools from
general class field theory.
This prompts us
to generalize
the congruence~\eqref{eqn:51}
as follows:

\textit{
Let $\pfrak$ be
a prime divisor
in $\Sigma$,
in which
the prime number~$p$
divisible by $\pfrak$
is larger than~12.
Suppose
$\pfrak$
is not a factor of
$\Nm\,\mfrak$.
Then one has
\[
  \tau_\Rssf(\kfrak_\mfrak\,\pfrak^{-1})
  \ \equiv\ 
  \tau_\Rssf(\kfrak_\mfrak)^{\Nm\pfrak}
  \mod
  \Pfrak
  \tag{52}
  \label{eqn:52}
\]
for every prime factor
$\Pfrak$ of $\pfrak$
in $\Omega_\mfrak$.
}

Proof.
If
$p = \pfrak\,\pfrakbar$,
$\pfrak \neq \pfrakbar$,
that is
the assertion
proven in~\S\textbf{3}
(even without
the assumption
$p{\,>\,}12$).
In the case of
$p = \pfrak^2$,
the proof
in~\S\textbf{3}
can be taken
word-for-word,
except
instead of
the theorem of~\S\textbf{13}
one has to use
the theorem of~\S\textbf{23}.
There remains
the case of
$p = \pfrak$.
We again set
$\tau_\Rssf(\kfrak_\mfrak)
 =
 \tau_\Rssf(1,\mfrak_\Rssf\,\rfrak_\Rssf^{-1})
$
by means of
a basis
$\alpha_1,\alpha_2$
of $\mfrak_\Rssf\,\rfrak_\Rssf^{-1}$
and the representation
$1 = N^{-1}(x_1\alpha_1 + x_2\alpha_2)$,
$(x_1,x_2,N) = 1$,
$p \nmid N$
in the form
\[
  \tau_\Rssf(\kfrak_\mfrak)
  \ =\ 
  \tau_\Rssf
  \Bigl(\ 
  N^{-1}
  (x_1,x_2)
  \begin{pmatrix}
    \alpha_1
    \\
    \alpha_2
  \end{pmatrix}
  \,,\,
  \begin{pmatrix}
    \alpha_1
    \\
    \alpha_2
  \end{pmatrix}
  \ \Bigr)
  .
\]
Let $P$ denote
an arbitrary
(but chosen fixed)
matrix
of determinant~$p$.
Then
\[
  \begin{gathered}
  \tau_\Rssf(\kfrak_\mfrak)^p
  \ -\ 
  \tau_\Rssf
  \Bigl(\ 
  p
  \,
  N^{-1}
  (x_1,x_2)
  \begin{pmatrix}
    \alpha_1
    \\
    \alpha_2
  \end{pmatrix}
  \,,\,
  P
  \begin{pmatrix}
    \alpha_1
    \\
    \alpha_2
  \end{pmatrix}
  \ \Bigr)
  \\
  \ =\ 
  \tau_\Rssf
  \Bigl(\ 
  N^{-1}
  (x_1,x_2)
  \begin{pmatrix}
    \alpha_1
    \\
    \alpha_2
  \end{pmatrix}
  \,,\,
  \begin{pmatrix}
    \alpha_1
    \\
    \alpha_2
  \end{pmatrix}
  \ \Bigr)^p
  \ -\ 
  \tau_\Rssf
  \Bigl(\ 
  p
  \,
  N^{-1}
  (x_1,x_2)
  \begin{pmatrix}
    \alpha_1
    \\
    \alpha_2
  \end{pmatrix}
  \,,\,
  P
  \begin{pmatrix}
    \alpha_1
    \\
    \alpha_2
  \end{pmatrix}
  \ \Bigr)
  \\
  \ =\ 
  \delta_P
  \bigl(\,
  (x_1,x_2)
  \,{;}\,
  \alpha
  \,\bigr)
  ,
  \qquad
  \text{
  (
  $\alpha
   =
   \dfrac{\alpha_1}{\alpha_2}
  $
  )
  }
  \end{gathered}
\]
is a zero of
${\displaystyle
 S_P(X,\alpha)
 =
 \sum_{\nu}
 D_P^{(\nu)}(\alpha)
 X^\nu
}$,
whose coefficients
$D_P^{(\nu)}(\alpha)$
are divisible by~$\pfrak$
by~\S\textbf{I.7}
and~\S\textbf{III.15},
and hence
for every prime factor
$\Pfrak^*$ of $\pfrak$
in
${\displaystyle
 \Omega_\mfrak
 \left(
  \tau_\Rssf
  \Bigl(
  p
  \,
  N^{-1}
  (x_1,x_2)
  \begin{pmatrix}
    \alpha_1
    \\
    \alpha_2
  \end{pmatrix}
  \,,\,
  P
  \begin{pmatrix}
    \alpha_1
    \\
    \alpha_2
  \end{pmatrix}
  \Bigr)
 \right)
}$
one has
the congruence
\[
  \tau_\Rssf(\kfrak_\mfrak)^p
  \ \equiv\ 
  \tau_\Rssf
  \Bigl(\ 
  p
  \,
  N^{-1}
  (x_1,x_2)
  \begin{pmatrix}
    \alpha_1
    \\
    \alpha_2
  \end{pmatrix}
  \,,\,
  P
  \begin{pmatrix}
    \alpha_1
    \\
    \alpha_2
  \end{pmatrix}
  \ \Bigr)
  \mod
  \Pfrak^*
  .
  \tag{53}
  \label{eqn:53}
\]
Now
$P'
 =
 p
 \,
 P^{-1}
$
is a matrix
of determinant~$p$.
We set
${\displaystyle
 P
 \begin{pmatrix}
   \alpha_1
   \\
   \alpha_2
 \end{pmatrix}
 =
 \begin{pmatrix}
   \alpha_1'
   \\
   \alpha_2'
 \end{pmatrix}
}$,
so
${\displaystyle
 P'
 \begin{pmatrix}
   \alpha_1'
   \\
   \alpha_2'
 \end{pmatrix}
 =
 \begin{pmatrix}
   \alpha_1
   \\
   \alpha_2
 \end{pmatrix}
 \,
 p
}$;
furthermore
$(x_1,x_2)
 \,
 P'
 =
 (x_1',x_2')
$,
so
$p
 \,
 (x_1,x_2)
 =
 (x_1',x_2')
 \,
 P
$,
and
one obviously has
$(x_1',x_2',N)
 =
 1
$.
We compute
${\displaystyle
  \tau_\Rssf
  \Bigl(\ 
  p
  \,
  N^{-1}
  (x_1,x_2)
  \begin{pmatrix}
    \alpha_1
    \\
    \alpha_2
  \end{pmatrix}
  \,,\,
  P
  \begin{pmatrix}
    \alpha_1
    \\
    \alpha_2
  \end{pmatrix}
  \ \Bigr)
}$
in terms of
$(\alpha_1',\alpha_2')$
and
$(x_1',x_2')$
as:
\[
  \tau_\Rssf
  \Bigl(\ 
  p
  \,
  N^{-1}
  (x_1,x_2)
  \begin{pmatrix}
    \alpha_1
    \\
    \alpha_2
  \end{pmatrix}
  \,,\,
  P
  \begin{pmatrix}
    \alpha_1
    \\
    \alpha_2
  \end{pmatrix}
  \ \Bigr)
  \ =\ 
  \tau_\Rssf
  \Bigl(\ 
  N^{-1}
  (x_1',x_2')
  \begin{pmatrix}
    \alpha_1'
    \\
    \alpha_2'
  \end{pmatrix}
  \,,\,
  \begin{pmatrix}
    \alpha_1'
    \\
    \alpha_2'
  \end{pmatrix}
  \ \Bigr)
  .
\]
On the other hand,
one has
\[
  \begin{gathered}
  \tau_\Rssf
  (\kfrak_\mfrak\,\pfrak^{-1})
  \ =\ 
  \tau_\Rssf
  (\kfrak_\mfrak\,p^{-1})
  \ =\ 
  \tau_\Rssf
  (\,
  1
  \,,\,
  \mfrak_\Rssf\,\rfrak_\Rssf^{-1}\,p^{-1}
  \,)
  \ =\ 
  \tau_\Rssf
  (\,
  p^2
  \,,\,
  \mfrak_\Rssf\,\rfrak_\Rssf^{-1}\,p
  \,)
  \\
  \ =\ 
  \tau_\Rssf
  \Bigl(\ 
  p^2
  \,
  N^{-1}
  (x_1,x_2)
  \begin{pmatrix}
    \alpha_1
    \\
    \alpha_2
  \end{pmatrix}
  \,,\,
  \begin{pmatrix}
    \alpha_1
    \\
    \alpha_2
  \end{pmatrix}
  \,
  p
  \ \Bigr)
  \ =\ 
  \tau_\Rssf
  \Bigl(\ 
  p
  \,
  N^{-1}
  (x_1',x_2')
  \begin{pmatrix}
    \alpha_1'
    \\
    \alpha_2'
  \end{pmatrix}
  \,,\,
  P'
  \begin{pmatrix}
    \alpha_1'
    \\
    \alpha_2'
  \end{pmatrix}
  \ \Bigr)
  ,
  \end{gathered}
\]
thus
\[
  \tau_\Rssf
  \Bigl(\ 
  p
  \,
  N^{-1}
  (x_1,x_2)
  \begin{pmatrix}
    \alpha_1
    \\
    \alpha_2
  \end{pmatrix}
  \,,\,
  P
  \begin{pmatrix}
    \alpha_1
    \\
    \alpha_2
  \end{pmatrix}
  \ \Bigr)^p
  \ -\ 
  \tau_\Rssf
  (\kfrak_\mfrak\,\pfrak^{-1})
  \ =\ 
  \delta_{P'}
  \bigl(\,
  (x_1',x_2')
  \,;\,
  \alpha'
  \,\bigr)
  ,
\]
with
$\alpha'
 =
 \dfrac{\alpha_1'}{\alpha_2'}
$.
As above
this now yields
\[
  \tau_\Rssf
  \Bigl(\ 
  p
  \,
  N^{-1}
  (x_1,x_2)
  \begin{pmatrix}
    \alpha_1
    \\
    \alpha_2
  \end{pmatrix}
  \,,\,
  P
  \begin{pmatrix}
    \alpha_1
    \\
    \alpha_2
  \end{pmatrix}
  \ \Bigr)^p
  \ \equiv\ 
  \tau_\Rssf
  (\kfrak_\mfrak\,\pfrak^{-1})
  \mod
  \Pfrak^*
  ;
\]
together with~\eqref{eqn:53},
one has
\[
  \tau_\Rssf(\kfrak_\mfrak)^{p^2}
  \ \equiv\ 
  \tau_\Rssf(\kfrak_\mfrak\,\pfrak^{-1})
  \mod
  \Pfrak
  ,
\]
and
since
$\Nm\,\pfrak = p^2$,
that is
the assertion.

It should be noted that
the congruence
\[
  j(\kfrak_\Rssf\,\pfrak^{-1})
  \ \equiv\ 
  j(\kfrak_\Rssf)^{\Nm\pfrak}
  \mod
  \Pfrak
\]
for all prime factors
$\Pfrak$ of $\pfrak$
in $\Omega_\Rssf$
can also be proven
in this way
for degree-two prime divisors
(under the assumption
$p{\,>\,}12$);
but of course
it follows directly from
the Artin reciprocity law
for $\Omega_\Rssf/\Sigma$
proven in~\S\textbf{III.9}.


\paragraph{29. The generation of the ray class field
by a single ray class invariant.}

We now prove,
with the reasoning
already used in~\S\textbf{III.7},
that
one has
$\Omega_\mfrak
 =
 \Omega_\Rssf
 (\tau_\Rssf(\kfrak_\mfrak))
$.
In any case
$\Omega_\mfrak$
is normal over
$\Omega_\Rssf$,
because
the coefficients of
the polynomial
$S_{\mfrak,\kfrak_\Rssf}(X)$
lie in $\Omega_\Rssf$.
Let $\lambda$ be
an automorphism of
$\Omega_\mfrak
 /
 \Omega_\Rssf
 (\tau_\Rssf(\kfrak_\mfrak))
$
(where
$\kfrak_\mfrak$
is a \emph{fixed}
ray class).
We have to show that
for \emph{all}
ray classes
$\hfrak_\mfrak$,
one has
\[
  \tau_\Rssf
  (\hfrak_\mfrak)
  \,
  \lambda
  \ =\ 
  \tau_\Rssf
  (\hfrak_\mfrak)
  .
\]
This is true
for a certain
ray class
$\hfrak_\mfrak$.
If $\pfrak$ is
a prime divisor of
$\Sigma$
not dividing
$\Nm\,\mfrak$,
unramified in
$\Omega_\mfrak$,
one has
\[
  \tau_\Rssf
  (\hfrak_\mfrak\,\pfrak^{-1})
  \ \equiv\ 
  \tau_\Rssf
  (\hfrak_\mfrak)^{\Nm\pfrak}
  \mod
  \pfrak
  ;
\]
and since
$\pfrak\,\lambda
 =
 \pfrak
$,
thus
\[
  \begin{gathered}
  \tau_\Rssf
  (\hfrak_\mfrak\,\pfrak^{-1})
  \,
  \lambda
  \ \equiv\ 
  (\,
  \tau_\Rssf
  (\hfrak_\mfrak)
  \,
  \lambda
  \,)^{\Nm\pfrak}
  \ =\ 
  \tau_\Rssf
  (\hfrak_\mfrak)^{\Nm\pfrak}
  \mod
  \pfrak
  ,
  \\
  \text{and so}
  \qquad
  \tau_\Rssf
  (\hfrak_\mfrak\,\pfrak^{-1})
  \,
  \lambda
  \ \equiv\ 
  \tau_\Rssf
  (\hfrak_\mfrak\,\pfrak^{-1})
  \mod
  \pfrak
  .
  \end{gathered}
\]
But since
$\Omega_\Rssf$
remains
elementwise fixed
under $\lambda$,
$\tau_\Rssf
 (\hfrak_\mfrak\,\pfrak^{-1})
 \,
 \lambda
$
is also
a zero of
$S_{\mfrak,\hfrak_\Rssf\pfrak^{-1}}(X)$
just as
$\tau_\Rssf
 (\hfrak_\mfrak\,\pfrak^{-1})
$
is,
where
$\hfrak_\Rssf$
denotes
the absolute class
to which
$\hfrak_\mfrak$
belongs.
As it is assumed that
$\pfrak$
is prime to
the discriminant of
$S_{\mfrak,\hfrak_\Rssf\pfrak}(X)$,
it follows that
\[
  \tau_\Rssf
  (\hfrak_\mfrak\,\pfrak^{-1})
  \,
  \lambda
  \ =\ 
  \tau_\Rssf
  (\hfrak_\mfrak\,\pfrak^{-1})
  .
\]
From
$\tau_\Rssf
 (\kfrak_\mfrak)
 \,
 \lambda
 =
 \tau_\Rssf
 (\kfrak_\mfrak)
$
we can now
go step-by-step to
$\tau_\Rssf
 (\hfrak_\mfrak)
 \,
 \lambda
 =
 \tau_\Rssf
 (\hfrak_\mfrak)
$
for an arbitrary
$\hfrak_\mfrak$.


\paragraph{30. The reciprocity law
for the ray class field.}

We still do not know
whether
$\Omega_\mfrak$
is normal over
$\Sigma$.
Let $\Omega^*$ be
an extension field of
$\Omega_\mfrak$
normal over
$\Sigma$.
Let $\pfrak$ be
a prime divisor of
$\Sigma$
not dividing
$\Nm\,\mfrak$,
let $\Pfrak^*$ be
a prime factor of $\pfrak$
in $\Omega^*$
and
let $F_{\Pfrak^*}$ be
the Frobenius automorphism of
$\Pfrak^*/\Sigma$.
Then
one has
\[
  j(\kfrak_\Rssf)^{\Nm\pfrak}
  \ \equiv\ 
  j(\kfrak_\Rssf)
  \,
  F_{\Pfrak^*}
  ,
  \qquad
  \tau_\Rssf(\kfrak_\mfrak)^{\Nm\pfrak}
  \ \equiv\ 
  \tau_\Rssf(\kfrak_\mfrak)
  \,
  F_{\Pfrak^*}
  \mod
  \Pfrak^*
  ,
\]
and since
\[
  j(\kfrak_\Rssf)^{\Nm\pfrak}
  \ \equiv\ 
  j(\kfrak_\Rssf\,\pfrak_\Rssf^{-1})
  ,
  \qquad
  \tau_\Rssf(\kfrak_\mfrak)^{\Nm\pfrak}
  \ \equiv\ 
  \tau_\Rssf(\kfrak_\mfrak\,\pfrak^{-1})
  \mod
  \Pfrak^*
  ,
\]
it follows that
\[
  j(\kfrak_\Rssf)
  \,
  F_{\Pfrak^*}
  \ \equiv\ 
  j(\kfrak_\Rssf\,\pfrak_\Rssf^{-1})
  ,
  \qquad
  \tau_\Rssf(\hfrak_\mfrak)
  \,
  F_{\Pfrak^*}
  \ \equiv\ 
  \tau_\Rssf(\hfrak_\mfrak\,\pfrak^{-1})
  \mod
  \Pfrak^*
\]
for all
absolute classes
$\kfrak_\Rssf$
and all
ray classes
$\kfrak_\mfrak$.
We now assume that
$\pfrak$
is prime to
the discriminant of
$H_\Rssf(X)$
and
to all the differences of
the (algebraic) numbers
which are
conjugate over $\Sigma$
to the ray class invariants
$\tau_\Rssf(\kfrak_\mfrak)$,
as long as
they are
distinct from~0;
that excludes
only finitely many~$\pfrak$.
Then
it follows that
\[
  j(\kfrak_\Rssf)
  \,
  F_{\Pfrak^*}
  \ =\ 
  j(\kfrak_\Rssf\,\pfrak_\Rssf^{-1})
  ,
  \qquad
  \tau_\Rssf(\kfrak_\mfrak)
  \,
  F_{\Pfrak^*}
  \ =\ 
  \tau_\Rssf(\kfrak_\mfrak\,\pfrak^{-1})
  .
\]
Let $\hfrak_\mfrak$ be
a ray class,
let $\hfrak_\Rssf$ be
the absolute class
to which
$\hfrak_\mfrak$
belongs,
let
$\pfrak_1\pfrak_2\cdots\pfrak_\mfrak$
be an integral divisor
in $\hfrak_\mfrak$
with
none of its prime factors
$\pfrak_\nu$
occuring among
the excluded
prime divisors,
and let
$\Pfrak_\nu^*$
be a prime factor of
$\pfrak_\nu$
in $\Omega^*$.
The automorphism
\[
  \sigma_\mfrak(\hfrak_\mfrak)
  \ =\ 
  F_{\Pfrak_1^*}
  \,
  F_{\Pfrak_2^*}
  \,
  \cdots
  \,
  F_{\Pfrak_\mfrak^*}
\]
of $\Omega^*/\Sigma$
satisfies the equations
\[
  j(\kfrak_\Rssf)
  \,
  \sigma_\mfrak(\hfrak_\mfrak)
  \ =\ 
  j(\kfrak_\Rssf\,\hfrak_\Rssf^{-1})
  ,
  \qquad
  \tau_\Rssf(\kfrak_\mfrak)
  \,
  \sigma_\mfrak(\hfrak_\mfrak)
  \ =\ 
  \tau_\Rssf(\kfrak_\mfrak\,\hfrak_\mfrak^{-1})
  .
  \tag{54}
  \label{eqn:54}
\]
The automorphism
$\sigma(\hfrak_\mfrak)$
is not uniquely determined
by these;
rather,
the most general
automorphism $\lambda$
of $\Omega^*/\Sigma$
with
\[
  j(\kfrak_\Rssf)
  \,
  \lambda
  \ =\ 
  j(\kfrak_\Rssf\,\hfrak_\Rssf^{-1})
  ,
  \qquad
  \tau_\Rssf(\kfrak_\mfrak)
  \,
  \lambda
  \ =\ 
  \tau_\Rssf(\kfrak_\mfrak\,\hfrak_\mfrak)
\]
is of the form
\[
  \lambda
  \ =\ 
  \varrho
  \,
  \sigma_\mfrak(\hfrak_\mfrak)
\]
with an arbitrary
$\varrho$
from the Galois group of
$\Omega^*/\Omega_\mfrak$.
One has
\[
  \varrho
  \,
  \sigma_\mfrak(\hfrak_\mfrak)
  \ =\ 
  \varrho'
  \,
  \sigma_\mfrak(\hfrak_\mfrak')
\]
if and only if
one has
$\varrho
 =
 \varrho'
$
and
$\hfrak_\mfrak
 =
 \hfrak_\mfrak'
$.
This is because
from
$\varrho
 \,
 \sigma_\mfrak(\hfrak_\mfrak)
 =
 \varrho'
 \,
 \sigma_\mfrak(\hfrak_\mfrak')
$
it follows that
$j(\kfrak_\Rssf\,\hfrak_\Rssf^{-1})
 =
 j(\kfrak_\Rssf\,{\hfrak_\Rssf'}^{-1})
$,
thus
$\hfrak_\Rssf
 =
 \hfrak_\Rssf'
$,
so
$\hfrak_\mfrak$
and
$\hfrak_\mfrak'$
lie in the same
absolute class,
whence
$\tau_\Rssf(\kfrak_\mfrak\,\hfrak_\mfrak^{-1})
 =
 \tau_\Rssf(\kfrak_\mfrak\,{\hfrak_\mfrak'}^{-1})
$,
and as
$\hfrak_\mfrak$
and
$\hfrak_\mfrak'$
lie in the same
absolute class,
one has
$\hfrak_\mfrak
 =
 \hfrak_\mfrak'
$,
and
$\varrho = \varrho'$.

The
$\varrho
 \,
 \sigma_\mfrak(\hfrak_\mfrak)
$
together
are thus
$h_\mfrak
 \,
 [\,
 \Omega^*
 \,{:}\,
 \Omega_\mfrak
 \,]
$
elements of
the Galois group of
$\Omega^*/\Sigma$.
But the degree
$[\,
 \Omega^*
 \,{:}\,
 \Sigma
 \,]
$
is at most equal to
$h_\mfrak
 \,
 [\,
 \Omega^*
 \,{:}\,
 \Omega_\mfrak
 \,]
$,
because
from
$\Omega_\mfrak
 =
 \Omega_\Rssf(\tau_\Rssf(\kfrak_\mfrak))
$,
it follows that
$[\,
 \Omega_\mfrak
 \,{:}\,
 \Omega_\Rssf
 \,]
$
is at most equal to
the degree
$h_\mfrak/h_\Rssf$
of $S_{\mfrak,\kfrak_\mfrak}(X)$.
Consequently
the
$\varrho
 \,
 \sigma_\mfrak(\hfrak_\mfrak)
$
are elements of
the \emph{whole}
Galois group of
$\Omega^*/\Sigma$.
From this
we can conclude that
$\Omega_\mfrak$
is normal over
$\Sigma$.
For a conjugate
$(\varrho
 \,
 \sigma(\hfrak_\mfrak)
 )^{-1}
 \,
 \varrho_0
 \,
 (\varrho
 \,
 \sigma(\hfrak_\mfrak)
 )
$
of an automorphism
$\varrho_0$
of $\Omega^*/\Omega_\mfrak$,
one has
namely
\[
  \begin{gathered}
  \begin{aligned}
  j(\kfrak_\Rssf)
  \,
  (\varrho
  \,
  \sigma(\hfrak_\mfrak)
  )^{-1}
  \,
  \varrho_0
  \,
  (\varrho
  \,
  \sigma(\hfrak_\mfrak)
  )
  &
  \ =\ 
  j(\kfrak_\Rssf)
  \,
  \sigma(\hfrak_\mfrak)^{-1}
  \,
  \varrho^{-1}
  \,
  \varrho_0
  \,
  \varrho
  \,
  \sigma(\hfrak_\mfrak)
  \\
  &
  \ =\ 
  j(\kfrak_\Rssf\hfrak_\Rssf)
  \,
  \varrho^{-1}
  \,
  \varrho_0
  \,
  \varrho
  \,
  \sigma(\hfrak_\mfrak)
  \ =\ 
  j(\kfrak_\Rssf\hfrak_\Rssf)
  \,
  \sigma(\hfrak_\mfrak)
  \ =\ 
  j(\kfrak_\Rssf)
  ,
  \end{aligned}
  \\[2ex]
  \begin{aligned}
  \tau_\Rssf(\kfrak_\Rssf)
  \,
  (\varrho
  \,
  \sigma(\hfrak_\mfrak)
  )^{-1}
  \,
  \varrho_0
  \,
  (\varrho
  \,
  \sigma(\hfrak_\mfrak)
  )
  &
  \ =\ 
  \tau_\Rssf(\kfrak_\mfrak)
  \,
  \sigma(\hfrak_\mfrak)
  \,
  \varrho^{-1}
  \,
  \varrho_0
  \,
  \varrho
  \,
  \sigma(\hfrak_\mfrak)
  \\
  &
  \ =\ 
  \tau_\Rssf(\kfrak_\mfrak\,\hfrak_\mfrak)
  \,
  \varrho^{-1}
  \,
  \varrho_0
  \,
  \varrho
  \,
  \sigma(\hfrak_\mfrak)
  \ =\ 
  \tau_\Rssf(\kfrak_\mfrak\,\hfrak_\mfrak)
  \,
  \sigma(\hfrak_\mfrak)
  \ =\ 
  \tau_\Rssf(\kfrak_\mfrak)
  .
  \end{aligned}
  \end{gathered}
\]
It is thus also
an automorphism of
$\Omega^*/\Omega_\mfrak$;
the Galois group of
$\Omega^*/\Omega_\mfrak$
is a normal subgroup of
the Galois group of
$\Omega^*/\Sigma$,
as claimed.

We can therefore
assume
$\Omega^* = \Omega_\mfrak$;
then
$\sigma_\mfrak(\hfrak_\mfrak)$
is uniquely determined by
the equations~\eqref{eqn:54},
from which
it follows that
\[
  \sigma_\mfrak(\hfrak_\mfrak^{(1)})
  \,
  \sigma_\mfrak(\hfrak_\mfrak^{(2)})
  \ =\ 
  \sigma_\mfrak(\hfrak_\mfrak^{(1)}\,\hfrak_\mfrak^{(2)})
  ,
\]
so
$\sigma_\mfrak$
is an isomorphism
from the ray class group of
$\Sigma$ modulo~$\mfrak$
onto
the Galois group of
$\Omega_\mfrak/\Sigma$,
and the definition of
$\sigma_\mfrak(\hfrak_\mfrak)$
by means of
the Frobenius automorphisms
shows that
$\sigma_\mfrak$
is exactly
the isomorphism asserted by
the Artin reciprocity law.
For that, however,
finitely many
prime divisors of $\Sigma$
are excluded.


\paragraph{31. The ray class polynomial.}

It still remains
to carry out
the proof that
$S_{\mfrak,\kfrak_\Rssf}(X)$
has coefficients in
$\Omega_\Rssf$.
\footnote{
\textit{Hasse},
loc. cit.${}^{8)}$,
78 to~81.
}

We introduce
the polynomial
\[
  T_{\mfrak,\kfrak_\Rssf}(X)
  \ =\ 
  \prod_
  {\substack{
      \gamma \bmod \afrak
      ,
      \\
      \gamma \equiv 0
      \bmod
      \afrak\,\mfrak^{-1}
      ,
      \\
      \gamma \not\equiv 0
      \bmod
      \afrak
  }}
  (X - \tau_\Rssf(\gamma,\afrak_\Rssf))
  ,
\]
whose zeros
are thus
the $\mfrak$-th division values of
the ideal
$\afrak_\Rssf$ of $\Rssf$;
in the notation
it is expressed that
$T_{\mfrak,\kfrak_\Rssf}(X)$
depends only on
the class $\kfrak_\Rssf$
of the divisor $\afrak$.

The
$\tau_\Rssf(\gamma,\afrak_\Rssf)$
are exactly
\emph{the} division values of
$\afrak_\Rssf$
whose orders
$\tfrak$
divide $\mfrak$.
The
$\tau_\Rssf(\gamma,\afrak_\Rssf)$
whose order
is a \emph{given} divisor
$\tfrak$
of $\mfrak$
are exactly
the invariants
$\tau_\Rssf(\kfrak_\tfrak)$
of those ray classes
of $\Sigma$
modulo~$\tfrak$
which lie in
the absolute class
$\kfrak_\Rssf$,
and since
two invariants
$\tau_\Rssf(\gamma',\afrak_\Rssf),
 \tau_\Rssf(\gamma,\afrak_\Rssf)
$
of order~$\tfrak$
are equal
if and only if
the prinicipal divisor
$(\gamma'\,\gamma^{-1})$
lies in the ray
modulo~$\tfrak$,
it follows that
every invariant
$\tau(\kfrak_\tfrak)$
occurs as
a zero of
$T_{\mfrak,\kfrak_\Rssf}(X)$
exactly
$e_\tfrak$~times,
where
$e_\tfrak$
denotes
the number of
residue classes of $\Rssf$
modulo~$\tfrak$
represented by units.
If we set
$S_{1,\kfrak_\Rssf}(X) = 1$,
one then has
\[
  T_{\mfrak,\kfrak_\Rssf}(X)
  \ =\ 
  \prod_{\tfrak \mid \mfrak}
  S_{\tfrak,\kfrak_\Rssf}(X)^{e_\tfrak}
  ,
\]
and this holds
also for
$\mfrak = 1$
if we set
$T_{1,\kfrak_\Rssf}(X) = 1$.
The M\"obius inversion formula
yields
\[
  S_{\mfrak,\kfrak_\Rssf}(X)^{e_\tfrak}
  \ =\ 
  \prod_{\tfrak \mid \mfrak}
  T_{\tfrak,\kfrak_\Rssf}(X)^{\scriptstyle \mu\left(\frac{\mfrak}{\tfrak}\right)}
  .
\]
It thus suffices
to prove that
the coefficients of
$T_{\mfrak,\kfrak_\Rssf}(X)$
lie in $\Omega_\Rssf$.

The singular elliptic function
$\tau_\Rssf(z, \afrak_\Rssf\,\mfrak_\Rssf^{-1})$
belongs to
the field
$K_0(\afrak_\Rssf)$,
because
the period lattice
$\afrak_\Rssf$
of $K(\afrak_\Rssf)$
is contained in
$\afrak_\Rssf\,\mfrak_\Rssf^{-1}$
and
$\tau_\Rssf(z, \afrak_\Rssf\,\mfrak_\Rssf^{-1})$
is invariant
under the
$\varrho_\varepsilon$.
Hence
$\tau_\Rssf(z, \afrak_\Rssf\,\mfrak_\Rssf^{-1})$
is a rational function of
$\tau_\Rssf(z, \afrak_\Rssf)$:
\[
  \tau_\Rssf(z, \afrak_\Rssf\,\mfrak_\Rssf^{-1})
  \ =\ 
  \dfrac
  {Z(\tau_\Rssf(z, \afrak_\Rssf))}
  {N(\tau_\Rssf(z, \afrak_\Rssf))}
  ,
  \tag{55}
  \label{eqn:55}
\]
where
$Z(X),N(X)$
are relatively prime
polynomials,
$N$ with highest coefficient~1.
We will show that
$N(X)
 =
 T_{\mfrak,\kfrak_\Rssf}(X)
$;
for that,
we compare
the poles of
both sides of~\eqref{eqn:55}.
The poles of
$\tau_\Rssf(z, \afrak_\Rssf\,\mfrak_\Rssf^{-1})$
are at
\[
  z
  \ \equiv\ 
  0
  \mod
  \afrak_\Rssf\,\mfrak_\Rssf^{-1}
  ;
\]
they are of order~$e$.
The right hand side of~\eqref{eqn:55}
has order
$e\,(\deg N {-} \deg Z)$
at
$z \equiv 0
 \mod
 \afrak_\Rssf
$,
whence
one has
$\deg Z
 =
 \deg N
 +
 1
$
and
$N(\tau_\Rssf(z,\afrak_\Rssf))$
must have
$e$-th order zeros
at
$z \equiv 0
 \mod
 \afrak_\Rssf\,\mfrak_\Rssf
$,
$z \not\equiv 0
 \mod
 \afrak_\Rssf
$,
but is otherwise~$\neq 0$.
The zeros of
$N(X)$
are thus
the numbers
$\tau_\Rssf(\gamma,\afrak_\Rssf)$,
$\gamma \equiv 0
 \mod
 \afrak_\Rssf\,\mfrak_\Rssf^{-1}
$,
$\gamma \not\equiv 0
 \mod
 \afrak_\Rssf
$,
thus
exactly the zeros of
$T_{\mfrak,\kfrak_\Rssf}(X)$,
and
it remains to show that
$\tau_\Rssf(\gamma,\afrak_\Rssf)$
as a zero of
$N(X)$
has the same multiplicity
as it does
as a zero of
$T_{\mfrak,\kfrak_\Rssf}(X)$.
If $n_\gamma$ is
the multiplicity of
$\tau_\Rssf(\gamma,\afrak_\Rssf)$
as a zero of
$N(X)$,
then
$(\tau_\Rssf(z,\afrak_\Rssf)
  -
  \tau_\Rssf(\gamma,\afrak_\Rssf)
 )^{n_\gamma}
$
has a zero
of order~$e$
at $z = \gamma$,
and thus
$z = \gamma$
is a zero
of order~$e/n_\gamma$
of
$\tau_\Rssf(z,\afrak)
 -
 \tau_\Rssf(\gamma,\afrak_\Rssf)
$.
On the other hand,
the multiplicity
$n_\gamma'$
of
$\tau_\Rssf(\gamma,\afrak_\Rssf)$
as a zero of
$T_{\mfrak,\kfrak_\Rssf}(X)$
is equal to
the number of
zeros of
$\tau_\Rssf(z,\afrak_\Rssf)
 -
 \tau_\Rssf(\gamma,\afrak_\Rssf)
$
which are incongruent
modulo~$\afrak$,
because
from
$\tau_\Rssf(\gamma',\afrak_\Rssf)
 -
 \tau_\Rssf(\gamma,\afrak_\Rssf)
 =
 0
$
it follows that
$\gamma'
 \equiv
 \varepsilon
 \,
 \gamma
 \mod
 \afrak
$,
where
$\varepsilon$
is a unit of $\Rssf$,
so
$\gamma'
 \equiv
 0
 \mod
 \afrak\,\mfrak^{-1}
$.
The zeros
$z
 \equiv
 \varepsilon
 \,
 \gamma
 \mod
 \afrak
$
of
$\tau_\Rssf(z,\afrak_\Rssf)
 -
 \tau_\Rssf(\gamma,\afrak_\Rssf)
$
are mapped to
one another
under the transformations
$z \to \varepsilon\,z + \alpha$
with
$\alpha$ from $\afrak_\Rssf$,
whence
they all have
the same order
$e/n_\gamma$;
thus
the order of
the elliptic function
$\tau_\Rssf(z,\afrak_\Rssf)
 -
 \tau_\Rssf(\gamma,\afrak_\Rssf)
$
which is~$e$
is equal to
$\dfrac{e}{n_\gamma}
 \cdot
 n_\gamma'
$,
i.e.~one has
$n_\gamma
 =
 n_\gamma'
$.

We thus have
\[
  \tau_\Rssf(z,\afrak_\Rssf\,\mfrak_\Rssf^{-1})
  \ =\ 
  \dfrac
  {Z(\tau_\Rssf(z,\afrak_\Rssf))}
  {T_{\mfrak,\kfrak_\Rssf}(\tau_\Rssf(z,\afrak_\Rssf))}
  .
\]

Now
if $N$ denotes
a natural number
divisible by~$\mfrak$,
then
\[
  T_{\mfrak,\kfrak_\Rssf}(X)
  \ =\ 
  \prod_
  {\substack{
      \gamma \bmod \afrak
      ,
      \\
      \gamma \equiv 0
      \bmod
      \afrak\,\mfrak^{-1}
      ,
      \\
      \gamma \not\equiv 0
      \bmod
      \afrak
  }}
  (X - \tau_\Rssf(\gamma,\afrak_\Rssf))
\]
is a factor of
\[
  L(X)
  \ =\ 
  \prod_
  {\substack{
      \gamma \bmod \afrak
      ,
      \\
      \gamma \equiv 0
      \bmod
      \afrak\,N^{-1}
      ,
      \\
      \gamma \not\equiv 0
      \bmod
      \afrak
  }}
  (X - \tau_\Rssf(\gamma,\afrak_\Rssf))
  ,
\]
but
$L(X)$
can be expressed via
the division polynomials:
\[
  L(X)
  \ =\ 
  \prod_
  {t \mid N
   \,,\,
   t > 1
  }
  T_t(X,j(\kfrak_\Rssf))
  ,
\]
so
its coefficients
lie in $\Omega_\Rssf$.
If
$L(X)
 =
 T_{\mfrak,\tfrak_\Rssf}(X)
 \,
 Q(X)
$
and
$Z(X)
 \,
 Q(X)
 =
 M(X)
$,
then
the coefficients of
\[
  T_{\mfrak,\tfrak_\Rssf}(X)
  \ =\ 
  L(X)
  \,\big/\,
  (\,
  M(X)
  \,,\,
  L(X)
  \,)
  \qquad
  \text{lie in $\Omega_\Rssf$}
  ,
\]
if that holds
for the coefficients of
$M(X)$.

Now to see this,
we substitute
for $\tau_\Rssf(z,\afrak_\Rssf)$
in $M(\tau_\Rssf(z,\afrak_\Rssf))$
the series expansion
\[
  \tau_\Rssf(z,\afrak_\Rssf)
  \ =\ 
  \sum_{\nu=-1}^{\infty}
  \eta_\nu
  \,
  u^\nu
  ,
  \qquad
  u
  \ =\ 
  z^e
  \,
  g^{(e)}(\afrak_\Rssf)^{-1}
  ;
\]
since the $\eta_\nu$
lie in $\Omega_\Rssf$,
it suffices
to show that
the coefficients of
the expansion of
$M(\tau_\Rssf(z,\afrak_\Rssf))$
in powers of $u$
lie in $\Omega_\Rssf$,
because
the coefficients of
the polynomial
$M(X)$
are then obtained as
the uniquely determined solutions of
a system of linear equations
in the field
$\Omega_\Rssf$.
But the expansion of
$M(\tau_\Rssf(z,\afrak_\Rssf))$
is obtained when,
in the right hand side of
the equation
\[
  M(\tau_\Rssf(z,\afrak_\Rssf))
  \ =\ 
  L(\tau_\Rssf(z,\afrak_\Rssf))
  \,
  \tau_\Rssf(z,\afrak_\Rssf\,\mfrak_\Rssf^{-1})
\]
we substitute
the expansions
\[
  \tau_\Rssf(z,\afrak_\Rssf)
  \ =\ 
  \sum_{\nu=-1}^{\infty}
  \eta_\nu
  \,
  u^\nu
  ,
  \qquad
  u
  \ =\ 
  z^e
  \,
  g^{(e)}(\afrak_\Rssf)^{-1}
\]
and
\[
  \tau_\Rssf(z,\afrak_\Rssf\,\mfrak_\Rssf)
  \ =\ 
  \sum_{\nu=-1}^{\infty}
  \eta_\nu'
  \,
  {u'}^\nu
  ,
  \qquad
  u'
  \ =\ 
  z^e
  \,
  g^{(e)}(\afrak_\Rssf\,\mfrak_\Rssf^{-1})^{-1}
\]
i.e.
\[
  \tau_\Rssf(z,\afrak_\Rssf\,\mfrak_\Rssf^{-1})
  \ =\ 
  \sum_{\nu=-1}^{\infty}
  \eta_\nu'
  \,
  \left(\ 
  \dfrac
  {g^{(e)}(\afrak_\Rssf)}
  {g^{(e)}(\afrak_\Rssf\,\mfrak_\Rssf^{-1})}
  \ \right)^\nu
  \,
  u^\nu
  .
\]
The coefficients of $L$
lie in $\Omega_\Rssf$,
and
the $\eta_\nu$
and
the $\eta_\nu'$
lie in $\Omega_\Rssf$.
But
$\dfrac
 {g^{(e)}(\afrak_\Rssf)}
 {g^{(e)}(\afrak_\Rssf\,\mfrak_\Rssf^{-1})}
$
also lie in $\Omega_\Rssf$,
because
$\alpha_1,\alpha_2$
is a basis of
$\afrak_\Rssf\,\mfrak_\Rssf^{-1}$,
so
a basis of $\afrak_\Rssf$
can be taken
in the form
$S
 \begin{pmatrix}
 \alpha_1
 \\
 \alpha_2
 \end{pmatrix}
 \,
 t
$
with a primitive matrix
$S$
and a natural number
$t$,
whence
\[
  \dfrac
  {g^{(e)}(\afrak_\Rssf)}
  {g^{(e)}(\afrak_\Rssf\,\mfrak_\Rssf^{-1})}
  \ =\ 
  t^e
  \,
  \dfrac
  {g^{(e)}
   \left(\ 
   S
   \begin{pmatrix}
     \alpha_1
     \\
     \alpha_2
   \end{pmatrix}
   \ \right)
  }
  {g^{(e)}
   \begin{pmatrix}
     \alpha_1
     \\
     \alpha_2
   \end{pmatrix}
  }
\]
and
$\dfrac
 {g^{(e)}
  \left(\ 
  S
  \begin{pmatrix}
    \alpha_1
    \\
    \alpha_2
  \end{pmatrix}
  \ \right)
 }
 {g^{(e)}
  \begin{pmatrix}
    \alpha_1
    \\
    \alpha_2
  \end{pmatrix}
 }
$
lies in $\Omega_\Rssf$
by~\S\textbf{III.12}.
With that
the proof
is finished.


\section*{E. Remarks}
\addcontentsline{toc}{section}{\hfill
E. Remarks
\hfill}

We have
concisely but completely
depicted
the theory of
the interrelations between
elliptic functions,
modular functions
and algebraic numbers,
usually called
``complex multiplication''
in short;
in the process
we have assumed
as known
the elementary parts of
the theory of
algebraic numbers.
It remains for us
to add some references
to the literature.

Our theory emerged
largely before
the completion of
general class field theory;
it was an important reason
for the latter's
development.
Therefore,
the classical results,
for which
we refer to
\textit{Weber}
\footnote{
\textit{Weber},
\ 
see Monographs.
},
seen from
today's point of view,
are incomplete
and
convoluted in reasoning.
Nevertheless,
our presentation of
the theory of singular moduli
(the ring class field),
as summarized
under the heading
``First Main Theorem'',
proceeds basically as in
\textit{Weber}
{\ \ \llap{${}^{15)}$}},
\textit{Fricke}
or
\textit{Fueter}
\footnote{
\textit{Fueter},
\ 
see Monographs.
}.
However,
the fundamental congruence~\eqref{eqn:33}
was first given by
\textit{Hasse}
{\ \ \llap{${}^{2)}$}}
(previously
only~\eqref{eqn:30}
was known)
and
the establishment of
the isomorphism
$\sigma
 :
 \Rfrak_\Rssf
 \to
 \Gfrak^{\Omega_\Rssf}_{\Sigma}
$
has been greatly simplified
here.

The theory of
the singular division values of
$\tau_\Rssf(z,\afrak_\Rssf)$
(second main theorem)
was developed from
an idea of \textit{Hasse}.
It consists of
introducing
the division values
$\tau_\Rssf
 \left(
  \dfrac{x_1\,\omega_1 + x_2\,\omega_2}{N}
  \,,\,
  \begin{pmatrix}
    \omega_1
    \\
    \omega_2
  \end{pmatrix}
 \right)
$
as higher level modular functions
and thus
largely avoiding
the transformation theory of
elliptic functions,
which only appears
in~\S\textbf{III.8}.
The establishment of
the fundamental congruence~\eqref{eqn:52}
also for
degree-two
prime divisors $\pfrak$
is new,
as is
the preparatory theorem
in~\S\textbf{III.15}.

\textit{Weber}
and
\textit{Fueter}
after him
have established
the theory of
ray class fields,
in particular
the congruence~\eqref{eqn:51}
(or the corresponding statement),
based more on
the transformation theory of
elliptic functions
\footnote{
\textit{Weber}
and
\textit{Fueter},
\ 
see Monographs.
};
however,
instead of
the function
$\tau_\Rssf$
which is
an elliptic function
of level~one
(invariant under $\Mfrak$),
they were
forced to use
functions of
higher levels
(\textit{Jacobi} functions),
which are easier
to handle.
That
it is possible
to manage
with
the function
$\tau_\Rssf$
of level~one
alone
was first shown by
\textit{Hasse}
\footnote{
\textit{Hasse}
loc. cit.${}^{2)}$
and ${}^{8)}$.
}.
However,
a method
was later given
by \textit{Deuring}
\footnote{
\textit{Deuring},
Abhandl. Akad. Mainz,
Nat.-Math. Kl.
1954.
}
to carry out
the approach of \textit{Weber}
also with
functions of level~one
such as
$\tau_\Rssf$.

But
from another point of view,
it is advantageous
to use
other elliptic functions
besides
$\tau_\Rssf$.
\textit{Fueter}
succeeds
in this way
to treat
the decomposition behavior
in the ray class field
$\Omega_\mfrak$
of the finitely many
prime divisors
of the field
$\Sigma$
which for us
must remain
excluded,
and
to compute
the discriminants of
the intermediate
class fields.
For the goal of
\textit{Fueter},
to prove
the \emph{completeness theorem}
asserting that
all abelian number fields
over $\Sigma$
are contained in
the ray class fields
$\Omega_\mfrak$,
that was unavoidable.
Today
one would simply take
this result
from Takagi's converse theorem of
general class field theory
(as \textit{Takagi} himself
has already done)
\footnote{
\textit{Takagi},
J. of the Coll. of Sci. Tokyo
\textbf{44},
no.~5 (1922).
},
because
the arguments
which
\textit{Fueter}
uses
for his proof of
the completeness theorem
are basically
the same ones
that appear
in the general considerations
of \textit{Takagi}.

It should also be noted that
it is actually
quite inappropriate
to study
a number field
by fixing
one (or some)
generating element
(as we have done
in the case of
$\Omega_\Rssf$
and
$\Omega_\mfrak$).
The exclusion of
finitely many
prime ideals
has its reason.
The purely algebraic theory,
which springs from
the theory of
algebraic function fields,
is largely free from
this defect,
and
will be reported
in this encyclopedia~\textup{I~2,~26}.

On specific questions
which
we did not address here:

\textit{Weber},
\textit{Fricke},
\textit{Fueter},
\textit{Watson}
treat the problem of
explicit computation of
class invariants
\ \ 
\footnote{
cf.~the monographs of
\textit{Weber},
\textit{Fricke}
and
\textit{Fueter},
furthermore
\textit{Watson},
J. London Math. Soc.
\textbf{6},
65 to~70
and
126 to~132
(1930);
J. f. Math.
\textbf{169},
238 to~251
(1933);
Proc. London Math. Soc. II.
\textbf{42},
398 to~409
(1936);
\textit{Mitra},
Indian phys.-math. J.
\textbf{2},
7 to~10
(1931);
Bull. Calcutta Math. Soc.
\textbf{24},
135 to~136
(1933).
}.

On the question of
whether
$\Omega_\mfrak
 =
 \Sigma(j(\kfrak_\Rssf),\tau_\Rssf(\kfrak_\mfrak))
$
is generated
over
$\Sigma$
by $\tau_\Rssf(\kfrak_\mfrak)$
alone,
see
\ \ 
\footnote{
\textit{Sugawara},
Proc. Phys.-Math. Soc. Japan III,
99 to~107
(1933);
J. f. Math.
\textbf{174},
189 to~191
(1936).
}.

On the number fields
which can be generated by
the singular values of
certain modular functions
of higher level,
see
\ \ 
\footnote{
\textit{S\"ohngen},
Math. Ann.
\textbf{111},
102 to~328
(1935).
}.

On the division values of
$\tau_\Rssf
 \left(
  z
  \,,\,
  \begin{pmatrix}
    \omega_1
    \\
    \omega_2
  \end{pmatrix}
 \right)
$
with respect to
non-maximal orders,
see
\ \ 
\footnote{
\textit{Franz},
J. f. Math.
\textbf{173},
60 to~64
(1935).
}

A new proof of
the first main theorem
which,
essentially,
gets by
without using
$\Delta(\omega)$,
is given by
\textit{Eichler}
\ \ 
\footnote{
\textit{Eicher},
Math. Zeitschr.
\textbf{64},
229 to~242
(1956).
}.


\addcontentsline{toc}{section}{}

\addcontentsline{toc}{section}{\hfill
Monographs
\hfill}

\addcontentsline{toc}{subsection}{
{\sl H.~Weber},~
Lehrbuch der Algebra. III. 3.~Aufl.,~
Braunschweig 1908.
}

\addcontentsline{toc}{subsection}{
{\sl R.~Fricke},~
Lehrbuch der Algebra. III.~
Braunschweig 1928.
}

\addcontentsline{toc}{subsection}{
{\sl R.~Fueter},~
Vorlesungen \"uber die singl\"aren Moduln
und die komplexe Multiplikation
der elliptischen Funktionen. I u.~II.~
Leipzig 1924 u.~1928.
}


\end{document}